\documentclass{amsproc}

\usepackage{amsmath,amssymb}
\newtheorem{theorem}{Theorem}[section]
\newtheorem{lemm}[theorem]{Lemma}
\newtheorem{prop}[theorem]{Proposition}

\theoremstyle{definition}
\newtheorem{defi}[theorem]{Definition}

\newtheorem{coro}[theorem]{Corollary}
\theoremstyle{remark}
\newtheorem{remark}[theorem]{Remark}

\newcommand{\ad}{\mathrm{ad}}
\allowdisplaybreaks[4]

\numberwithin{equation}{section}

\def\lg{\langle}
\def\rg{\rangle}

\def\al{\alpha}

\def\ep{\epsilon}

\def\p{\partial}

\def\dim{\hbox{dim}}
\def\ad{\hbox{ad}}

\def\a{\alpha}

\def\mod{\hbox{mod}}

\newfont{\df}{eufm10}
\def\vep{\varepsilon}
\def\ep{\epsilon}

\def\dim{\hbox{\rm dim}\,}

\def\ad{\hbox{\rm ad}\,}

\usepackage{txfonts}

\begin{document}
\allowdisplaybreaks
\title[Modular Quantization of Lie Algebra of Cartan Type $\mathbf{K}$]
{Modular Quantizations of Lie Algebras of Cartan Type $K$ \\ via Drinfeld Twists of Jordanian Type}

\author[Tong]{Zhaojia Tong}
\address{Department of Mathematics, Shanghai University,
Baoshan Campus, Shangda Road 99, Shanghai 200444, PR China; and
Shanghai Dianji University} \email{tongzhaojia@gmail.com}

\thanks{$^\star$N.H., corresponding author,
supported in part by the NNSF of China (No. 11271131).}

\author[Hu]{Naihong Hu$^\star$}
\address{Department of Mathematics, Shanghai Key Laboratory of Pure Mathematics and Mathematical Practice, East China Normal University,
Minhang Campus, Dong Chuan Road 500, Shanghai 200241, PR China}
\email{nhhu@math.ecnu.edu.cn}

\date{Version on Apr. 6nd, 2006}

\begin{abstract}
We construct explicit Drinfel'd twists of Jordanian type for the generalized Cartan
type $K$ Lie algebras in characteristic $0$ and obtain the corresponding quantizations, especially their integral forms. By making
modular reductions including modulo $p$ and modulo $p$-restrictedness
reduction, and base changes, we derive certain modular
quantizations of the restricted universal enveloping algebra
$\mathbf u(\mathbf{K}(2n{+}1;\underline{1}))$ for the restricted simple Lie algebra of Cartan type $K$ in characteristic $p$.
They are
new Hopf algebras of noncommutative and noncocommutative and with
dimension $p^{p^{2n+1}+1}$ (if $2n+4\not\equiv0 \; (\mod
p)$) or $p^{p^{2n+1}}$ (if $2n+4\equiv0 \; (\mod p))$ over a truncated $p$-polynomials ring, which also contain the well-known Radford algebras as Hopf subalgebras.
Some open questions are proposed.
\end{abstract}

\maketitle

This paper is a continuation of \cite{HW1, HW2,HTW} in which
quantizations of Cartan type  Lie algebras of types $W$, $S$ and $H$ were
studied. In the present paper, we continue to treat the same
questions both for the generalized Cartan type $K$ Lie algebras in
characteristic $0$ (for the definition, see \cite{OZ1}) and for the
restricted simple contact algebra $\mathbf{K}(2n{+}1;\underline{1})$ of Cartan type $K$  in the
modular case (for the definition, see \cite{H}, \cite{HR}). This work, together with papers \cite{HW1, HW2,HTW},
completes the work of this kinds of modular quantizations of Jordanian type for the Cartan
type Lie algebras over a field $\mathcal{K}$ with $\text{char}\,(\mathcal{K})
\geq 7$.

We survey some previous related work. After the works on quantum
groups which were introduced by Drinfel'd \cite{D} and Jimbo
\cite{J}, Drinfel¡¯d in \cite{Dr1} raised the question of the
existence of a universal quantization for Lie bialgebras.
Etingof-Kazhdan gave a positive answer to this question in
\cite{EK1, EK2}, where the Lie bialgebras they considered including
finite- and infinite-dimensional ones are those defined by
generalized Cartan matrices. Enriquez-Halbout showed that any
coboundary Lie bialgebra, in principle, can be quantized via a
certain Etingof-Kazhdan quantization functor \cite{EH}, and Geer
\cite{G} further extended Etingof-Kazhdan's work from Lie bialgebras
to the setting of Lie superbialgebras. After the work \cite{EK1,
EK2}, it is natural to consider the quantizations of Lie algebras of
Cartan type that are defined by differential operators. Grunspan
\cite{CG} obtained the quantization of the (infinite-dimensional)
Witt algebra $\bold W$ in characteristic $0$ using the twist due to
Giaquinto-Zhang \cite{AJ}, however, his approach didn't work for the
quantum version of its simple modular Witt algebra $\bold W(1;
\underline{1})$ in characteristic $p$. Hu-Wang in \cite{HW1}
obtained the quantizations of the generalized Witt algebra $\bold W$
in characteristic $0$ and the Jacobson-Witt algebra $\bold W(n;
\underline{1})$ in characteristic $p$; these are new families of
noncommutative and noncocommutative Hopf algebras of dimension
$p^{1+np^n}$ in characteristic $p>0$, while in rank $1$ case, the
work not only recovered Grunspan's one in characteristic $0$, but
also gave the correct modular quantum version. Note that the concept
of ``modular" quantization and relevant methods were settled in the
work \cite{HW1} and \cite{HW2} (also see \cite{HTW}).

Although, in principle, the possibility to quantize an arbitrary Lie
bialgebra has been proved (\cite{EK1, EK2, EH, ES, G}, etc.), an
explicit formulation of Hopf operations remains  nontrivial. In
particular, only a few kinds of twists were known with explicit
expressions, see \cite{AJ, KL, KLS, VO, NR} and the references
therein. In this research, we start with an explicit Drinfel'd twist
due to \cite{AJ,CG} and in fact, this Drinfel'd twist  is
essentially a variation (see the proof in \cite{HW2}) of the
Jordanian twist which first appeared, using different expression, in
Coll-Gerstenhaber-Giaquinto's paper \cite{CGG}, and recently used
extensively by Kulish et al (see \cite{KL, KLS}, etc.). Using this
explicit Drinfel'd twist we obtain \textit{vertical basic} twists
and \textit{horizontal basic} twists for the generalized Cartan type
$K$ Lie algebras and the corresponding quantizations in
characteristic $0$. These basic twists can afford many more
Drinfel'd twists, likewise on types $W,\,S,\,H$. To study the modular case,
what we discuss first involves the arithmetic property of
quantizations, for working out their quantization integral forms. To
this end, we have to work over the so-called {\it ``positive"} part
subalgebra $\mathbf{K}^+$ of the generalized Cartan type $K$ Lie
algebra. This is the crucial observation here. It is an
infinite-dimensional simple Lie algebra defined over a field of
characteristic $0$, while,  over a field of characteristic $p$, it
contains a maximal ideal $J_{\underline{1}}$ and the corresponding
quotient is exactly the algebra $\mathbf{K}'(2n{+}1;\underline{1})$. Its
derived subalgebra
$\mathbf{K}(2n{+}1;\underline{1})=\mathbf{K}'(2n{+}1;\underline{1})^{(1)}$
is a Cartan type restricted simple modular Lie algebra of $K$ type.
Secondly, in order to yield the {\it expected} finite-dimensional
quantizations of the restricted universal enveloping algebra of the
Contact  algebra $\mathbf{K}(2n{+}1;\underline{1})$, we need to
carry out the modular reduction process: {\it modulo $p$ reduction}
and {\it modulo ``$p$-restrictedness" reduction}, during which we
have to take the suitable {\it base changes}. These are the other
two crucial technical points. Our work gives a new class of
noncommutative and noncocommutative Hopf algebras of prime-power
dimension in characteristic $p$, which is significant to recognize
the Kaplansky's problem.

The paper is organized as follows. In Section 1, we recall  some
definitions and basic facts  related to  the  Cartan type $K$
Lie algebra and Drinfel'd twist. In Section 2, we construct the  Drinfel'd twists  for  the generalized Cartan type $K$ Lie
algebra, including
 \textit{vertical basic} twists and {\it horizontal basic} twists.
In Section 3, we quantize explicitly Lie bialgebra
structures of the generalized Cartan type $K$ Lie algebra by
the {\it vertical basic} Drinfel'd twists, and using the
similar methods as in type $H$, we obtain the quantizations of the restricted
universal enveloping algebra of the Contact algebra
$\mathbf{K}(2n{+}1;\underline{1})$. In Section 4, using the {\it
horizontal} twists, we get some new modular quantizations of horizontal type
of $\mathbf u(\mathbf{K}(2n{+}1;\underline{1}))$.  In Section 5, we give another two kinds of Drinfel'd twists and get some new modular quantizations of horizontal type
of $\mathbf u(\mathbf{K}(2n{+}1;\underline{1}))$. Finally, we present some open questions.

\section{Preliminaries}

\subsection{The generalized Cartan type $K$ Lie algebra and its subalgebra}
We recall the definitions of the generalized Cartan type $K$ Lie
algebras and the restricted simple Lie algebras of Cartan type from \cite{OZ1, HR} and some basics about their structures.

Let $\mathbb{F}$ be a field of characteristic $0$.
$\mathbb{Q}_{2n+1}=\mathbb{F}[x_{-n}^{\pm 1},{\cdots},x_{-1}^{\pm
1},x_{0}^{\pm 1},x_{1}^{\pm 1},{\cdots},x_{n}^{\pm 1}]$
($n\in\mathbb N$) be Laurent polynomial algebra, and $\partial_i$
coincide with the degree operator $x_i \frac{\partial}{\partial
x_i}$. Set
$T=\bigoplus\limits_{i=1}^{n}\mathbb{Z}\partial_{-i}\oplus
\mathbb{Z}\partial_i \bigoplus \mathbb{Z} \partial_{0}$, $x^{\al}=
x_{-n}^{\al_{-n}} {\cdots} x_{-1}^{\al_{-1}}x_0^{\al_0}x_1^{\al_1}
{\cdots} x_{n}^{\al_n}$, for
$\al=(\al_{-n},{\cdots},\al_{-1},\al_0,\al_1,{\cdots},\al_n)$.
Define ${\bf
W}=\mathrm{Der}(\mathbb{Q}_{2n+1})=\mathrm{Span}_{\,\mathbb{F}} \{
x^{\al} \partial \mid \al \in \mathbb{Z}^{2n+1},~\partial \in
T\,\}$. Then ${\bf W}$ is a Lie algebra of generalized Witt type
under the bracket
$$
[x^\al \partial,x^\beta \partial']=x^{\al+\beta} (\partial(\beta)\partial'-\partial'(\al)\partial),\qquad \forall \al,~ \beta \in \mathbb{Z}^{2n+1},~\partial,~\partial' \in T,
$$
where $\partial(\beta)=\langle\partial,\beta\rangle=\langle\beta,\partial\rangle=\sum\limits_{i=1}^{n} a_{-i}\beta_{-i}+a_i\beta_{i}+a_0\beta_0$, for $\partial=\sum\limits_{i=1}^{n}a_{-i}\partial_{-i}+a_i \partial_i+a_0\partial_0$ and $\beta=(\beta_{-n},{\cdots},\beta_{-1},\beta_0,\beta_1,{\cdots},\beta_{n})$. The bilinear map $\langle~,~\rangle:~T \times \mathbb{Z}^{2n+1} \longrightarrow ~\mathbb{Z}$ is nondegenerate in the sense that
\begin{gather*}
\partial (\alpha)= \langle \partial,~ \alpha \rangle =0  \quad (\forall
~\partial \in T), \ \Longrightarrow \ \alpha=\underline{0},\\
 \partial (\alpha)=\langle \partial,~\alpha \rangle =0 \quad
(\forall~\alpha \in \mathbb{Z}^{2n+1}), \ \Longrightarrow \ \partial=0.
\end{gather*}
where $\underline{0}=(0,\cdots,0)$.
${\bf W}$ is an infinite dimensional Lie algebra over $\mathbb{F}$.

Consider the linear map $\mathcal D_K
:~\mathbb{Q}_{2n{+}1}~\longrightarrow {\bf W}$ defined by
\begin{equation*}
\begin{split}
\mathcal D_K (x^{\al})&=\Big(2{-}\sum\limits_{i=1}^{n}(\al_i{+}\al_{-i})\Big)x^\al \frac{\p}{\p x_0}\\
&\ +\sum\limits_{i=1}^{n} \Big((\al_0x^{\al+\ep_i-\ep_0}+\al_{-i}x^{\al-\ep_{-i}})\frac{\p}{\p x_i}+(\al_0 x^{\al+\ep_{-i}-\ep_0}-\al_i x^{\al-\ep_i})\frac{\p}{\p x_{-i}}\Big).
\end{split}
\end{equation*}
So we get
$$
\mathcal D_K (x^{\al})=\big(2-\sum\limits_{i=1}^{n}(\al_i{+}\al_{-i})\big)x^{\al-\ep_0}\p_0
+\sum\limits_{i=1}^{n}\al_0x^{\al-\ep_0}(\p_i{+}\p_{-i})+\sum\limits_{i=1}^{n}x^{\al-\ep_i-\ep_{-i}} (\al_{-i}\p_i{-}\al_i \p_{-i}).
$$
It is easy to see that $\mathcal D_K $ is injective. Then for any $x^{\al},x^{\beta} \in \mathbb{Q}_{2n+1}$, the Lie bracket becomes
\begin{equation*}
\begin{split}
[\mathcal D_K& (x^\al),\mathcal D_K (x^\beta)]
=\mathcal D_K \Big[\Big(2x^\al{-}\sum\limits_{i=1}^{n}(\al_i+\al_{-i})x^\al\Big)\beta_0x^{\beta-\ep_0}-\Big(2x^\beta{-}\sum\limits_{i=1}^{n}(\beta_i+\beta_{-i})x^\beta\Big)\al_0 x^{\al-\ep_0}\\
&\quad+\sum\limits_{i=1}^{n}(\al_{-i} \beta_{i}-\al_i\beta_{-i})x^{\al+\beta-\ep_i-\ep_{-i}}\Big]
\end{split}
\end{equation*}
\begin{equation*}
\begin{split}
&=\mathcal D_K \bigg(\Big(\big(2{-}\sum\limits_{i=1}^{n}(\al_i{+}\al_{-i})\big)\beta_0{-}\big(2{-}\sum\limits_{i=1}^{n}
(\beta_i{+}\beta_{-i})\big)\al_0\Big)x^{\al+\beta-\ep_0}{+}\sum\limits_{i=1}^{n}(\al_{-i}\beta_i{-}\al_i\beta_{-i})x^{\al+\beta-\ep_i-\ep_{-i}}\bigg).
\end{split}
\end{equation*}

It follows that ${\bf K}=\mathcal D_K (\mathbb{F}[x_{-n}^{\pm 1},{\cdots}, x_{-1}^{\pm 1}, x_{0}^{\pm1}, x_1^{\pm 1},
{\cdots}, x_n^{\pm 1}])$, and ${\bf K}^+=\mathcal D_K (\mathbb{F}[x_{-n},{\cdots},
x_{-1}$, $x_0, x_1, {\cdots}, x_n])$ are Lie subalgebras of ${\bf W}$.
${\bf K}$, ${\bf K}^+$ are simple algebras with basis $\mathcal D_K (x^\al),~\al$ $\in \mathbb{Z}^{2n+1}$, $\mathcal D_K (x^\al),~\al \in \mathbb{Z}_+^{2n+1}$, respectively.
${\bf K}$ is the Lie algebra of generalized Cartan type $K$ (\cite{OZ1}).

\subsection{The Contact algebra ${\bf K}(2n{+}1;\underline{1})$}
Assume that char$(\mathcal{K})=p$, then by definition (see
\cite{H}), the Jacobson-Witt algebra $\mathbf{W}(2n{+}1;\underline{1})$
is a restricted simple Lie algebra over a field $\mathcal{K}$. Its
structure of $p$-Lie algebra is given by $D^{[p]}=D^p,\; \forall\, D
\in \mathbf{W}(2n{+}1;\underline{1})$ with a basis $\{\,x^{(\alpha)}D_j
\mid -n \leq j\leq n,\, \ 0 \leq \alpha \leq \tau \}$,
where $\tau=(p-1,\cdots,p-1) \in \mathbb{N}^{2n+1}$;
$\ep_i=(\delta_{-n,i},{\cdots},\delta_{-1,i},\delta_{0,i},\delta_{1,i},{\cdots},\delta_{n,i})$ with $x^{(\ep_i)}=x_i$, $x^{(\al)}\in
\mathcal {O}(2n{+}1;\underline{1})=\mathrm{Span}_\mathcal{K}\{x^{(\alpha)} \mid 0 \leq \alpha
\leq \tau \}$, the restricted divided power algebra with
$x^{(\alpha)} x^{(\beta)}=\binom{\alpha+\beta}{\alpha} x^{\alpha
+\beta}$ and a convention: $x^{(\alpha)}=0$ if $\alpha$ has a
component $\alpha_j < 0$ or $\geq p$,
where
 $\binom{\alpha +\beta }{\alpha}=
 \prod\limits_{i=1}^{n}\binom{\alpha_i + \beta_i}{\alpha_i}\binom{\alpha_{-i} + \beta_{-i}}{\alpha_{-i}}
 \binom{\alpha_{0}+\beta_{0}}{\alpha_{0}}$.
Define $\mathcal D_K :~\mathcal{O}(2n{+}1;\underline{1})\longrightarrow {\bf W}(2n{+}1;\underline{1})$, where $\mathcal D_K (x^{(\al)})=(2-\sum\limits_{i=1}^{n}(\al_i+\al_{-i}))x^{(\al)}D_0+
\sum\limits_{i=1}^{n}x^{(\al-\ep_0)}(x^{(\ep_i)}D_i+x^{(\ep_{-i})}D_{-i})+\sum\limits_{i=1}^{n}
(x^{(\alpha-\epsilon_{-i})}D_i-x^{(\alpha-\epsilon_{i})}D_{-i})$.
Then the subspace ${\bf K}'(2n{+}1;~\underline{1})=\mathcal D_K (\mathcal{O}(2n{+}1;\underline{1}))$ is a Lie subalgebra of ${\bf
W}(2n{+}1;\underline{1})$. Its derived algebra ${\bf K}(2n{+}1;\underline{1})$ is called the Contact algebra.
Define $|| \al ||=|\al|+\al_0-2$, Then ${\bf K}(2n{+}1;\underline{1})=\bigoplus\limits_{i=-2}^{s}{\bf K}(2n{+}1;\underline{1})_i$ is graded,
where ${\bf K}(2n{+}1;\underline{1})_i=\text{Span}_{\mathcal{K}}\{\,\mathcal D_K (x^{(\al)})\mid 0 \leq \al \leq \tau,~||\al||=i\,\} ,~  s=(2n+2)(p-1)-2$ when $2n{+}4 \nequiv 0 \; (\mod\; p);~{\bf K}(2n{+}1;\underline{1})_i=\text{Span}_{\mathcal{K}}\{\,\mathcal D_K (x^{(\al)})\mid 0 \leq \al < \tau,~||\al||=i\,\}, s=(2n{+}2)(p{-}1){-}3$ when $2n{+}4 \equiv 0 \;(\mod\; p)$. Then by Theorem 5.5 of \cite{H}, we have
\begin{equation*}
\begin{split}
{\bf K}(2n{+}1;\underline{1})=\begin{cases}\text{Span}_{\mathcal{K}}\{\,\mathcal D_K (x^{(\al)})\mid x^{(\al)} \in \mathcal{O}(2n{+}1;\underline{1}),~0 \leq \al \leq \tau\,\},& \mathrm{if} \ 2n{+}4 \nequiv 0 \ (\mod\; p),\\
\text{Span}_{\mathcal{K}}\{\,\mathcal D_K (x^{(\al)})\mid x^{(\al)} \in \mathcal{O}(2n{+}1;\underline{1}),~0 \leq \al < \tau\,\},& \mathrm{if} \ 2n{+}4 \equiv 0 \ (\mod\; p),\\
\end{cases}
\end{split}
\end{equation*}
is a Lie $p$-subalgebra of ${\bf W}(2n{+}1;\underline{1})$ with restricted gradation.

By definition (cf. \cite{HR}), the restricted universal enveloping
algebra ${\bf u(K}(2n{+}1;\underline{1}))$ is isomorphic to $U({\bf
K}(2n{+}1;\underline{1}))/I$, where $I$ is the Hopf ideal of $U({\bf
K}(2n{+}1;\underline{1}))$ generated by $(\mathcal D_K (x^{(\ep_k +\ep_{-k})}))^p
- \mathcal D_K ( x^{(\ep_k+\ep_{-k})}),~(\mathcal D_K (x^{(\ep_0)}))^p-\mathcal D_K (x^{(\ep_0)}),~(\mathcal D_K ( x^{(\alpha)}))^p$ with $\alpha \neq \ep_k +\ep_{-k},~\ep_0$, $1 \leq k \leq n$.
Since
\begin{equation*}
\begin{split}
\dim _{\mathcal{K}}(({\bf K}(2n{+}1;\underline{1})))=\begin{cases}p^{2n+1} ,&~\mathrm{if}~~2n+4 \nequiv 0 \;(\mod \;p),\\
p^{2n+1}-1 ,&~\mathrm{if}~~2n+4 \equiv 0 \;(\mod\; p),\\
\end{cases}
\end{split}
\end{equation*}
 we have
\begin{equation*}
\begin{split}
\dim _{\mathcal{K}}({\bf u}({\bf K}(2n{+}1;\underline{1})))=\begin{cases}p^{p^{2n+1}} ,&~\mathrm{if}~~2n+4 \nequiv 0 \;(\mod\; p),\\
p^{p^{2n+1}-1} ,&~\mathrm{if}~~2n+4 \equiv 0 \;(\mod\; p).\\
\end{cases}
\end{split}
\end{equation*}

\subsection{Quantization by Drinfel'd twists} The following result is well-known (\cite{D}).

\begin{lemm}{ \label{twist1}
Let $(A, m, \iota, \Delta_0, \vep, S_0)$ be a Hopf algebra over a
commutative ring. A Drinfel'd twist $\mathcal {F}$ on $A$ is an
invertible element of $A \otimes A$ such that
\begin{gather*}
(\mathcal{F} \otimes 1)(\Delta_0 \otimes \text{\rm
Id})(\mathcal{F})=(1
\otimes \mathcal{F})(\text{\rm Id} \otimes \Delta_0)(\mathcal{F}),\\
(\varepsilon \otimes \text{\rm Id})(\mathcal{F})=1=(\text{\rm Id} \otimes
\varepsilon)(\mathcal{F}).
\end{gather*}

Then, $w=m(\text{\rm Id} \otimes S_0)(\mathcal{F})$ is invertible
in $A$ with $w^{-1}=m(S_0 \otimes \text{\rm Id})(\mathcal{F}^{-1})$.

Moreover, if we define $\Delta: A \longrightarrow A \otimes
A$ and $S: A \longrightarrow A$ by
$$\Delta (a)=\mathcal{F} \Delta_0 \mathcal{F}^{-1},~~~~S(a)=w S_0 (a)
w^{-1},$$
then $(A, m, \iota, \Delta, \vep,S)$ is a new Hopf algebra, called
the twisting of A by Drinfel'd twist $\mathcal{F}$.}
\end{lemm}

Let $\mathbb{F}[[t]]$ be a ring of formal power series over a
field $\mathbb{F}$ with $\textit{char}\,(\mathbb{F})=0$. Assume that $L$ is a
triangular Lie algebra over $\mathbb{F}$ with a classical
Yang-Baxter $r$-matrix $r$. (see
\cite{D},~\cite{EK2}). Let $U(L)$ denote the universal enveloping
algebra of $L$, with the standard Hopf algebra
 $(U(L),~m,~\iota,~\Delta_0,~\ep_0,~S_0)$.

Let us consider the topologically free
 $\mathbb{F}[[t]]$-algebra $U(L)[[t]]$ (for definition, see
\cite{EK2}, p.4), which can be viewed as an associative
 $\mathbb{F}$-algebra of formal power series with coefficients
in $U(L)$. Naturally, $U(L)[[t]]$ equips with an induced Hopf
algebra structure arising from that on $U(L)$. By abuse of
notation, we denote it
by $(U(L)[[t]],~m,~\iota,~\Delta_0,~\ep_0,~S_0)$.

\begin{defi}{\upshape \cite{HW1}
For a triangular Lie
bialgebra $L$ over $\mathbb{F}$ with char($\mathbb{F}$)$=0$,
$U(L)[[t]]$ is called a \textit{quantization} of $U(L)$ by a
Drinfel'd
twist $\mathcal{F}$ over $U(L)[[t]]$ if
$U(L)[[t]]/tU(L)[[t]]$ $\cong U(L)$, and $\mathcal{F}$ is determined
by its $r$-matrix $r$ (namely, its Lie bialgebra structure).}
\end{defi}

\section{Drinfel'd twist in U({\bf K})[[t]]}
\subsection{Construction of Drinfel'd twists}
Let $L$ be a Lie algebra containing linearly independent elements
$h$ and $e$ satisfying $[h,e]=e$; then the classical Yang-Baxter
$r$-matrix $r=h \otimes e -e \otimes h$ equips $L$ with the
structure of a triangular coboundary Lie bialgebra (see \cite{M}).
To describe a quantization of $U(L)$ by a Drinfel'd twist
$\mathcal{F}$ over $U(L)[[t]]$, we need an explicit construction for
such a Drinfel'd twist. In what follows, we shall see that such a
Drinfel'd twist depends on the choice of two distinguished elements
$h$  and $e$ arising from its $r$-matrix $r$.

For any element of a unital $R$-algebra ($R$ a ring) and $a \in R$,
we set
\begin{gather*}
x_{a}^{\langle m \rangle}:=(x+a)(x+a+1)\cdots(x+a+m-1),\\
x_{a}^{[m]}:=(x+a)(x+a-1)\cdots (x+a-m+1),
\end{gather*}
and then denote $x^{\langle m \rangle}:=x_{0}^{\langle m \rangle},
x^{[m]}:=x_{0}^{[m]}$.

Note that $h$ and $e$ satisfy the following equalities
\begin{eqnarray}
e^s\cdot h_a^{[m]}=h_{a-s}^{[m]}\cdot e^s, \\
e^s\cdot h_a^{\langle m \rangle }=h_{a-s}^{\langle m \rangle }\cdot
e^s,
\end{eqnarray} where $m,\,s$ are non-negative
integers, $a \in \mathbb{F}$.

For $a \in \mathbb{F}$, we set
$\mathcal{F}_a=\sum\limits_{r=0}^{\infty}\frac{(-1)^r}{r!}h_a^{[r]}\otimes
e^rt^r, F_a=\sum\limits_{r=0}^{\infty}\frac{1}{r!}h_a^{\lg
r\rg}\otimes e^rt^r,$ $ u_a=m\cdot(S_0\otimes \text{\rm Id})(F_a),
v_a=m\cdot(\text{\rm Id}\otimes S_0)(\mathcal{F}_a).$ Write
$\mathcal{F}=\mathcal{F}_0,\, F=F_0,\,u=u_0,\,v=v_0$. Since
$S_0(h_a^{\lg r\rg})=(-1)^rh_{-a}^{[r]}$ and $S_0(e^r)=(-1)^re^r$,
one has $ v_a=\sum\limits_{r=0}^{\infty}\frac{1}{r!}h_a^{[r]}
e^rt^r, \ u_b=\sum\limits_{r=0}^{\infty}\frac{(-1)^r}{r!}h_{-b}^{[r]} e^rt^r$.

\begin{lemm} \label{1.6}$($\cite{CG}$)$
For $a,\, b \in \mathbb{F}$, one has
$$
\mathcal{F}_a F_b=1\otimes(1-et)^{a-b}, \quad\text{and }\quad v_a
u_b=(1-et)^{-(a+b)}.
$$
\end{lemm}

\begin{coro}\label{1.7}{
For $a \in \mathbb{F}$, $\mathcal{F}_a$ and $u_a$ are invertible
with $\mathcal{F}_a^{-1}=F_a$ and $u_a^{-1}=v_{-a}$.  In particular,
$\mathcal{F}^{-1}=F$ and $u^{-1}=v$.}
\end{coro}

\begin{lemm}\label{1.8} {
$($\cite{HW1}$)$ For any positive integers $r$, we have
$$\Delta_0(h^{[r]})=\sum\limits_{i=0}^r \dbinom{r}{i}h^{[i]}\otimes
h^{[r-i]}.$$  Furthermore, $\Delta_0(h^{[r]})=\sum\limits_{i=0}^r
\dbinom{r}{i}h^{[i]}_{-s}\otimes h^{[r-i]}_s$ for any $s \in
\mathbb{F}$.}
\end{lemm}

\begin{prop} \label{twist}$($\cite{CG, HW1}$)$ {If a Lie algebra $L$ contains a two-dimensional solvable
Lie subalgebra with a basis $\{h, e\}$ satisfying $[h, e]=e$, then
$\mathcal{F}=\sum\limits_{r=0}^{\infty}\frac{(-1)^r}{r!}h^{[r]}$
$\otimes\, e^rt^r$ is a Drinfel'd twist on $U(L)[[t]]$.}
\end{prop}

\begin{remark}(\cite{HW2})
 Kulish et al used early the so-called
{\it Jordanian twist} (see \cite{KL}) with the two-dimensional
carrier subalgebra $B(2)$ such that $[H,E]=E$, defined by the
canonical twisting element
$$\mathcal{F}_{\mathcal{J}}^{c}=\mathrm{exp}(H\otimes\sigma(t)), \quad
\sigma(t)=\mathrm{ln}(1+Et),$$ where
$\exp(X)=\sum_{i=0}^{\infty}\frac{X^n}{n!}$ and
$\ln(1+X)=\sum_{n=1}^{\infty}\frac{(-1)^{n+1}}{n}X^n$.

By the proof of \cite{HW2}, we can rewrite the twist $\mathcal{F}$
in Proposition \ref{twist} as
$$\mathcal{F}=\sum_{r=0}^{\infty}
\frac{(-1)^r}{r!}H^{[r]}\otimes E^rt^r=\exp(H\otimes\sigma'(t)),
\quad \sigma'(t)=\ln(1-Et),$$ where $[H,-E]=-E$. So there is no
difference between the twists $\mathcal{F}$ and
$\mathcal{F}_{\mathcal{J}}^{c}$. They are essentially the same up to
an isomorphism on the carrier subalgebra $B(2)$.
\end{remark}
\subsection{Basic Drinfel'd twists}
Take two distinguished elements $h=\mathcal D_K  (x^{\ep_k+\ep_{-k}}),~e=\mathcal D_K (x^{\al})$, such that $[h,e]=e$, where $1 \leq k \leq n$. It is easy to see that $\al_k-\al_{-k}=1$. Using the result of
\cite{M}, we have the following:

\begin{prop}\label{rmatrix}
There is a triangular Lie bialgebra structure on $\mathbf{K}$,
given by the classical Yang-Baxter $r$-matrix
$$
r:=\mathcal D_K (x^{{\ep_k+\ep_{-k}}}) \otimes \mathcal D_K ( x^{\alpha})-\mathcal D_K ( x^{\alpha})
\otimes \mathcal D_K ( x^{{\ep_k+\ep_{-k}}}),\quad 1 \leq k
\leq n,
$$
where $\alpha \in \mathbb{Z}^{2n+1},~\alpha_k-\alpha_{-k}
=1$, and $[\mathcal D_K (x^{\ep_k+\ep_{-k}}),\mathcal D_K (x^{\alpha})]=\mathcal D_K (x^{\alpha})$.\hfill\qed
\end{prop}

Fix two distinguished elements $h= \mathcal D_K ( x^{{\ep_k+\ep_{-k}}}),~e=\mathcal D_K (x^{\alpha})$, with $\alpha_k-\alpha_{-k}=1$, then
$\mathcal{F}=\sum\limits_{r=0}^{\infty}\frac{(-1)^r}{r!}h^{[r]}\otimes
e^rt^r$ is a Drinfel'd twist on $U(\mathbf{K})[[t]]$. But the
coefficients of the quantizations of standard Hopf structure
$(U(\mathbf{K})[[t]], m, \iota, \Delta_0$, $S_0, \varepsilon_0)$ by
$\mathcal{F}$ may be not integral. In order to get integral forms,
it suffices to consider what conditions are needed for those
coefficients to be integers.

\begin{lemm}\label{relation3} $($\cite{CG}$)$
For any $a,\,k,\,\ell\in\mathbb{Z}$,
$a^\ell\prod\limits_{j=0}^{\ell-1}(k{+}ja)/\ell!$ is an
integer.\hfill\qed
\end{lemm}

From this Lemma, we are interested in the following two simple cases:

\smallskip
$\mathrm{(i)}\ h=\mathcal D_K (x^{\ep_k+\ep_{-k}}),~e=\mathcal D_K (x^{2
\ep_{k}+\ep_{-k}}) \quad (1 \leq k \leq n);$

\smallskip
$\mathrm{(ii)}\ h=\mathcal D_K (x^{\ep_k + \ep_{-k}}),~e=\mathcal D_K (x^{\ep_k +
\ep_m}) \quad (1 \leq k \neq |m| \leq n).$

Let $\mathcal{F}(k)$ be the corresponding Drinfel'd thist in case (i)
and $\mathcal{F}(k;m)$ the corresponding Drinfel'd twist in case (ii).

\begin{defi}\label{VH twist}$\mathcal{F}(k)$ $(1\leq k \leq
n)$ are called {\it  vertical basic Drinfel'd twists};  $\mathcal{F}(k;m)$ $(1 \leq k \neq |m| \leq n)$
 are called  {\it horizontal basic Drinfel'd twists}.
\end{defi}

\begin{remark}\label{2.3}
In case (i): we get $n$ {\it vertical basic Drinfel'd twists}
$\mathcal{F}(1)$, $\cdots,\mathcal{F}(n)$, over
 $U(\mathbf{K}^+_{\mathbb{Z}})[[t]]$. It is interesting to consider the products of some {\it basic
Drinfel'd twists},  one can get many more new Drinfel'd twists which
will lead to many more new complicated quantizations not only over
the $U(\mathbf{K}_{\mathbb{Z}}^+)[[t]]$, but over the
$\mathbf{u}(\mathbf{K}(2n{+}1;\underline{1}))$ as well, via our modulo
reduction approach developed in the next section.

In case (ii): according to the parametrization of twists $\mathcal
F(k;m)$, we obtain $2n(n{-}1)$ {\it horizontal basic Drinfel'd twists}
over $U(\mathbf{K}^+_{\mathbb{Z}})[[t]]$. We will discuss these
twists and corresponding quantizations in Section 4.
\end{remark}

\subsection{More Drinfel'd twists}
We consider the products of pairwise different and mutually
commutative basic Drinfel'd twists and can get many more new
complicated quantizations not only over the
$U(\mathbf{K}_{\mathbb{Z}}^+)[[t]]$, but over the
$\mathbf{u}(\mathbf{K}(2n+1;\underline{1}))$ as well. Note that
$[\mathcal{F}(k), \mathcal{F}(k^{\prime})]=0$ for $1 \leq k \neq
k^{\prime} \leq n$. This fact, according to the definition of
$\mathcal{F}(k)$, implies the commutative relations in the case
when $1 \leq k \neq k^{\prime} \leq n$:
\begin{equation}\label{relation0}
\begin{split}
(\mathcal{F}(k)\otimes 1)(\Delta_0\otimes\text{\rm Id})
(\mathcal{F}(k^{\prime}))&=(\Delta_0\otimes\text{\rm Id})
(\mathcal{F}(k))(\mathcal{F}(k^{\prime})\otimes 1),\\
(1\otimes \mathcal{F}(k))(\text{\rm Id}\otimes\Delta_0)
(\mathcal{F}(k^{\prime}))&=(\text{\rm Id}\otimes\Delta_0)
(\mathcal{F}(k))(1\otimes\mathcal{F}(k^{\prime})),
\end{split}
\end{equation}
which give rise to the following property.
\begin{theorem}\label{twist3}
$\mathcal{F}(k)\mathcal{F}(k^\prime) \ (1 \leq k \neq k' \leq n )$ is
still a Drinfel'd twist on $U(\mathbf{K}_{\mathbb{Z}}^+)[[t]]$.
\end{theorem}
\begin{proof} Similar to the proof of Theorem 2.9 in \cite{HTW}.
\end{proof}
%

More generally, we have the following
\begin{coro}\label{2.5}
Let $\mathcal{F}(k_1),\cdots, \mathcal{F}(k_m)$ be $m$ pairwise
different basic Drinfel'd twists and
 $[\mathcal{F}(k_i), \mathcal{F}(k_s)]=0$ for all $1\leq i\neq s\leq n$.
 Then $\mathcal{F}(k_1) \cdots \mathcal{F}(k_m)$ is still a Drinfel'd twist.
\end{coro}
We denote $\mathcal{F}_m=\mathcal{F}(k_1)\cdots \mathcal{F}(k_m)$
and its length as $m$. We shall show that the twisted structures
given by Drinfel'd twists with different product-length are
nonisomorphic.

\begin{defi}\label{com}(\cite{MG}, \cite{HW2}) {\upshape
A Drinfel'd twist $\mathcal{F} \in A\otimes A$ on any Hopf algebra
$A$ is called \textit{compatible} if $\mathcal{F}$ commutes with the
coproduct $\Delta_0$.}
\end{defi}

In other words, twisting a Hopf algebra $A$ with a
\textit{compatible} twist $\mathcal{F}$ gives exactly the same Hopf
structure, that is, $\Delta_{\mathcal{F}}=\Delta_0$. The set of
\textit{compatible} twists on $A$ thus forms a group.

\begin{lemm}\label{2.7}$($\cite{MG}$)$
Let $\mathcal{F} \in A\otimes A$ be a  Drinfel'd twist on a Hopf
algebra $A$. Then the twisted structure induced by $\mathcal{F}$
coincides with the structure on $A$  if and only if $\mathcal{F}$ is
a compatible twist.
\end{lemm}

Using the same proof as in Theorem \ref{twist3}, we obtain

\begin{lemm}\label{twist4}$($\cite{HW2}$)$
Let $\mathcal{F}, \mathcal{G} \in A\otimes A$ be Drinfel'd twists on
a Hopf algebra $A$ with
$\mathcal{F}\mathcal{G}=\mathcal{G}\mathcal{F}$ and $\mathcal{F}\neq
\mathcal{G}$. Then $\mathcal{F}\mathcal{G}$ is a Drinfel'd twist.
 Furthermore, $\mathcal{G}$ is a Drinfel'd twist on
$A_{\mathcal{F}}$, $\mathcal{F}$ is a Drinfel'd twist on
$A_{\mathcal{G}}$ and $\Delta_{\mathcal{F}\mathcal{G}}=(\Delta_{\mathcal{F}})_{\mathcal{G}}
=(\Delta_{\mathcal{G}})_{\mathcal{F}}$.
\end{lemm}

\begin{prop}\label{2.9}
Drinfel'd twists
$\mathcal{F}^{\zeta(i)}=\mathcal{F}(1)^{\zeta_1}{\cdots}\mathcal{F}(n)^{\zeta_n}$
$($where $\zeta(i)=(\zeta_1,{\cdots},\zeta_{n})=(\underbrace{1,{\cdots},1}_i,0,{\cdots},0)\in
\mathbb Z_2^{n})$ lead to $n$ different twisted Hopf algebra
structures on $U(\mathbf{K}_{\mathbb{Z}}^{+})[[t]]$.
\end{prop}


\section{Quantizations of vertical type for Lie bialgebra of Cartan type $K$}
In this subsection, we explicitly quantize the Lie bialgebra ${\bf K}$ by the
 vertical basic Drinfel'd twists, and obtain certain quantizations of the restricted
 universal enveloping algebra ${\bf u}({\bf K}(2n{+}1;\underline{1}))$ by the modular reduction and base changes.
\subsection{Quantization integral forms of the $\mathbb{Z}$-form ${\bf K}_{\mathbb{Z}}^{+}$ in characteristic $0$}

For the universal enveloping algebra $U(\mathbf{K})$ for the Lie
algebra $\mathbf{K}$ over $\mathbb{F}$, denote by $(U(\mathbf{K}), m,\iota, \Delta_0, S_0$, $\varepsilon_0)$ the standard
Hopf algebra structure. We can perform the process of twisting the standard Hopf structure by the vertical Drinfel'd twist $\mathcal{F}$. We need to give some commutative relations, which are important to the quantizations  of Lie bialgebra structure of $\mathbf{K}$ in the sequel.
\begin{lemm}
For $h=\mathcal D_K (x^{\ep_k+\ep_{-k}})$ in ${\bf K}$, $a \in \mathbb{F},~m$~a non-negative integer, the following equalities hold in $U({\bf K})$:
\begin{gather}
\mathcal D_K ( x^{\alpha})\cdot
h_a^{[m]}=h_{a+(\alpha_{-k}-\alpha_k)}^{[m]}\cdot \mathcal D_K (
x^\alpha),\label{vrelation1}\\
\mathcal D_K ( x^\alpha) \cdot h_{a}^{\langle m \rangle}
=h_{a+(\alpha_{-k}-\alpha_k)}^{\langle m \rangle}\cdot \mathcal D_K ( x^\alpha).\label{vrelation2}
\end{gather}
\end{lemm}

To simplify the formulas, let us introduce the operator $d^{(\ell)}$~on~$U({\bf K})$ defined by $d^{(\ell)}=\frac{1}{\ell !}(\ad e)^\ell$. So we can get
\begin{lemm}
For $\mathcal D_K (x^\al) \in U({\bf K})$, the following equalities hold
\begin{gather}
d^{(\ell)}\Big(\mathcal D_K (x^\al)\Big)=\sum\limits_{j=0}^{\ell} A_j B_{\ell-j}\mathcal D_K (x^{\al+\ell(2\ep_k+\ep_{-k})-j(\ep_k+\ep_{-k})-(\ell-j)\ep_0}),\label{vrelation3}
\\
d^{(\ell)}(a_1 \cdots
a_s)=\sum\limits_{\ell_1+\cdots+\ell_s=l}d^{(\ell_1)}(a_1)\cdots
d^{(\ell_s)}(a_s),\label{vrelation4}\\
\mathcal D_K (x^\al) \cdot e^m=\sum\limits_{\ell=0}^{m}(-1)^\ell \binom{m}{\ell}\,\ell !\,e^{m-\ell}\cdot\sum\limits_{j=0}^{\ell} A_j B_{\ell-j}\mathcal D_K (x^{\al+\ell(2\ep_k+\ep_{-k})-j(\ep_k+\ep_{-k})-(\ell-j)\ep_0}),\label{vrelation5}\\
(\ad \mathcal D_K (x^\al))^\ell \big( e \big)=\sum\limits_{j=0}^{\ell}\binom{\ell}{j}C_{\ell-j}D(j,k) \mathcal D_K  (x^{\ell\al+2\ep_k+\ep_{-k}-j(\ep_k+\ep_{-k})-(\ell-j)\ep_0}),\label{vrelation6}
\end{gather}
where $A_j=\frac{1}{j!}\prod\limits_{i=0}^{j-1}(\al_k-2\al_{-k}+i),
~B_{\ell-j}=\frac{1}{(\ell-j)!}\prod\limits_{i=0}^{\ell-j-1}(i-\al_0),
~C_{\ell-j}=\prod\limits_{i=0}^{\ell-j-1}(i||\al||+\al_0),~D(j,k)=\prod\limits_{i=0}^{j-1}\big((2{-}i)\al_{-k}{-}(1{-}i)\al_k\big)$, $A_0=B_0=C_0=D(0,k)=1$.
\end{lemm}

\begin{proof}
For (\ref{vrelation3}), use induction on $\ell$.
When $\ell=1$,
\begin{equation*}
\begin{split}
d\Big(\mathcal D_K (x^\al)\Big)=[\mathcal D_K (x^{2\ep_k+\ep_{-k}}),\mathcal D_K (x^\al)]=(-\al_0)\mathcal D_K (x^{\al+2\ep_k+\ep_{-k}-\ep_0})+(\al_k-2\al_{-k})\mathcal D_K (x^{\al+\ep_k}).
\end{split}
\end{equation*}

For $\ell \geq 1$, we have
\begin{eqnarray*}\begin{split}
d^{(\ell+1)}\Big(\mathcal D_K (x^{\al})\Big)&=\frac{d}{\ell{+}1}d^{(\ell)}\Big(\mathcal D_K (x^\al)\Big) \\
&=\frac{1}{\ell{+}1}\sum\limits_{j=0}^{\ell}A_j B_{\ell-j} d \Big(\mathcal D_K (x^{\al+\ell(2\ep_k+\ep_{-k})-j(\ep_k+\ep_{-k})-(\ell-j)\ep_0})\Big)\\
&=\frac{1}{\ell{+}1}\sum\limits_{j=0}^{\ell}A_j B_{\ell-j}\Big((\ell{-}j{-}\al_0)\mathcal D_K  (x^{\al+(\ell+1)(2\ep_k+\ep_{-k})-j(\ep_k+\ep_{-k})-(\ell-j+1)\ep_0})\\
&\quad+\big(\al_k{+}2\ell{-}j{-}2(\al_{-k}{+}\ell{-}j)\big)\mathcal D_K (x^{\al+(\ell+1)(2\ep_k+\ep_{-k})-(j+1)(\ep_k+\ep_{-k})-(\ell-j)\ep_0})\Big)\\
&=\frac{1}{\ell{+}1}\sum\limits_{j=0}^{\ell}A_j (\ell{-}j{+}1)B_{\ell-j+1}\mathcal D_K (x^{\al+(\ell+1)(2\ep_k{+}\ep_{-k}){-}j(\ep_k+\ep_{-k})-(\ell{-}j{+}1)\ep_0})\\
&\quad+\frac{1}{\ell{+}1}\sum\limits_{j=0}^{\ell}(j{+}1)A_{j+1} B_{\ell-j} \mathcal D_K (x^{\al+(\ell+1)(2\ep_k+\ep_{-k})-(j{+}1)(\ep_k+\ep_{-k})-(\ell{-}j)\ep_0})
\\
&=\frac{1}{\ell{+}1}\sum\limits_{j=0}^{\ell}A_j (\ell{-}j{+}1)B_{\ell-j+1}\mathcal D_K (x^{\al+(\ell+1)(2\ep_k+\ep_{-k})-j(\ep_k+\ep_{-k})-(\ell-j+1)\ep_0})\\
&\quad+\frac{1}{\ell{+}1}\sum\limits_{j=1}^{\ell+1}jA_{j} B_{\ell-j+1} \mathcal D_K (x^{\al+(\ell+1)(2\ep_k+\ep_{-k})-j(\ep_k+\ep_{-k})-(\ell-j+1)\ep_0})\\
&=\sum\limits_{j=0}^{\ell{+}1} A_j B_{\ell-j+1}\mathcal D_K (x^{\al+(\ell+1)(2\ep_k+\ep_{-k})-j(\ep_k+\ep_{-k})-(\ell+1-j)\ep_0}).
\end{split}
\end{eqnarray*}

(\ref{vrelation4}) follows from the derivation property of $d^{(\ell)}$.

(\ref{vrelation5}):  recalling that for any elements $a$, $e$ in an associative algebra, one has
\begin{equation}
ca^m=\sum\limits_{\ell=0}^{m} (-1)^\ell \binom{m}{\ell}a^{m-\ell} (\ad a)^\ell (c). \label{vrelation7}
\end{equation}
Combing this with (\ref{vrelation4}), we can get (\ref{vrelation5}).

(\ref{vrelation6}): Use induction on $\ell$. When $\ell=1$, we have
$$
\ad \mathcal D_K (x^\al) \big( e \big)=\al_0 \mathcal D_K (x^{\al+2\ep_k+\ep_{-k}-\ep_0})+(2\al_{-k}-\al_k) \mathcal D_K (x^{\al+\ep_k}).
$$
When $\ell> 1$,
\begin{equation*}
\begin{split}
(\ad &\mathcal D_K (x^\al))^{\ell+1} \big( e \big)=\sum\limits_{j=0}^{\ell}\binom{\ell}{j} C_{\ell-j}
D(j,k)[\mathcal D_K (x^\al),\mathcal D_K (x^{\ell \al +2\ep_k+\ep_{-k}-j(\ep_k+\ep_{-k})-(\ell-j)\ep_0})]\\
&=\sum\limits_{j=0}^{\ell}\binom{\ell}{j} C_{\ell-j} D(j,k)
\Bigg(\Big((2{-}\sum\limits_{i=1}^{n}(\al_i{+}\al_{-i}))(\ell \al_0{-}(\ell{-}j))-\\
&\quad(2{-}\ell\sum\limits_{i=1}^{n}(\al_i{+}\al_{-i}){-}3{+}2j)\al_0\Big)\mathcal D_K (x^{(\ell+1)\al+2\ep_k+\ep_{-k}-j(\ep_k+\ep_{-k})-(\ell+1-j)\ep_0})\\
&\quad+\Big(\al_{-k}(\ell \al_k{+}2{-}j){-}\al_k(\ell \al_{-k}{+}1{-}j)\Big) \mathcal D_K (x^{(\ell+1)\al+2\ep_k+\ep_{-k}-(j+1)(\ep_k+\ep_{-k})-(\ell-j)\ep_0})\Bigg)\\
&=\sum\limits_{j=0}^{\ell}\binom{\ell}{j}C_{\ell-j+1} D(j,k)\mathcal D_K (x^{(\ell+1)
\al{+}2\ep_k{+}\ep_{-k}{-}j(\ep_k{+}\ep_{-k}){-}(\ell{+}j{-}1)\ep_0})\\
&\quad +\sum\limits_{j=0}^{\ell}\binom{\ell}{j} C_{\ell-j} D(j{+}1,k)\mathcal D_K (x^{(\ell+1)\al+2\ep_k+\ep_{-k}-(j+1)(\ep_k+\ep_{-k})-(\ell-j)\ep_0})\\
&=\sum\limits_{j=0}^{\ell}\binom{\ell}{j}C_{\ell-j+1} D(j,k)\mathcal D_K (x^{(\ell+1)\al+2\ep_k+\ep_{-k}-j(\ep_k+\ep_{-k})-
(\ell-j+1)\ep_0})\\
&\quad +\sum\limits_{j=1}^{\ell+1}\binom{\ell}{j{-}1} C_{\ell-j+1} D(j,k)\mathcal D_K (x^{(\ell+1)\al+2\ep_k+\ep_{-k}-j(\ep_k+\ep_{-k})-(\ell+1-j)\ep_0})\\
&=\sum\limits_{j=0}^{\ell+1}\binom{\ell{+}1}{j}C_{\ell-j+1} D(j,k)
\mathcal D_K (x^{(\ell+1)\al+2\ep_k+\ep_{-k}-j(\ep_k+\ep_{-k})-(\ell+1-j)\ep_0}).
\end{split}
\end{equation*}
This completes the proof.
\end{proof}

\begin{lemm}\label{vrelation8-10}
For $a \in \mathbb{F},~\al \in \mathbb{Z}^{2n+1}$, and $\mathcal D_K (x^\al) \in K$, the following equalities hold
\begin{eqnarray}
&&((\mathcal D_K (x^\al))^s \otimes 1) \cdot F_a=F_{a+s(\al_{-k}-\al_{k})} \cdot \big((\mathcal D_K (x^\al))^s \otimes 1\big),\label{vrelation8}\\
&&(\mathcal D_K (x^\al))^s \cdot u_a=u_{a+s(\al_k-\al_{-k})}\sum\limits_{\ell=0}^{\infty} d^{(\ell)}(\mathcal D_K (x^\al))^s h_{1-a}^{\lg \ell \rg} t^\ell,\label{vrelation9}\\
&&(1 \otimes (\mathcal D_K (x^\al))^s) \cdot F_a =\sum\limits_{\ell=0}^{\infty} (-1)^{\ell} F_{a+\ell} \big(h_{a}^{\lg\ell \rg} \otimes d^{(\ell)}(\mathcal D_K (x^{\al}))^s t^\ell \big).\label{vrelation10}
\end{eqnarray}
\end{lemm}
\begin{proof}
By (\ref{vrelation1}) and (\ref{vrelation2}), we have
\begin{equation*}
\begin{split}
\big((\mathcal D_K (x^\al))^s& \otimes 1\big) \cdot F_a = \big((\mathcal D_K (x^\al))^s \otimes 1\big) \Big(\sum\limits_{r=0}^{\infty} \frac{1}{r!} h_a^{\lg r \rg} \otimes e^r t^r\Big)\\
&=\sum\limits_{r=0}^{\infty} \frac{1}{r!} h_{a-s(\al_k-\al_{-k})}^{\lg m\rg } (\mathcal D_K (x^\al))^s \otimes e^rt^r=F_{a+s(\al_{-k}-\al_k)}\cdot \big( (\mathcal D_K (x^\al))^s \otimes 1\big).
\end{split}
\end{equation*}
Thus, we can obtain  (\ref{vrelation8}).

(\ref{vrelation9}): Use induction on $s$. When $s=1$, we have
\begin{equation*}
\begin{split}
\mathcal D_K &(x^\al) \cdot u_a=\mathcal D_K (x^\al) \sum\limits_{r=0}^{\infty}\frac{(-1)^r}{r!}h_{-a}^{[r]} e^r t^r\\
&\ =\sum\limits_{r=0}^{\infty} \frac{(-1)^r}{r!}h_{-a-(\al_k-\al_{-k})}^{[r]}\sum\limits_{\ell=0}^{r}(-1)^\ell \binom{r}{\ell}\ell!e^{r-\ell}\sum\limits_{j=0}^{\ell}A_j B_{\ell-j}\mathcal D_K (x^{\al +\ell(2\ep_k+\ep_{-k})-j(\ep_k+\ep_{-k})-(\ell-j)\ep_0}) t^{r}\\
&\ =\sum\limits_{r,\ell=0}^{\infty}\frac{(-1)^{r+\ell}}{(r+\ell)!} h_{-a-(\al_k-\al_{-k})}^{[r+\ell]} (-1)^\ell \binom{r{+}\ell}{\ell}\ell!e^r\sum\limits_{j=0}^{\ell}A_jB_{\ell-j} \mathcal D_K (x^{\al+\ell(2\ep_k+\ep_{-k})-j(\ep_k+\ep_{-k})-(\ell-j)\ep_0})t^{r+\ell}\\
&\ =\sum\limits_{r,\ell=0}^{\infty} \frac{(-1)^r}{r!} h_{-a-(\al_k-\al_{-k})}^{[r]} h_{-a-(\al_k-\al_{-k})-r}^{[\ell]} e^r \sum\limits_{j=0}^{\ell} A_j B_{\ell-j} \mathcal D_K (x^{\al+\ell(2\ep_k+\ep_{-k})-j(\ep_k+\ep_{-k})-(\ell-j)\ep_0})t^{r+\ell}\\
&\ =\sum\limits_{r,\ell=0}^{\infty} \frac{(-1)^r}{r!} h_{-a-(\al_k-\al_{-k})}^{[r]} e^r t^r h_{-a-(\al_k-\al_{-k})}^{[\ell]} \sum\limits_{j=0}^{\ell}A_j B_{\ell-j} \mathcal D_K (x^{\al+\ell(2\ep_k+\ep_{-k})-j(\ep_k+\ep_{-k})-(\ell-j)\ep_0}) t^\ell\\
&\ =u_{a+(\al_k-\al_{-k})}\sum\limits_{\ell=0}^{\infty} \sum\limits_{j=0}^{\ell}A_j B_{\ell-j}\mathcal D_K (x^{\al+\ell(2\ep_k+\ep_{-k})-j(\ep_k+\ep_{-k})-(\ell-j)\ep_0}) h_{-a+\ell}^{[\ell]} t^\ell\\
&\ =u_{a+(\al_k-\al_{-k})}\sum\limits_{\ell=0}^{\infty} d^{(\ell)}(\mathcal D_K (x^\al))~ h_{1-a}^{\lg \ell \rg} ~t^\ell.
\end{split}
\end{equation*}
When $s> 1$, we have
\begin{eqnarray*}
\begin{split}
(\mathcal D_K &(x^\al))^{s+1} \cdot u_a
=\mathcal D_K (x^\al) u_{a+s(\al_k-\al_{-k})}\sum\limits_{\ell=0}^{\infty} d^{(\ell)}(\mathcal D_K (x^\al))^s \cdot h_{1-a}^{\lg \ell \rg} t^\ell\\
&=u_{a+(s+1)(\al_k-\al_{-k})}\sum\limits_{\ell'=0}^{\infty} d^{(\ell')}\big(\mathcal D_K (x^\al)\big) h_{1-a-s(\al_k-\al_{-k})}^{\lg \ell' \rg} t^{\ell'} \sum\limits_{\ell=0}^{\infty} d^{(\ell)}(\mathcal D_K (x^\al))^s h_{1-a}^{\lg \ell \rg}t^\ell
\\
&=u_{a+(s+1)(\al_k-\al_{-k})}\sum\limits_{\ell'+\ell=0}^{\infty} d^{(\ell'+\ell)}(\mathcal D_K (x^\al))^{s+1} h_{1-a}^{\lg \ell+\ell' \rg} t^{\ell+\ell'}\\
&=u_{a+(s+1)(\al_k-\al_{-k})}\sum\limits_{\ell=0}^{\infty} d^{(\ell)}(\mathcal D_K (x^\al))^{s+1} h_{1-a}^{\lg \ell \rg} t^{\ell},
\end{split}
\end{eqnarray*}
where the third ``=" comes from the following equation, together with (\ref{vrelation3}) and (\ref{vrelation4}):
$$h_{1-a-s(\al_k-\al_{-k})}^{\lg \ell' \rg} \cdot d^{(\ell)}(\mathcal D_K (x^\al))^s=d^{(\ell)}(\mathcal D_K (x^\al))^s \cdot h_{1-a+\ell}^{\lg \ell' \rg}. $$

(\ref{vrelation10}):  If $s=1$, we have
\begin{equation*}
\begin{split}
(1 &\otimes \mathcal D_K (x^\al)) \cdot F_a =(1 \otimes \mathcal D_K (x^\al)) \sum\limits_{r=0}^{\infty} \frac{1}{r!}h_{a}^{\lg r \rg} \otimes e^r t^r\\
&\quad = \sum\limits_{r=0}^{\infty} \frac{1}{r!} h_a^{\lg r \rg} \otimes \sum\limits_{\ell=0}^{r} (-1)^\ell \binom{r}{\ell} \ell ! e^{r-\ell}\sum\limits_{j=0}^{\ell} A_j B_{\ell-j} \mathcal D_K (x^{\al+\ell(2\ep_k+\ep_{-k})-j(\ep_k+\ep_{-k})-(\ell-j)\ep_0})t^r\\
&\quad = \sum\limits_{r,\ell=0}^{\infty}(-1)^\ell \big(\frac{1}{r!} h_{a+\ell}^{\lg r \rg} \otimes e^r t^r  \big) \big(h_a^{\lg \ell \rg} \otimes \sum\limits_{j=0}^{\ell} A_j B_{\ell-j} \mathcal D_K (x^{\al+\ell(2\ep_k+\ep_{-k})-j(\ep_k+\ep_{-k})-(\ell-j)\ep_0} t^\ell)\big)\\
&\quad=\sum\limits_{\ell=0}^{\infty} (-1)^\ell F_{a+\ell} h_{a}^{\lg \ell \rg} \otimes d^{(\ell)} (\mathcal D_K (x^\al)) t^\ell.
\end{split}
\end{equation*}
If $s> 1$, we have
\begin{eqnarray*}
&&(1 \otimes (\mathcal D_K (x^\al))^{s+1}) \cdot F_a=(1 \otimes \mathcal D_K (x^\al))(1 \otimes (\mathcal D_K (x^\al))^s)\cdot F_a\\
&&\ =(1 \otimes \mathcal D_K (x^\al))\sum\limits_{\ell=0}^{\infty} F_{a+\ell}(-1)^\ell h_{a}^{\lg \ell \rg} \otimes d^{(\ell)}(\mathcal D_K (x^\al))^s) t^\ell\\
&&\ =\sum\limits_{\ell=0}^{\infty}\sum\limits_{\ell'=0}^{\infty}(-1)^\ell F_{a+\ell+\ell'}(-1)^{\ell'} (h_{a+\ell}^{\lg \ell' \rg} \otimes d^{(\ell')}(\mathcal D_K (x^\al)) t^{\ell'} )\big( h_{a}^{\lg \ell \rg} \otimes d^{(\ell)} (\mathcal D_K (x^\al))^s t^\ell\big)\\
&&\ =\sum\limits_{\ell,\ell'=0}^{\infty}(-1)^{\ell+\ell'}F_{a+\ell+\ell'} h_{a+\ell}^{\lg \ell' \rg} h_{a}^{\lg \ell \rg} \otimes d^{(\ell')}\big(\mathcal D_K (x^\al)\big) d^{(\ell)} (\mathcal D_K (x^\al))^s t^{\ell+\ell'}\\
&&\ =\sum\limits_{\ell=0}^{\infty}(-1)^\ell F_{a+\ell} h_{a}^{\lg \ell \rg} \otimes d^{(\ell)} (\mathcal D_K (x^\al))^{s+1} t^\ell.
\end{eqnarray*}
Thus, (\ref{vrelation10}) holds by induction on $s$.
This completes the proof.
\end{proof}

\begin{lemm}
For $s \geq 1$, we have
\begin{gather}
\Delta\big((\mathcal D_K (x^\al))^s\big)=\sum\limits_{0 \leq j \leq s \atop \ell \geq 0} (-1)^\ell \binom{s}{j} (\mathcal D_K (x^\al))^j h^{\lg \ell \rg} \otimes (1-et)^{j(\al_k-\al_{-k})-\ell}\big(d^{(\ell)}  (\mathcal D_K (x^\a))^{s-j}\big) t^\ell, \label{vrelation11}\\
S((\mathcal D_K (x^\al))^s)=(-1)^s (1-et)^{-s(\al_k-\al_{-k})}\sum\limits_{\ell=0}^{\infty} d^{(\ell)} (\mathcal D_K (x^\al))^s h_{1}^{\lg \ell \rg} t^\ell .\label{vrelation12}
\end{gather}
\end{lemm}

\begin{proof}
By Lemma \ref{vrelation8-10}, we have
\begin{equation*}
\begin{split}
&\Delta((\mathcal D_K (x^\al))^s)
=\mathcal{F}\sum\limits_{j=0}^{s} \binom{s}{j} \big(\mathcal D_K (x^\al) \otimes 1\big)^j F \mathcal{F}\big(1 \otimes (\mathcal D_K (x^\al))^{s-j}\big) F\\
&=\sum\limits_{j=0}^{s} \binom{s}{j} \mathcal{F} F_{j(\al_{-k}-\al_{k})} (\mathcal D_K (x^\al) \otimes 1)^j
\big( \mathcal{F} \sum\limits_{\ell=0}^{\infty} (-1)^\ell F_\ell h^{\lg \ell \rg} \otimes d^{(\ell)}(\mathcal D_K (x^\al))^{s-j}\big) t^\ell\\
&=\sum\limits_{j=0}^{s} \binom{s}{j}\big(1 {\otimes} (1{-}et)^{j(\al_k-\al_{-k})}) \big(\mathcal D_K (x^\al))^j {\otimes} 1 \big)
\big(\sum\limits_{\ell=0}^{\infty} ({-}1)^\ell (1 {\otimes} (1{-}et)^{-\ell}) (h^{\lg \ell \rg} {\otimes} d^{(\ell)} (\mathcal D_K (x^\al))^{s-j}) t^\ell\big)
\end{split}
\end{equation*}\begin{equation*}
\begin{split}
&=\sum\limits_{0 \leq j \leq s \atop \ell \geq 0} (-1)^\ell \binom{s}{j} (\mathcal D_K (x^\al))^j h^{\lg \ell \rg} \otimes (1{-}et)^{j(\al_k-\al_{-k})-\ell} d^{(\ell)}(\mathcal D_K (x^\al))^{s-j} t^\ell,\\
&S((\mathcal D_K (x^\al))^s)=v(-1)^s (\mathcal D_K (x^\al))^s u
=(-1)^s vu_{s(\al_k-\al_{-k})} \sum\limits_{\ell=0}^{\infty} d^{(\ell)}(\mathcal D_K (x^\al))^s) ~h_{1}^{\lg \ell \rg} t^\ell\\
&=(-1)^s (1{-}et)^{-s(\al_k-\al_{-k})}\sum\limits_{\ell=0}^{\infty} d^{(\ell)}(\mathcal D_K (x^\al))^s) ~h_{1}^{\lg \ell \rg} t^\ell.
\end{split}
\end{equation*}

This completes the proof.
\end{proof}

\begin{theorem}\label{vHopf1}
Fix two distinguished elements $ h=\mathcal D_K (x^{\ep_k+\ep_{-k}}),~e=\mathcal D_K (x^{2\ep_k+\ep_{-k}})$, such that $[h,e]$ $=e$ in the generalized Cartan type $K$ Lie algebra ${\bf K}$ over $\mathbb{F}$. There exists a structure of noncommutative and noncocommutative Hopf algebra $(U({\bf K})[[t]]~,m,~\iota,~\Delta,~S,~\varepsilon)$ 
which leaves the product of ~$(U({\bf K})[[t]]$ undeformed but with the deformed coproduct, the antipode and the counit defined by:
\begin{gather}
\Delta(\mathcal D_K (x^\al))=\mathcal D_K (x^\al) \otimes (1-et)^{\al_k-\al_{-k}}+\sum\limits_{\ell=0}^{\infty} (-1)^\ell h^{\lg \ell \rg} \otimes (1{-}et)^{-\ell} d^{(\ell)} (\mathcal D_K (x^\al)) t^\ell, \label{vHopf1d}\\
S(\mathcal D_K (x^\al))=-(1{-}et)^{\al_{-k}-\al_k}\sum\limits_{\ell=0}^{\infty} d^{(\ell)}(\mathcal D_K (x^\al)) h_{1}^{\lg \ell \rg} t^\ell,\hskip4cm \label{vHopf1a}
\end{gather}
and $\varepsilon(\mathcal D_K (x^\al))=0$, for any $\mathcal D_K (x^\al) \in {\bf K}$.
\end{theorem}

\begin{proof} By Lemma \ref{vrelation8-10}, we have
\begin{eqnarray*}
\begin{split}
&\Delta(\mathcal D_K (x^\al))
=\mathcal{F}(\mathcal D_K (x^\al)\otimes 1) F+ \mathcal{F}(1 \otimes \mathcal D_K (x^\al)) F\\
&\quad=\mathcal{F}F_{\al_{-k}-\al_k}(\mathcal D_K (x^\al) \otimes 1)+\mathcal{F}\sum\limits_{\ell=0}^{\infty} (-1)^\ell F_\ell h^{\lg \ell \rg} \otimes d^{(\ell)}\Big(\mathcal D_K (x^\al)\Big) t^\ell\\
&\quad=\mathcal D_K (x^\al) \otimes (1{-}et)^{\al_k-\al_{-k}}+\sum\limits_{\ell=0}^{\infty} (-1)^\ell h^{\lg \ell \rg} \otimes (1{-}et)^{-\ell} d^{(\ell)}\big(\mathcal D_K (x^\al)\big) t^\ell,\\
&S(\mathcal D_K (x^\al))=u^{-1} S_0(\mathcal D_K (x^\al)) u=-v \mathcal D_K (x^\al) u=-vu_{\al_k-\al_{-k}} \sum\limits_{\ell=0}^{\infty} d^{(\ell)}(\mathcal D_K (x^\al)) h_{1}^{\lg \ell \rg}t^\ell\\
&\quad=-(1{-}et)^{-(\al_k-\al_{-k})} \sum\limits_{\ell=0}^{\infty} d^{(\ell)}(\mathcal D_K (x^\al)) h_{1}^{\lg \ell \rg}t^\ell.
\end{split}
\end{eqnarray*}

This completes the proof.\end{proof}

Note that $\{\, \mathcal D_K (x^{\al}) \mid \al \in \mathbb{Z}_+^{2n+1}\,\}$ is a $\mathbb{Z}$-basis of
${\bf K}_{\mathbb{Z}}^+$ as a subalgebra of ${\bf K}_{\mathbb{Z}}$ and ${\bf W}_{\mathbb{Z}}^+$. Consequently, we have

\begin{coro}\label{vHopf2}
Fix distinguished elements $h=\mathcal D_K (x^{\ep_k+\ep_{-k}}),~e=\mathcal D_K (x^{2\ep_k+\ep_{-k}}), ~1 \leq k \leq n$, the corresponding quantization of $U({\bf K}_{\mathbb{Z}}^{+})$ over $U({\bf K}_{\mathbb{Z}}^{+})[[t]]$ by Drinfeld'd twist $\mathcal{F}$ with the product undeformed is given by:
\begin{eqnarray}
&&\Delta(\mathcal D_K (x^\al))=\mathcal D_K (x^\al) \otimes (1{-}et)^{\al_k-\al_{-k}} \nonumber\\
&&\qquad\qquad\quad+ \sum\limits_{\ell=0}^{\infty}(-1)^\ell h^{\lg \ell \rg} \otimes (1{-}et)^{-\ell} \sum\limits_{j=0}^{\ell} A_j B_{\ell-j}\mathcal D_K (x^{\al+\ell(2\ep_k+\ep_{-k})-j(\ep_k+\ep_{-k})-(\ell-j)\ep_0}) t^\ell, \label{vHopf2d}\\
&&S(\mathcal D_K (x^\al))=-(1{-}et)^{\al_{-k}-\al_k}\sum\limits_{\ell=0}^{\infty} \sum\limits_{j=0}^{\ell} A_j B_{\ell-j}\mathcal D_K (x^{\al+\ell(2\ep_k+\ep_{-k})-j(\ep_k+\ep_{-k})-(\ell-j)\ep_0}) h_{1}^{\lg \ell \rg}
t^\ell ,\label{vHopf2a}
\end{eqnarray}
and $\varepsilon(\mathcal D_K (x^\al))=0$, where $A_j=\frac{1}{j!}\prod\limits_{i=0}^{j-1}(\al_k{-}2\al_{-k}{+}i),
~B_{\ell-j}=\frac{1}{(\ell-j)!}\prod\limits_{i=0}^{\ell-j-1}(i{-}\al_0)$, with $A_0=B_0=1$, $A_{-1}=B_{-1}=0$.
\end{coro}

\begin{proof}
By Theorem \ref{vHopf1} and formula (\ref{vrelation3}), we can get  formula (\ref{vHopf2d}), and (\ref{vHopf2a}).
By Lemma \ref{relation3}, the coefficients in the two formulas are all integers.
\end{proof}
\subsection{Quantization of the Contact algebra ${\bold K}(2n{+}1; \underline{1})$}

 In this subsection, firstly, we make
{\it modulo $p$ reduction and base change with $\mathcal K[[t]]$
replaced by $\mathcal K[t]$}, for the quantization of
$U(\mathbf{K}^+_{\mathbb{Z}})$ in characteristic $0$ (Corollary \ref{vHopf2}) to yield
the quantization of $U(\mathbf{K}(2n{+}1;\underline{1}))$, for the
restricted simple modular Lie algebra $\mathbf{K}(2n{+}1;\underline{1})$
in characteristic $p$. Secondly, we shall further make {\it
``$p$-restrictedness" reduction as well as base change with
$\mathcal K[t]$ replaced by $\mathcal K[t]_p^{(q)}$}, for the
quantization of $U(\mathbf{K}(2n{+}1;\underline{1}))$, which will lead
to the required quantization of
$\mathbf{u}(\mathbf{K}(2n{+}1;\underline{1}))$, the restricted universal
enveloping algebra of $\mathbf{K}(2n{+}1;\underline{1})$.

Let $\mathbb{Z}_p$ be the prime subfield of $\mathcal{K}$ with
$\text{char}\,(\mathcal{K})=p$. When considering
$\mathbf{W}_{\mathbb{Z}}^+$ as a $\mathbb{Z}_p$-Lie algebra, namely,
making modulo $p$ reduction for the defining relations of
$\mathbf{W}_{\mathbb{Z}}^+$, denoted by
$\mathbf{W}_{\mathbb{Z}_p}^+$, we see that
$(J_{\underline{1}})_{\mathbb{Z}_p}=\text{Span}_{\mathbb{Z}_p}\{x^\alpha
D_i \mid \exists\, j: \alpha_j\ge p\,\}$ is a maximal ideal of
$\mathbf{W}^+_{\mathbb{Z}_p}$, and
$\mathbf{W}^+_{\mathbb{Z}_p}/(J_{\underline{1}})_{\mathbb{Z}_p}
\cong \mathbf{W}(2n{+}1;\underline{1})_{\mathbb{Z}_p}
=\text{Span}_{\mathbb{Z}_p}\{x^{(\alpha)}D_{i}\mid 0\le \alpha\le
\tau, -n\le i\le n\}$. For the subalgebra
$\mathbf{K}_{\mathbb{Z}}^+$, we have
$\mathbf{K}^+_{\mathbb{Z}_p}/(\mathbf{K}^+_{\mathbb{Z}_p}\cap(J_{\underline{1}})_{\mathbb{Z}_p})
\cong \mathbf{K}'(2n{+}1;\underline{1})_{\mathbb{Z}_p}$. We simply denote
$\mathbf{K}^+_{\mathbb{Z}_p}\cap(J_{\underline{1}})_{\mathbb{Z}_p}$
as  $(J^+_{\underline{1}})_{\mathbb{Z}_p}$.

Moreover, we have $\mathbf{K}'(2n{+}1;\underline{1})
=\mathcal{K}\otimes_{\mathbb{Z}_p}\mathbf{K}'(2n{+}1;\underline{1})_{\mathbb{Z}_p}
=\mathcal{K}\mathbf{K}'(2n{+}1;\underline{1})_{\mathbb{Z}_p}$, and
$\mathbf{K}^+_{\mathcal{K}}=\mathcal{K}\mathbf{K}^+_{\mathbb{Z}_p}$.

Observe that the ideal
$J^+_{\underline{1}}:=\mathcal{K}(J^+_{\underline{1}})_{\mathbb{Z}_p}$
generates an ideal of $U(\mathbf{K}^+_{\mathcal{K}})$ over
$\mathcal{K}$, denoted by
$J:=J^+_{\underline{1}}U(\mathbf{K}^+_{\mathcal{K}})$, where
$\mathbf{K}^+_{\mathcal{K}}/J^+_{\underline{1}}\cong
\mathbf{K}'(2n{+}1;\underline{1})$. Based on the formulae of Corollary
\ref{vHopf2}, $J$ is a Hopf ideal of $U(\mathbf{K}^+_{\mathcal{K}})$
satisfying $U(\mathbf{K}^+_{\mathcal{K}})/J\cong
U(\mathbf{K}'(2n{+}1;\underline{1}))$. Note that elements $\sum a_{
\alpha}\frac{1}{\alpha!}\mathcal D_K ( x^{\alpha})$ in
$\mathbf{K}^+_{\mathcal{K}}$ for $0\le\alpha\le\tau$ will be
identified with $\sum a_{\alpha}\mathcal D_K ( x^{(\alpha)})$ in
$\mathbf{K}'(2n{+}1;\underline{1})$ and those in $J_{\underline{1}}$
with $0$. Hence, by  Corollary \ref{vHopf2}, we get the quantization of
$U\big(\mathbf{K}'(2n{+}1;\underline{1})\big)$ over
$U_t(\mathbf{K}'(2n{+}1;\underline{1})):=U(\mathbf{K}'(2n{+}1;\underline{1}))\otimes_{\mathcal
K}\mathcal K[t]$ (not necessarily in
$U(\mathbf{K}'(2n{+}1;\underline{1}))[[t]]$, as seen in formulae (\ref{vHopf3d}) \& (\ref{vHopf3a})
 as follows.

\begin{theorem}\label{vHopf3}
Fix two distinguished elements $h:=\mathcal D_K ( x^{(\ep_k
+\ep_{-k})}),~e:=2\mathcal D_K ( x^{(2\ep_k+\ep_{-k})})$ $(1 \leq k \leq n)$; the
corresponding quantization of $U(\mathbf{K}^\prime(2n{+}1;\underline{1}))$ over $U_t(\mathbf{K}^\prime
(2n{+}1;\underline{1}))$ with the product undeformed is given by
\begin{eqnarray}
&&\Delta\big(\mathcal D_K ( x^{(\alpha)})\big)=\mathcal D_K ( x^{(\alpha)}) \otimes
(1{-}et)^{\alpha_k-\alpha_{-k}}+\sum\limits_{\ell=0}^{p-1}\sum\limits_{j=0}^{\ell}(-1)^\ell\overline{A}_{\ell} \overline{B}_{\ell-j}
h^{\lg \ell \rg}\otimes(1{-}et)^{-\ell}\cdot \label{vHopf3d}\\
&&\hskip6cm
\cdot\mathcal D_K (x^{(\alpha +\ell(2\ep_k+\ep_{-k})-j(\ep_k+\ep_{-k})-(\ell-j)\ep_0)}) t^\ell, \nonumber\\
&&S\big( \mathcal D_K ( x^{(\alpha)})\big)=-(1{-}et)^{\alpha_{-k} -\alpha_{k}}\sum\limits_{\ell=0}^{p-1} \sum\limits_{j=0}^{\ell} \overline{A}_j \overline{B}_{\ell-j} \mathcal D_K (x^{(\alpha +\ell(2\ep_k+\ep_{-k})-j(\ep_k+\ep_{-k})-(\ell-j)\ep_0)}) h_1^{\lg \ell \rg } t^\ell, \label{vHopf3a}
\end{eqnarray}
and $\varepsilon\big(\mathcal D_K ( x^{(\alpha)})\big)=0$, where $\overline{A}_j=(2\ell{-}j)!\binom{\al_k{+}(2\ell{-}j)}{2\ell{-}j}A_j$, with $A_j=\frac{1}{j!}\prod\limits_{i=0}^{j-1} (\al_k{-}2\al_{-k}{+}i),~ A_0=1,~A_{-1}=0$. For $0 \leq \ell{-}j \leq \al_0,\overline{B}_{\ell-j}=(-1)^{\ell-j}\binom{\al_{-k}+(\ell-j)}{\ell-j}$, otherwise, $\overline{B}_{\ell-j}=0$.
\end{theorem}
\begin{proof}
Note that the elements $\frac{1}{\al !} \mathcal D_K
(x^{\al})$ in $\mathbf{K}_{\mathcal{K}}^{+}$ will be identified with
 $\mathcal D_K ({x^{(\al)}})$ in $\mathbf{K}(2n{+}1;1)$ and those
in $J_{\underline{1}}$ with 0. Hence, by Corollary \ref{vHopf2}, we can get
\begin{equation*}
\begin{split}
&\Delta(\mathcal D_K (x^{(\al)}))=\frac{1}{\al !}\Delta(\mathcal D_K (x^\al))\\
&=\mathcal D_K (x^{(\al)}) \otimes (1{-}et)^{\al_k-\al_{-k}}+\sum\limits_{\ell=0}^{p-1}(-1)^\ell h^{\lg \ell \rg} \otimes (1{-}et)^{-\ell}\times\\
&\qquad \sum\limits_{j=0}^{\ell}A_j B_{\ell-j}\frac{ (\alpha {+}\ell(2\ep_k{+}\ep_{-k}){-}j(\ep_k{+}\ep_{-k}){-}(\ell{-}j)\ep_0)!}{\al!}\mathcal D_K (x^{(\alpha +\ell(2\ep_k+\ep_{-k})-j(\ep_k+\ep_{-k})-(\ell-j)\ep_0)}) t^\ell,
\end{split}
\end{equation*}
where
\begin{equation*}
\begin{split}
&A_j B_{\ell-j}\frac{ (\alpha {+}\ell(2\ep_k{+}\ep_{-k}){-}j(\ep_k{+}\ep_{-k}){-}(\ell{-}j)\ep_0)! }{\al!}\\
&=A_j  B_{\ell-j} \frac{ (\al_{-k}{+}\ell{-}j)!(\al_k{+}2\ell{-}j)!(\al_0{-}(\ell{-}j))!}{\al_{-k}! \al_k!\al_0!}\\
&=A_j \frac{(\al_k{+}2\ell{-}j)!}{\al_k !}\cdot \frac{\prod\limits_{i=0}^{\ell-j-1}(i{-}\al_0)}{(\ell{-}j)!}\cdot\frac{ (\al_{-k}{+}\ell{-}j)!(\al_0{-}(\ell{-}j))!}{\al_{-k}!\al_0!}\\
&=\overline{A}_j \cdot (-1)^{\ell-j} \binom{\al_{-k}{+}\ell{-}j}{\ell{-}j}=\overline{A}_j ~\overline{B}_{\ell-j}.
\end{split}
\end{equation*}

Hence, we can get (\ref{vHopf3d}). Another formula can be obtained similarly.
\end{proof}

Note that for $\al\le\tau$, if there exist $0< \ell \leq p-1,~0 \leq j \leq \ell$, such that
 $\al {+}\ell(2\ep_k{+}\ep_{-k}){-}j(\ep_k$ ${+}\ep_{-k}){-}(\ell{-}j)\ep_0=\tau$, that is,
$$ \al_0{-}(\ell{-}j)=p{-}1;\quad \al_k{+}2\ell{-}j=p{-}1;\quad \al_{-k}{+}(\ell{-}j)=p{-}1.$$
This means that
$\al_0=p{-}1,\ j=\ell,\ \al_k{+}\ell=p{-}1,\ \al_{-k}=p{-}1$.
Thus, the coefficients of $\mathcal D_K (x^{(\tau)})$ is
$\overline{A}_\ell \overline{B}_0=\ell!\binom{\al_k+\ell}{\ell}A_\ell
=\binom{p-1}{\ell}\prod\limits_{i=0}^{\ell-1}(\al_k{-}2\al_{-k}{+}i)\equiv 0 \;(\mod\; p)$.
So whenever $2n{+}4 \equiv 0\; (\mod\; p)$, or $\nequiv 0\; (\mod\; p)$, Theorem \ref{vHopf3} gives the quantization of $U({\bf K}(2n{+}1;\underline{1}))$ over $U_t(\mathbf{K}(2n{+}1;\underline{1}))=U(\mathbf{K}(2n{+}1;\underline{1}))\otimes_{\mathcal{K}}
\mathcal{K}[t]$ and also over $U(\mathbf{K}(2n{+}1;\underline{1}))[[t]]$. It should be noticed
that in this step including from the quantization integral form of $U(\mathbf{K}_{\mathbb{Z}}^{+})$
and making the modulo $p$ reduction, we used the first base change with $\mathcal{K}[[t]]$ replaced by $\mathcal{K}[t]$, and
the objects from $U(\mathbf{K}(2n{+}1;\underline{1}))[[t]]$ turning to $U_t(\mathbf{K}(2n{+}1;\underline{1}))$.

Assume $0 \leq \alpha \leq
\tau $ when $2n{+}4 \nequiv 0 \;(\mod\, p)$, and $0 \leq \al <\tau $ when $2n{+}4 \equiv 0 \;(\mod\, p)$. Denote by $I$ the ideal
of $U(\mathbf{K}(2n{+}1;\underline{1}))$ over $\mathcal{K}$ generated by $\big(\mathcal D_K (x^{(\alpha)})\big)^p$ and
 $\big(\mathcal D_K (x^{(\ep_k +\ep_{-k})})\big)^p-\mathcal D_K ( x^{(\ep_k
+\ep_{-k})}),~\big(\mathcal D_K (x^{(\ep_0)})\big)^p-\mathcal D_K (x^{(\ep_0)})$ with $\alpha \neq \ep_k +\ep_{-k},~\al \neq \ep_0$ for $1 \leq k \leq n$. Thus $\mathbf{u}\big(\mathbf{K}(2n{+}1;\underline{1})\big)=U\big(\mathbf{K}(2n{+}1;\underline{1})\big)/I$ is
of prime-power dimension $p^{p^{2n+1}}$ when $2n{+}4 \nequiv 0\;(\mod\, p)$; and of prime-power dimension $p^{p^{2n+1}-1}$ otherwise. In order to get a
reasonable quantization of prime-power dimension
for $\mathbf{u}\big(\mathbf{K}(2n{+}1;\underline{1})\big)$ in
characteristic $p$, at first, it is necessary to clarify in
concept what is the underlying vector space where the required $t$-deformed object exists. According to our modular reduction
approach, it should start to be induced from the
 $\mathcal{K}[t]$-algebra $U_t(\mathbf{K}(2n{+}1;\underline{1}))$ in
Theorem \ref{vHopf3}.

Firstly, we observe the following facts (for the proof of \cite{HW1} or \cite{HW2}):

\begin{lemm}\label{relation} $(\text{\rm i})$ \ $(1-et)^p\equiv 1 \quad (\text{\rm mod}\,p, I)$.

\smallskip $(\text{\rm ii})$ \ $(1-et)^{-1}\equiv
1+et+\cdots+e^{p-1}t^{p-1} \quad (\text{\rm mod}\,
 p, I)$.

\smallskip $(\text{\rm iii})$ \ $h_a^{\lg \ell\rg} \equiv 0 \quad
(\text{\rm mod} \, p, I)\;$ for $\ell \geq p$, and $a\in\mathbb{Z}_p$.
\end{lemm}

The above Lemma, together with Theorem \ref{vHopf3}, shows
that the required $t$-deformation of
$\mathbf{u}(\mathbf{K}(2n{+}1;\underline{1}))$ (if it exists) in fact
only happens in a $p$-truncated polynomial ring (with degrees of $t$
less than $p$) with coefficients in
$\mathbf{u}(\mathbf{K}(2n{+}1;\underline{1}))$, i.e.,
$$\mathbf{u}_{t,q}(\mathbf{K}(2n{+}1;\underline{1})):=
\mathbf{u}(\mathbf{K}(2n{+}1;\underline{1}))\otimes_{\mathcal K}
\mathcal{K}[t]_p^{(q)}$$ (rather than in
$\mathbf{u}_t(\mathbf{K}(2n{+}1;\underline{1})):=\mathbf{u}(\mathbf{K}(2n{+}1;\underline{1}))
\otimes_{\mathcal K}\mathcal K[t]$), where $\mathcal{K}[t]_p^{(q)}$
is conveniently taken to be a $p$-truncated polynomial ring which is
a quotient of $\mathcal K[t]$ defined as
$$
\mathcal{K}[t]_p^{(q)}= \mathcal{K}[t]/(t^p-qt), \qquad\text{\it for
}\ q\in\mathcal{K}.\label{truncated}
$$
Thereby, we obtain the underlying ring for our required
$t$-deformation of $\mathbf{u}(\mathbf{K}(2n{+}1;\underline{1}))$ over
$\mathcal{K}[t]_p^{(q)}$, and
$$\dim_{\mathcal{K}}\mathbf{u}_{t,q}(\mathbf{K}(2n{+}1;\underline{1}))
=p\cdot\dim_{\mathcal{K}}\mathbf{u}(\mathbf{K}(2n{+}1;\underline{1}))
=\begin{cases} p^{p^{2n+1}+1},~&2n{+}4\nequiv 0\;(\mod\, p),\\
 p^{p^{2n+1}},~&2n{+}4\equiv 0\;(\mod\, p).
 \end{cases}
$$

Via modulo ``$p$-restrictedness" reduction, it is
necessary for us to work over the objects from $U_t(\mathbf
K(2n{+}1;\underline{1}))$ passage to $U_{t,q}(\mathbf K(2n{+}1;\underline1))$
first, and then to $\mathbf u_{t,q}(\mathbf K(2n{+}1;\underline1))$ (see
the proof of Theorem \ref{vHopf4} below), here we used the second
base change with $\mathcal K[t]_p^{(q)}$ instead of $\mathcal K[t]$.

As the definition of \cite{HW2}, we gave the following description:

\begin{defi}\label{fdquan}  A Hopf algebra
$(\mathbf{u}_{t,q}(\mathbf{K}(2n{+}1;\underline{1}))$, $m,
\iota,\Delta,S,\varepsilon)$ over a ring $\mathcal K[t]_p^{(q)}$ of
characteristic $p$ is said to be a finite-dimensional quantization
of $\mathbf{u}(\mathbf{K}(2n{+}1;\underline{1}))$ if it is obtained from the standard Hopf algebra $U(\mathbf {K}^+_\mathbb Z)[[t]]$
by means of a Drinfel'd twisting, modular reductions and shrinking of base rings, and if there is
an isomorphism as algebras
$\mathbf{u}_{t,q}(\mathbf{K}(2n{+}1;\underline{1}))/t\mathbf{u}_{t,q}(\mathbf{K}(2n{+}1;\underline{1}))
$ $\cong \mathbf{u}(\mathbf{K}(2n{+}1;\underline{1}))$.
\end{defi}

To describe $\mathbf{u}_{t,q}(\mathbf{K}(2n{+}1;\underline{1}))$
explicitly, we still need an auxiliary Lemma.

\begin{lemm}\label{vrelation13}
Let $e=2\mathcal D_K (x^{(2\ep_k+\ep_{-k})})$, and $d^{(\ell)}=\frac{1}{\ell!} \ad e$, then:
\begin{enumerate}
    \item[(i)~~] $d^{(\ell)} \Big(\mathcal D_K (x^{(\al)})\Big)=\sum\limits_{j=0}^{\ell} \overline{A}_j \overline{B}_{\ell-j} \mathcal D_K (x^{(\al+\ell(2\ep_k+\ep_{-k})-j(\ep_k+\ep_{-k})-(\ell-j)\ep_0)})$, where $\overline{A}_j,~ \overline{B}_{\ell-j}$ as in Theorem \ref{vHopf3}.
    \item[(ii)~] $d^{(\ell)}  \Big(\mathcal D_K (x^{(\ep_i+\ep_{-i})})\Big)=\delta_{\ell 0}\mathcal D_K (x^{(\ep_i+\ep_{-i})})-\delta_{\ell 1}\delta_{ik}e, ~1 \leq i \leq n$,\\
                  $d^{(\ell)}  \Big(\mathcal D_K (x^{(\ep_0)})\Big)=\delta_{\ell 0}\mathcal D_K (x^{(\ep_0)})-\delta_{\ell 1}e$.
     \item[(iii)~]  $d^{(\ell)} \Big(\mathcal D_K (x^{(\al)})\Big)^p=\delta_{\ell 0} \Big(\mathcal D_K (x^{(\al)})\Big)^p-\delta_{\ell1}\big(\delta_{\al,\ep_k+\ep_{-k}}+\delta_{\al,\ep_0}\big)e$.
\end{enumerate}
\end{lemm}

\begin{proof}
(i) By (\ref{vrelation3}) and the proof of Theorem \ref{vHopf3}, we can get (i).

(ii) When $\ell=1,~d^{(\ell)}  \big(\mathcal D_K (x^{(\ep_i+\ep_{-i})})\big)=\overline{A}_0 \overline{B}_1 \mathcal D_K (x^{(\ep_i+\ep_{-i}+2\ep_k+\ep_{-k}-\ep_0)})+\overline{A}_1 \overline{B}_0 \mathcal D_K (x^{(\ep_i+\ep_{-i}+\ep_k)})$.

For $\al=\ep_i{+}\ep_{-i},~\al_0=0$, we have: $\overline{B}_0=1,\overline{B}_1=0$; $\overline{A}_0=2\binom{\delta_{ik}+2}{2},~\overline{A}_1=\binom{\delta_{ik}+1}{1} A_1=(\delta_{ik}+1)\big((\ep_i+\ep_{-i})_k-2(\ep_i+\ep_{-i})_{-k}\big)=-2\delta_{ik}.$
 Thus $d  \big(\mathcal D_K (x^{(\ep_i+\ep_{-i})})\big)=-2\delta_{ik} x^{(\ep_i+\ep_{-i}+\ep_k)}=-\delta_{ik}e$. Hence, we prove the first equation of (ii).

For $\al=\ep_0,~\al_0=1$, we can see that $\overline{A}_0=2,~\overline{A}_1=0$; $\overline{B}_0=1,~\overline{B}_1=(-1)^1 \binom{0+1}{1} =-1$. Thus, $d \big(\mathcal D_K (x^{(\ep_0)})\big)=-2\mathcal D_K (x^{(\ep_0+2\ep_k+\ep_{-k}-\ep_0)})=-e$. This completes the second equation of (ii).

(iii) By Lemma \ref{vrelation3}, (\ref{vrelation6}) and (\ref{vrelation7}), we have
\begin{equation*}
\begin{split}
&d\Big((\mathcal D_K (x^{(\al)}))^p\Big)=[e,(\mathcal D_K (x^{(\al)}))^p]\\
&=\sum\limits_{\ell=1}^{p}(-1)^\ell \binom{p}{\ell} (\mathcal D_K (x^{(\al)}))^{p-\ell} (\ad \mathcal D_K (x^{(\al)}))^\ell( e )\\
& \equiv (-1)^p  (\ad \mathcal D_K (x^{(\al)}))^p ( e )\quad(\mod\; p)\\
&\equiv -\frac{1}{(\al !)^p}  (\ad \mathcal D_K (x^{\al}))^p \Big(\mathcal D_K (x^{(2\ep_k+\ep_{-k})})\Big) \quad(\mod\; p)\\
&\equiv -\frac{1}{(\al !)^p} \sum\limits_{j=0}^{p} \binom{p}{j}C_{p-j} D(j,k) \mathcal D_K (x^{p\al+2\ep_k+\ep_{-k}-j(\ep_k+\ep_{-k})-(p-j)\ep_0}) \quad(\mod\; p)\\
&\equiv -\frac{1}{(\al !)^p}C_p\mathcal D_K (x^{p\al+2\ep_k+\ep_{-k}-p\ep_0})
 -\frac{1}{(\al !)^p} D(p, k)\mathcal D_K (x^{p\al+2\ep_k+\ep_{-k}-p(\ep_k+\ep_{-k})}) \quad(\mod\; p)\\
&\equiv \begin{cases} -e,  &~~\al=\ep_k+\ep_{-k},~\textit{or }~\ep_0, \\
                       0,  &~~\textit{otherwise},
        \end{cases}(\text{mod\; }p,\,J).
\end{split}
\end{equation*}

This completes the proof.
\end{proof}

Based on Theorem \ref{vHopf3}, Definition \ref{fdquan} and
Lemma \ref{vrelation13}, we arrive at
\begin{theorem}\label{vHopf4}
Fix two distinguished elements
$h:=\mathcal D_K (x^{(\epsilon_k+\epsilon_{-k})})$,
$e:=2\mathcal D_K (x^{(2\epsilon_k+\epsilon_{-k})})$ $(1\leq k\leq n)$,
there is a noncommutative and noncocommutative Hopf algebra structure
$(\mathbf{u}_{t,q}(\mathbf{K}(2n{+}1;\underline{1})),m,\iota,\Delta,S,\varepsilon)$
over $\mathcal{K}[t]_p^{(q)}$ with its algebra structure undeformed,
whose coalgebra structure is given by
\begin{gather}
\Delta(\mathcal D_K (x^{(\al)}))=\mathcal D_K (x^{(\al)}) \otimes (1{-}et)^{\al_k-\al_{-k}}+\sum\limits_{\ell=0}^{p-1} (-1)^\ell h^{\lg \ell \rg} \otimes (1{-}et)^{-\ell} d^{(\ell)}\Big(\mathcal D_K (x^{(\al)})\Big) t^\ell,\\
S(\mathcal D_K (x^{(\al)}))=-(1{-}et)^{\al_{-k}-\al_k}\sum\limits_{\ell=0}^{p-1} d^{(\ell)}\Big(\mathcal D_K (x^{(\al)})\Big) h_{1}^{\lg \ell \rg} t^\ell,\hskip4cm\\
\varepsilon(\mathcal D_K (x^{(\al)}))=0,\hskip9cm
\end{gather}
for $0 \leq \al \leq \tau$, which is finite dimensional with $\dim_{\mathcal{K}}\mathbf{u}_{t,q}(\mathbf{K}(2n{+}1;\underline{1}))=p^{p^{2n{+}1}+1}$ if $2n+4 \nequiv 0 \;(\mathrm{mod}\, p)$, and  for $0 \leq \al < \tau$, which is a finite dimensional with $\dim_{\mathcal{K}}\mathbf{u}_{t,q}(\mathbf{K}(2n{+}1;\underline{1}))=p^{p^{2n+1}}$ if $2n+4 \equiv 0\; (\mathrm{mod}\, p)$.
\end{theorem}

\begin{proof}
Set $U_{t,q}(\mathbf K(2n{+}1;\underline 1)):=U(\mathbf K(2n{+}1;\underline
1))\otimes_{\mathcal K}\mathcal K[t]_p^{(q)}$. Note the result
of Theorem \ref{vHopf3}, via the base change with $\mathcal K[t]$
instead of $\mathcal K[t]_p^{(q)}$, is still valid over
$U_{t,q}(\mathbf K(2n{+}1;\underline 1))$.
 Denote by $I_{t,q}$ the ideal of $U_{t,q}(\mathbf K(2n{+}1;\underline
1))$ over the ring $\mathcal K[t]_p^{(q)}$ generated by the same
generators of the ideal $I$ in $U(\mathbf K(2n{+}1;\underline 1))$ via
the base change with $\mathcal K$ replaced by $\mathcal
K[t]_p^{(q)}$. We shall show that $I_{t,q}$ is a Hopf ideal of
$U_{t,q}(\mathbf K(2n{+}1;\underline 1))$. It suffices to verify that
$\Delta$ and $S$ preserve the generators in $I_{t,q}$ of
$U_{t,q}(\mathbf K(2n{+}1;\underline 1))$.

(I) By \ref{vrelation11}, Lemmas \ref{relation}, and \ref{vrelation13}, we have
\begin{equation}\label{vHopf4d}
\begin{split}
&\Delta((\mathcal D_K (x^{(\al)}))^p)=\sum\limits_{0 \leq j \leq p \atop \ell \geq 0}(-1)^\ell \binom{p}{j}(\mathcal D_K (x^{(\al)}))^j h^{\lg \ell \rg} \otimes (1{-}et)^{j(\al_k-\al_{-k})-\ell} d^{(\ell)}  (\mathcal D_K (x^{(\al)}))^{p-j} t^\ell\\
&\equiv \sum\limits_{\ell=0}^{p-1} (-1)^\ell h^{\lg \ell \rg} \otimes (1{-}et)^{-\ell}  d^{(\ell)} \cdot (\mathcal D_K (x^{(\al)}))^{p} t^\ell+\\
&\quad \sum\limits_{\ell=0}^{p-1} (-1)^\ell (\mathcal D_K (x^{(\al)}))^p h^{\lg \ell \rg} \otimes (1{-}et)^{p(\al_k-\al_{-k})-\ell} d^{(\ell)} \cdot 1 t^\ell\\
&\equiv 1 \otimes (\mathcal D_K (x^{(\al)}))^p+(-1)h \otimes (1{-}et)^{-1} (\delta_{\al,\ep_k+\ep_{-k}}+\delta_{\al,\ep_0})(-e)t+(\mathcal D_K (x^{(\al)}))^p \otimes 1\\
&=1 \otimes (\mathcal D_K (x^{(\al)}))^p+(\mathcal D_K (x^{(\al)}))^p \otimes 1\\
&\quad +h \otimes (1{-}et)^{-1} \delta_{\al,\ep_i+\ep_{-i}}\delta_{ik}et+h \otimes (1{-}et)^{-1} \delta_{\al,\ep_0}et.
\end{split}
\end{equation}
If $\al \neq \ep_i+\ep_{-i}$, $\ep_0$ for $1 \leq i \leq n$, we get
\begin{equation*}
\begin{split}
\Delta( (\mathcal D_K (x^{(\al)}))^p)&= (\mathcal D_K (x^{(\al)}))^p \otimes 1+ 1 \otimes  (\mathcal D_K (x^{(\al)}))^p\\
&\in I_{t,q} \otimes U_{t,q}({\bf K}(2n{+}1;\underline{1}))+ U_{t,q}({\bf K}(2n{+}1;\underline{1})) \otimes I_{t,q}.
\end{split}
\end{equation*}
When $\al=\ep_i+\ep_{-i}$, by Lemma  \ref{vrelation13}, we have
$$
\Delta(\mathcal D_K (x^{(\ep_i+\ep_{-i})}))=\mathcal D_K (x^{(\ep_i+\ep_{-i})}) \otimes 1+1 \otimes \mathcal D_K (x^{(\ep_i+\ep_{-i})})+h \otimes (1{-}et)^{-1} \delta_{ik}et.
$$
Combing this with \ref{vHopf4d}, we have
\begin{equation*}
\begin{split}
\Delta\big((\mathcal D_K (x^{(\ep_i+\ep_{-i})}))^p-\mathcal D_K (x^{(\ep_i+\ep_{-i})})\big)
&= \big((\mathcal D_K (x^{(\ep_i+\ep_{-i})}))^p-\mathcal D_K (x^{(\ep_i+\ep_{-i})})\big) \otimes 1\\
&\quad+1 \otimes \big((\mathcal D_K (x^{(\ep_i+\ep_{-i})}))^p-\mathcal D_K (x^{(\ep_i+\ep_{-i})})\big)\\
&\in I_{t,q} \otimes U_{t,q}({\bf K}(2n{+}1;\underline{1}))+ U_{t,q}({\bf K}(2n{+}1;\underline{1})) \otimes I_{t,q}.
\end{split}
\end{equation*}
When $\al=\ep_0$, by Lemma  \ref{vrelation13}, we obtain
$$
\Delta(\mathcal D_K (x^{(\ep_0)}))=\mathcal D_K (x^{(\ep_0)}) \otimes 1+ 1 \otimes \mathcal D_K (x^{(\ep_0)})+h \otimes (1{-}et)^{-1}et.
$$
Thus, we have
\begin{equation*}
\begin{split}
\Delta\big(\mathcal D_K (x^{(\ep_0)})^p-\mathcal D_K (x^{(\ep_0)})\big)
&=\big(\mathcal D_K (x^{(\ep_0)})^p-\mathcal D_K (x^{(\ep_0)}) \big) \otimes 1\\
&\quad+1 \otimes \big(\mathcal D_K (x^{(\ep_0)})^p-\mathcal D_K (x^{(\ep_0)})\big)\\
&  \in I_{t,q} \otimes U_{t,q}({\bf K}(2n{+}1;\underline{1}))+ U_{t,q}({\bf K}(2n{+}1;\underline{1})) \otimes I_{t,q}.
\end{split}
\end{equation*}

So the ideal $I_{t,q}$ is a coideal of the Hopf algebra $U_{t,q}({\bf K}(2n{+}1;\underline{1}))$.

(II) By \ref{vrelation12}, Lemmas \ref{relation} \& \ref{vrelation13}, we have
\begin{equation*}
\begin{split}
S(\mathcal D_K (x^{(\al)})^p)&=(-1)^p (1{-}et)^{-p(\al_k-\al_{-k})}\sum\limits_{\ell=0}^{p-1} d^{(\ell)}\mathcal D_K (x^{(\al)})^p) h_{1}^{\lg \ell \rg} t^\ell\\
&=-\big(\mathcal D_K (x^{(\al)})^p-\delta_{\al,\ep_k+\ep_{-k}}e h_{1}^{\lg 1 \rg}t-\delta_{\al,\ep_0}eh_{1}^{\lg 1 \rg}t\big)\\
&=-\mathcal D_K (x^{(\al)})^p+\delta_{\al,\ep_i+\ep_{-i}}\delta_{ik}eh_{1}^{\lg 1 \rg}t+\delta_{\al,\ep_0}eh_{1}^{\lg 1 \rg} t .
\end{split}
\end{equation*}

When $\al \neq \ep_i+\ep_{-i}$, $\al \neq \ep_0$, we have $S(\mathcal D_K (x^{(\al)})^p)=-S(\mathcal D_K (x^{(\al)})^p) \in I_{t,q}$. If $\al=\ep_i+\ep_{-i}$, by Lemma \ref{vrelation13},
$S(\mathcal D_K (x^{(\ep_i+\ep_{-i})}))=-\big(\mathcal D_K (x^{(\ep_i+\ep_{-i})})-\delta_{ik}eh_{1}^{\lg 1 \rg} t\big)$.
Hence,
\begin{equation*}
\begin{split}
S\big(\mathcal D_K (x^{(\ep_i+\ep_{-i})})^p-\mathcal D_K (x^{(\ep_i+\ep_{-i})})\big)&=-\mathcal D_K (x^{(\ep_i+\ep_{-i})})^p+\delta_{ik}eh_{1}^{\lg 1 \rg} t+\big(\mathcal D_K (x^{(\ep_i+\ep_{-i})})-\delta_{ik}eh_{1}^{\lg 1 \rg} t\big)\\
&=-\big(\mathcal D_K (x^{(\ep_i+\ep_{-i})})^p-\mathcal D_K (x^{(\ep_i+\ep_{-i})})\big) \in I_{t,q}.
\end{split}
\end{equation*}

Using a similar argument, we have
\begin{equation*}
\begin{split}
S\big((\mathcal D_K (x^{(\ep_0)}))^p-\mathcal D_K (x^{(\ep_0)})\big)&=-(\mathcal D_K (x^{(\ep_0)}))^p+eh_{1}^{\lg 1 \rg}t+\mathcal D_K (x^{(\ep_0)})-eh_{1}^{\lg 1 \rg}t\\
&=-\big((\mathcal D_K (x^{(\ep_0)}))^p-\mathcal D_K (x^{(\ep_0)})\big) \in I_{t,q}.
\end{split}
\end{equation*}

Thus, the ideal $I_{t,q}$ is preserved by the antipode $S$ of the quantization $U_{t,q}({\bf K}(2n{+}1;\underline{1}))$.

(III) It is obvious that $\ep((\mathcal D_K (\al))^p)=0$ for all $\al$.

Therefore, by (I), (II), (III), we proved that $I_{t,q}$ is a Hopf ideal of $U_{t,q}({\bf K}(2n{+}1;\underline{1}))$. Thus we obtain the required
 $t$-deformation on
 $\mathbf{u}_{t,q}(\mathbf{K}(2n{+}1;\underline{1}))$, for the Cartan type
 simple modular restricted Lie algebra of $\mathbf{K}$ type
 --- the Contact algebra $\mathbf{K}(2n{+}1;\underline{1})$.
\end{proof}

\subsection{More quantizations}

We consider the modular reduction process for the quantizations of
$U(\mathbf{K}^+)[[t]]$ arising from those products of some pairwise
different and mutually commutative basic Drinfel'd twists. We will
then get lots of new families of noncommutative and noncocommutative
Hopf algebras of dimension $p^{p^{2n+1}+1}$ (resp. $p^{p^{2n+1}})$ with
indeterminate $t$ or of dimension $p^{p^{2n+1}}$ (resp. $p^{p^{2n+1}-1})$ with
specializing $t$ into a scalar in $\mathcal{K}$ if $2n+4 \nequiv
0\; (\mod p)$ (resp. if $2n+4 \equiv 0\; (\mod p)$).

Let $A(k)_j$ and $A(k')_{j'}$ denote the coefficients of the
corresponding quantizations of $U(\mathbf{K}^+_{\mathbb{Z}})$ over
$U(\mathbf{K}^+_{\mathbb{Z}})[[t]]$ given by the Drinfel'd twists
$\mathcal{F}(k)$ and $\mathcal{F}(k')$ as in Corollary \ref{vHopf2},
respectively. Note that $A(k)_0=A(k')_0$ $=1$,
$A(k)_{-1}=A(k')_{-1}=0$.

\begin{lemm}\label{mHopf1}
Fix distinguished elements $h(k)=\mathcal D_K ( x^{\ep_k+\ep_{-k}})$, $e(k)=\mathcal D_K ( x^{2\ep_k+\ep_{-k}})$
$(1 \leq k \leq n)$
 and $h(k')=\mathcal D_K (x^{\ep_{k'}+\ep_{-k'}}),~e(k')=\mathcal D_K (x^{2\ep_{k'}+\ep_{-k'}})(1
\leq k' \leq n)$ with $k \neq k'$. the corresponding quantization of
 $U(\mathbf{K}_{\mathbb{Z}}^{+})[[t]]$ by the Drinfel'd twist $\mathcal{F}=\mathcal{F}(k)\mathcal{F}(k')$
with the product undeformed is given by
\begin{eqnarray*}
&&\Delta(\mathcal D_K (x^\al))=\mathcal D_K (x^\al) \otimes
(1{-}e(k)t)^{\al_k-\al_{-k}}(1{-}e(k')t)^{\al_{k'}-\al_{-k'}}\\
&&+ \sum\limits_{n,\ell=0}^{\infty}\sum\limits_{j'=0}^{\ell}  \sum\limits_{j=0}^{n} (-1)^{\ell+n}A(k)_j A(k')_{j'}B(k')_{\ell-j'}C_{n-j}^{\ell,j'} h(k')^{\lg
\ell \rg} h(k)^{\lg n\rg} \otimes (1{-}e(k')t)^{-\ell}(1{-}e(k)t)^{-n}
\\
&&\qquad \times
\mathcal D_K (x^{\al+\ell(2\ep_{k'}+\ep_{-k'})+n(2\ep_k+\ep_{-k})-j'(\ep_{k'}+\ep_{-k'})-j(\ep_k+\ep_{-k})-(\ell+n-j-j')\ep_0})t^{n+\ell},\\
&&S(\mathcal D_K (x^\al))=-(1{-}e(k')t)^{\al_{-k'}-\al_{k'}}(1{-}e(k)t)^{\al_{-k}-\al_{k}}
\sum\limits_{n,\ell=0}^{\infty} \sum\limits_{j'=0}^{\ell} \sum\limits_{j=0}^{n} A(k)_j A(k')_{j'} B(k')_{\ell-j'}C_{n-j}^{\ell,j'}\\
&&\qquad
\mathcal D_K (x^{\al+\ell(2\ep_{k'}+\ep_{-k'})+n(2\ep_k+\ep_k)-j'(\ep_{k'}+\ep_{-k'})-j(\ep_k+\ep_{-k})-(\ell+n-j-j')\ep_0})\times
h(k)_1^{\lg n \rg} h(k')_{1}^{\lg \ell \rg} t^{\ell+n},\\
&&\varepsilon(\mathcal D_K (x^\al))=0
\end{eqnarray*}
where $C_{n-j}^{\ell,j'}=\frac{\prod\limits_{i=0}^{n-j-1}(i-\al_0+\ell-j')}{(n-j)!}$ for $\mathcal D_K  (x^\al)\in {\bf K}^+_{\mathbb{Z}}$.
\end{lemm}
\begin{proof} First of all, let us consider the formula of comultiplication
\begin{equation*}
\begin{split}
&\Delta(\mathcal D_K (x^\al))=\mathcal{F}(k)\mathcal{F}(k')\Delta_0(\mathcal D_K (x^\al))\mathcal{F}(k')^{-1}\mathcal{F}(k)^{-1}\\
&=\mathcal{F}(k) \mathcal D_K (x^\al) \otimes
(1{-}e(k')t)^{\al_{k'}-\al_{-k'}}F(k)\\
&\quad+\mathcal{F}(k)\big(\sum\limits_{\ell=0}^{\infty}(-1)^\ell
h(k')^{\lg \ell \rg} \otimes (1{-}e(k')t)^{-\ell} d_{k'}^{(\ell)}
\mathcal D_K (x^{\al}) t^\ell \big)F(k).\\
\end{split}
\end{equation*}

By Corollary \ref{vHopf2}, we can get
\begin{equation*}
\begin{split}
\mathcal{F}(k)\big(\mathcal D_K (x^\al)& \otimes
\big(1{-}e(k')t\big)^{\al_{k'}-\al_{-k'}}\big)F(k)\\
&=\mathcal{F}(k)(\mathcal D_K (x^\al) \otimes 1)(1 \otimes
(1{-}e(k')t)^{\al_{k'}-\al_{\al_{-k'}}})F(k)\\
&=\mathcal{F}(k)(\mathcal D_K (x^\al) \otimes 1) F(k)( 1 \otimes
(1{-}e(k')t)^{\al_{k'}-\al_{-k'}})\\
&=\mathcal{F}(k)F(k)_{\al_{-k}-\al_k} (\mathcal D_K (x^\al) \otimes 1)(1
\otimes (1{-}e(k')t)^{\al_{k'}-\al_{-k'}})\\
&=\mathcal D_K (x^\al) \otimes
(1{-}e(k)t)^{\al_k-\al_{-k}}(1{-}e(k')t)^{\al_{k'}-\al_{-k'}},
\end{split}
\end{equation*}
\begin{equation*}
\begin{split}
\mathcal{F}(k)&\Big(\sum\limits_{\ell=0}^{\infty}(-1)^\ell h(k')^{\lg
\ell \rg} \otimes (1{-}e(k')t)^{-\ell} d_{k'}^{(\ell)}(\mathcal D_K (x^{\al}))
t^\ell \Big)F(k)\\
&=\mathcal{F}(k)\sum\limits_{\ell=0}^{\infty} (-1)^\ell h(k')^{\lg
\ell \rg} \otimes (1{-}e(k')t)^{-\ell}\times\\
&\quad \times \sum\limits_{j'=0}^{\ell}
A(k')_{j'}
B(k')_{\ell-j'}\mathcal D_K (x^{\al+\ell(2\ep_{k'}+\ep_{-k'})-j'(\ep_{k'}+\ep_{-k'})-(\ell-j')\ep_0})
t^\ell F(k) \\
&=\Big(\sum\limits_{\ell=0}^{\infty} (-1)^\ell h(k')^{\lg \ell \rg}
\otimes (1{-}e(k')t) ^{-\ell}\Big)\times\\
&\quad \times \Big( \sum\limits_{j'=0}^{\ell}
A(k')_{j'} B(k')_{\ell-j'} \mathcal{F}(k)\big(1 \otimes
\mathcal D_K (x^{\al+\ell(2\ep_{k'}+\ep_{k'})-j'(\ep_{k'}+\ep_{-k'})-(\ell-j')\ep_0)})\big)
F(k) t^\ell \Big)\\
&=\sum\limits_{\ell=0}^{\infty} (-1)^\ell h(k')^{\lg \ell \rg}
\otimes (1{-}e(k')t) ^{-\ell} \sum\limits_{j'=0}^{\ell} A(k')_{j'} B(k')_{\ell-j'}
\mathcal{F}(k)\times \\
&\quad \times \sum\limits_{n=0}^{\infty} (-1)^n F(k)_n
\Big(h(k)^{\lg n \rg} \otimes d_{k}^{(n)}(\mathcal D_K (x^{\al+\ell(2\ep_{k'}+\ep_{-k'})-j'(\ep_{k'}+\ep_{-k'})-(\ell-j')\ep_0}))
\Big) t^{n+\ell}.
\end{split}
\end{equation*}
Set $\al(\ell,k',j'):=\al{+}\ell(2\ep_{k'}{+}\ep_{-k'}){-}j'(\ep_{k'}{+}\ep_{-k'}){-}(\ell{-}j')\ep_0$.
It is easy to see
\begin{eqnarray*}
&&d_{k}^{(n)}
\mathcal D_K (x^{\al+\ell(2\ep_{k'}+\ep_{-k'})-j'(\ep_{k'}+\ep_{-k'})-(\ell-j')\ep_0})=d_{k}^{(n)} \mathcal D_K (x^{\al(\ell,k',j')}) \\
&&= \sum\limits_{j=0}^{n}
\frac{\prod\limits_{i=0}^{n-j-1}(i{-}\al(\ell,k',j')_0)}{(n{-}j)!}
\frac{\prod\limits_{i=0}^{j-1}(\al(\ell,k',j')_k{-}2\al(\ell,k',j')_{-k}{+}i)}{j!}\times\\
&&\hskip6.15cm\times
\mathcal D_K (x^{\al(\ell,k',j')+n(2\ep_k+\ep_{-k})-j(\ep_k+\ep_{-k})-(n-j)\ep_0})\\
&&=\sum\limits_{j=0}^{n} \frac{\prod\limits_{i=0}^{n-j-1}
(i{-}\al_0{+}\ell{-}j')}{(n{-}j)!}
\frac{\prod\limits_{i=0}^{j-1}(\al_k{-}2\al_{-k}{+}i)}{j!}\times\\
&&\hskip4cm\times \mathcal D_K (x^{\al+\ell(2\ep_{k'}+\ep_{-k'})
-j'(\ep_{k'}{+}\ep_{-k'})-(\ell-j')\ep_0 +n(2\ep_k+\ep_{-k})-j(\ep_k+\ep_{-k})-(n-j)\ep_0})\\
&&=\sum\limits_{j=0}^{n} C_{n-j}^{\ell,j'} A(k)_j
\mathcal D_K (x^{\al+\ell(2\ep_{k'}+\ep_{-k'})+n(2\ep_k+\ep_{-k})-j'(\ep_{k'}+\ep_{-k'})-j(\ep_k+\ep_{-k})-(\ell+n-j'-j)\ep_0}).
\end{eqnarray*}
On the other hand, we have
\begin{equation*}
\begin{split}
&\mathcal{F}(k)\Bigl(\sum\limits_{\ell=0}^{\infty}(-1)^\ell h(k')^{\lg
\ell \rg} \otimes (1{-}e(k')t)^{-\ell} d_{k'}^{(\ell)} \mathcal D_K (x^{\al})
t^\ell \Bigr)F(k)\\
&=\sum\limits_{\ell=0}^{\infty} (-1)^\ell h(k')^{\lg \ell \rg}
\otimes (1{-}e(k')t)^{-\ell}  \sum\limits_{j'=0}^{\ell}
A(k')_{j'}B(k')_{\ell-j'} \sum\limits_{n=0}^{\infty}(-1)^n (1
\otimes (1{-}e(k)t)^{-n}) \\
&\quad\Big(h(k)^{\lg n \rg} \otimes \sum\limits_{j=0}^{n}
C_{n-j}^{\ell,j'} A(k)_j
\mathcal D_K (x^{\al+\ell(2\ep_{k'}+\ep_{-k'})-j'(\ep_{k'}+\ep_{-k'})-(\ell-j')\ep_0
+n(2\ep_k+\ep_{-k})-j(\ep_k+\ep_{-k})-(n-j)\ep_0})\Big) t^{n+\ell}\\
&=\sum\limits_{n,\ell=0}^{\infty}\sum\limits_{j'=0}^{\ell}\sum\limits_{j=0}^{n} (-1)^{\ell+n} h(k')^{\lg \ell \rg}
h(k)^{\lg n \rg} \otimes (1{-}e(k')t)^{-\ell} (1{-}e(k)t)^{-n}A(k)_j A(k')_{j'} B(k')_{\ell-j'}C_{n-j}^{\ell,j'} \\
&\quad
\cdot
\mathcal D_K (x^{\al+\ell(2\ep_{k'}+\ep_{-k'})+n(2\ep_k+\ep_{-k})-j'(\ep_{k'}+\ep_{-k'})-j(\ep_k+\ep_{-k})-(\ell+n-j'-j)\ep_0})
t^{n+\ell}.
\end{split}
\end{equation*}
So we get the required formula for the comultiplication.

Next, we consider the second formula
\begin{equation*}
\begin{split}
&S(\mathcal D_K (x^\al))=-v(k)v(k') \mathcal D_K (x^\al)u(k') u(k)\\
&=-v(k) (1{-}e(k')t)^{\al_{-k'}-\al_{k'}}\sum\limits_{\ell=0}^{\infty}
\sum\limits_{j=0}^{\ell} A(k')_{j'} B(k')_{\ell-j'}
\mathcal D_K (x^{\al+\ell(2\ep_{k'}+\ep_{k'})-j'(\ep_{k'}+\ep_{-k'})-(\ell-j')\ep_0})
h(k')_1^{\lg \ell \rg} u(k) t^\ell\\
&=- (1{-}e(k')t)^{\al_{-k'}-\al_{k'}}\sum\limits_{\ell=0}^{\infty}
\sum\limits_{j'=0}^{\ell} A(k')_{j'}
B(k')_{\ell-j'}v(k)u(k)_{\al_k-\al_{-k}}\times\\
&\quad \times \sum\limits_{n=0}^{\infty} d_k^{(n)}
\mathcal D_K (x^{\al+\ell(2\ep_{k'}+\ep_{k'})-j'(\ep_{k'}+\ep_{-k'})-(\ell-j')\ep_0})
h(k)_1^{\lg n \rg} t^n h(k')_1^{\lg \ell \rg} t^\ell\\
&=-(1{-}e(k')t)^{\al_{-k'}-\al_{k'}}
(1{-}e(k)t)^{\al_{-k}-\al_{k}}\sum\limits_{\ell=0}^{\infty}
\sum\limits_{j'=0}^{\ell} A(k')_{j'}B(k')_{\ell-j'} \sum\limits_{n=0}^{\infty} \sum\limits_{j=0}^{n}
C_{n-j}^{\ell,j'}A(k)_j \times \\
&\quad \times
\mathcal D_K (x^{\al+\ell(2\ep_{k'}+\ep_{-k'})-j'(\ep_{k'}+\ep_{-k'})-(\ell-j')\ep_0
+n(2\ep_k+\ep_{-k})-j(\ep_k+\ep_{-k})-(n-j)\ep_0}) h(k)_1^{\lg n \rg}
 h(k')_1^{\lg \ell \rg}) t^{n+\ell}\\
\end{split}
\end{equation*}
\begin{equation*}
\begin{split}
&=-(1{-}e(k')t)^{\al_{-k'}-\al_{k'}}
(1{-}e(k)t)^{\al_{-k}-\al_{k}}\sum\limits_{n,\ell=0}^{\infty}\sum\limits_{j'=0}^{\ell}\sum\limits_{j=0}^{n}
 A(k)_j A(k')_{j'} B(k')_{\ell-j'} C_{n-j}^{\ell,j'}\times \\
&\quad \times
\mathcal D_K (x^{\al+\ell(2\ep_{k'}+\ep_{-k'})+n(2\ep_k+\ep_{-k})-j'(\ep_{k'}+\ep_{-k'})-j(\ep_k+\ep_{-k})-(\ell+n-j'-j)\ep_0})
h(k)_1^{\lg n \rg} h(k')_1^{\lg \ell \rg} t^{n+\ell}.
\end{split}
\end{equation*}

This completes the proof.
\end{proof}

\begin{lemm}\label{mHopf2}
Fix distinguished elements $h(k)=\mathcal D_K (x^{(\ep_k+\ep_{-k})}),~e(k)=2\mathcal D_K (x^{(2\ep_k+\ep_{-k})})$,
$h(k')=\mathcal D_K (x^{(\ep_{k'}+\ep_{-k'})}),~e(k')=2\mathcal D_K (x^{(2\ep_{k'}+\ep_{-k'})})$, with $1 \leq k
\neq k' \leq n$, then the associated quantization of
$U(\mathbf{K}(2n{+}1;\underline{1}))$ on $U_t(\mathbf{K}(2n{+}1;\underline{1}))$(also
on $U(\mathbf{K}(2n{+}1;\underline{1}))[[t]]$) with the product
undeformed is given by
\begin{equation*}
\begin{split}
&\Delta(\mathcal D_K (x^{(\al)}))=\mathcal D_K (x^{(\al)}) \otimes
(1{-}e(k)t)^{\al_k-\al_{-k}}(1{-}e(k')t)^{\al_{k'}-\al_{-k'}}\\
&+\sum\limits_{n,\ell=0}^{p-1}(-1)^{\ell+n} \sum\limits_{j'=0}^{\ell}
\sum\limits_{j=0}^{n} \overline{A(k',\ell)}_{j'} \overline{A(k,n)}_j
\overline{B(k,k')}_{n,j}^{\ell,j'}h(k')^{\lg
\ell \rg} h(k)^{\lg n\rg}\otimes (1{-}e(k')t)^{-\ell}(1{-}e(k)t)^{-n} \\
&\qquad \qquad  \qquad \times\mathcal D_K (x^{\big(\al+\ell(2\ep_{k'}+\ep_{-k'})+n(2\ep_k+\ep_{-k})-j'(\ep_{k'}+\ep_{-k'})-j(\ep_k+\ep_{-k})-(\ell+n-j-j')\ep_0\big)}) t^{n+\ell},\\
&S(\mathcal D_K (x^{(\al)}))=-(1{-}e(k')t)^{\al_{-k'}-\al_{k'}}(1{-}e(k)t)^{\al_{-k}-\al_{k}}\sum\limits_{n,\ell=0}^{p-1}
\sum\limits_{j'=0}^{\ell}
\sum\limits_{j=0}^{n} \overline{A(k',\ell)}_{j'} \overline{A(k,n)}_j
\overline{B(k,k')}_{n,j}^{\ell,j'}\\
&\qquad \qquad \qquad\times
\mathcal D_K (x^{\big(\al+\ell(2\ep_{k'}+\ep_{-k'})+n(2\ep_k+\ep_k)-j'(\ep_{k'}+\ep_{-k'})-j(\ep_k+\ep_{-k})-(\ell+n-j-j')\ep_0\big)})
h(k)_1^{\lg n \rg} h(k')_{1}^{\lg \ell \rg} t^{\ell+n},\\
&\varepsilon(\mathcal D_K (x^{(\al)}))=0,
\end{split}
\end{equation*}
where $0 \leq \al \leq \tau$ for $2n{+}4 \nequiv
0 \;(\mathrm{mod} p)$, and $0 \leq \al < \tau$ for $2n{+}4 \equiv 0\ (\mathrm{mod}
p)$. $\overline{A(k',\ell)}_{j'}=(2\ell{-}j')!\binom{\al_{k'}{+}2\ell{-}j'}{2\ell{-}j'},
~\overline{A(k,n)}_j=(2n{-}j)!\binom{\al_k{+}2n{-}j}{2n{-}j},~
\overline{B(k,k')}_{n,j}^{\ell,j'}=(-1)^{\ell+n-j-j'}\binom{\al_k+\ell-j}{\ell-j}
\binom{\al_{-k'}{+}\ell{-}j'}{\ell{-}j'}$ for $n{+}\ell{-}j'{-}j
\leq \al_0$, otherwise $\overline{B(k,k')}_{n,j}^{\ell,j'}=0$.
\end{lemm}
\begin{proof}
Note that the elements $\frac{1}{\al !} \mathcal D_K
(x^{\al})$ in $\mathbf{K}_{\mathcal{K}}^{+}$ will be identified with
 $\mathcal D_K ({x^{(\al)}})$ in $\mathbf{K}(2n{+}1;1)$ and those
in $J_{\underline{1}}$ with 0. Hence, by Lemma \ref{mHopf1}, we can get
\begin{equation*}
\begin{split}
&\Delta(\mathcal D_K (x^{(\al)}))=\frac{1}{\al
!}\Delta(\mathcal D_K (x^\al))=\mathcal D_K (x^{(\al)}) \otimes
(1{-}e(k)t)^{\al_k-\al_{-k}}(1{-}e(k')t)^{\al_{k'}-\al_{-k'}}\\
&\quad+ \sum\limits_{n,\ell=0}^{p-1} h(k')^{\lg \ell \rg} h(k)^{\lg
n \rg} \otimes
(1{-}e(k')t)^{-\ell}(1{-}e(k)t)^{-n}\sum\limits_{j'=0}^{\ell}\sum\limits_{j=0}^{n}
A(k')_{j'}B(k')_{\ell-j'} C(k)_{n-j}^{\ell,j'}A(k)_j\\
&\qquad \cdot
\frac{(\al{+}\ell(2\ep_{k'}{+}\ep_{-k'}){+}n(2\ep_k{+}\ep_{-k}){-}j'(\ep_{k'}{+}\ep_{-k'}){-}j(\ep_k{+}\ep_{-k}){-}(n{+}\ell{-}j{-}j')\ep_0)!}{\al
!}\\
&\qquad \cdot
\mathcal D_K (x^{(\al+\ell(2\ep_{k'}+\ep_{-k'})+n(2\ep_k+\ep_{-k})-j'(\ep_{k'}+\ep_{-k'})-j(\ep_k+\ep_{-k})-(n+\ell-j-j')\ep_0)})t^{n+\ell}.
\end{split}
\end{equation*}
We can see that
{\setlength{\arraycolsep}{0pt}
\begin{eqnarray*}
&&A(k')_{j'}B(k')_{\ell-j'} C_{n-j}^{\ell,j'}A(k)_j
\frac{(\al{+}\ell(2\ep_{k'}{+}\ep_{-k'}){+}n(2\ep_k{+}\ep_{-k}){-}j'(\ep_{k'}{+}\ep_{-k'}){-}j(\ep_k{+}\ep_{-k}){-}(n{+}\ell{-}j{-}j')\ep_0)!}{\al
!}\\
&&=A(k')_{j'}\frac{\prod\limits_{i=0}^{\ell-j'-1}(i{-}\al_0)}{(\ell{-}j')!}
\frac{\prod\limits_{i=0}^{n-j-1}(i{-}\al_0{+}\ell{-}j')}{(n{-}j)!}A(k)_j\cdot \\
&&\quad
\cdot\frac{(\al_{k'}{+}2\ell{-}j')!(\al_{-k'}{+}\ell{-}j')!(\al_k{+}2n{-}j)!(\al_{-k}{+}n{-}j)!(\al_0{-}(\ell{+}n{-}j{-}j'))!}{\al_{k'}!\al_{-k'}!\al_k!\al_{-k}!\al_0!}\\
&&=A(k')_{j'}\frac{(\al_{k'}{+}2\ell{-}j')!}{\al_{k'}!}A(k)_j
\frac{(\al_k{+}2n{-}j)!}{\al_k !}\frac{\prod\limits_{i=0}^{\ell-j'-1}(i{-}\al_0)}{(\ell{-}j')!}
\frac{\prod\limits_{i=0}^{n-j-1}(i{-}\al_0{+}\ell{-}j')}{(n{-}j)!}\cdot\\
&&\qquad\cdot\frac{(\al_{-k'}{+}\ell{-}j')!}{\al_{-k'}!}\frac{(\al_0{-}(\ell{+}n{-}j{-}j'))!}{\al_0
!}\frac{(\al_{-k}{+}n{-}j)!}{\al_{-k}!} \\
&&=A(k')_{j'}(2\ell{-}j')!\binom{\al_{k'}{+}2\ell{-}j'}{2\ell{-}j'} A(k)_j
(2n{-}j)! \binom{\al_k{+}2n{-}j}{2n{-}j}\binom{\al_{-k'}{+}\ell{-}j'}{\ell{-}j'}\binom{\al_{-k}{+}n{-}j}{n{-}j}\cdot \\
&& \cdot\frac{(-1)^{\ell-j'}\al_0 \cdots(\al_0{-}(\ell{-}j'){+}1)
(-1)^{n-j}(\al_0{-}(\ell{-}j')) \cdots (\al_0{-}(\ell{+}n{-}j{-}j'))!}{\al_0 !}\\
&&=\overline{A(k',\ell)}_{j'} ~\overline{A(k,n)}_j
\binom{\al_{-k'}{+}\ell{-}j'}{\ell{-}j'}\binom{\al_{-k}{+}n{-}j}{n{-}j}B,
\end{eqnarray*}
}
where
$$
B=\frac{(-1)^{\ell+n-j-j'}\al_0
\cdots(\al_0{-}(\ell{-}j'){+}1)(\al_0{-}(\ell{-}j')) \cdots
(\al_0{-}(\ell{+}n{-}j{-}j'))!}{\al_0 !}
=(-1)^{\ell+n-j-j'},
$$
when $\ell{-}j'\leq \al_0,$ and $n{-}j \leq \al_0{-}(\ell{-}j')$. In fact, it is easy to see  $0 \leq n{-}j \leq \al_0{-}(\ell{-}j')$ implies $ \ell{-}j'\leq
\al_0$. So we have $ B=(-1)^{\ell+n-j-j'}$, if $n{-}j \leq
\al_0{-}(\ell{-}j')$, i.e. $ n{+}\ell{-}j{-}j' \leq \al_0$. Otherwise, $B=0$.
\begin{equation*}
\begin{split}
&A(k')_{j'}B(k')_{\ell-j'} C_{n-j}^{\ell,j'}A(k)_j\times\\
&\quad \times
\frac{(\al{+}\ell(2\ep_{k'}{+}\ep_{-k'}){+}n(2\ep_k{+}\ep_{-k}){-}j'(\ep_{k'}{+}\ep_{-k'}){-}j(\ep_k{+}\ep_{-k}){-}(n{+}\ell{-}j{-}j')\ep_0)!}{\al
!}\\
&=\overline{A(k',\ell)}_{j'} ~\overline{A(k,n)}_j
\binom{\al_{-k'}{+}\ell{-}j'}{\ell{-}j'}\binom{\al_{-k}{+}n{-}j}{n{-}j} B\\
&=\overline{A(k',\ell)}_{j'} ~\overline{A(k,n)}_j
\overline{B(k,k')}_{n,j}^{\ell,j'}.
\end{split}
\end{equation*}
Hence, we get the first formula.

Applying a similar argument, we can get the other formulas.

This completes the proof.
\end{proof}

\begin{lemm}
For $s\geq 1$, one has
\begin{equation*}\label{mrelation1}
\begin{split}
&\Delta(\mathcal D_K (x^{(\al)})^s)=\sum_{0\le j'\le s\atop n, \ell\ge
0}\binom{s}{j'}(-1)^{n+\ell}\mathcal D_K (x^{(\al)})^{j'} h(k')^{\lg \ell \rg}
h(k)^{\lg n \rg} \otimes\\
&\qquad\qquad(1{-}e(k)t)^{j'(\al_k-\al_{-k})-n}
(1{-}e(k')t)^{j(\al_{k'}-\al_{-k'})-\ell}
d_k^{(n)}\big(d_{k'}^{(\ell)} \mathcal D_K (x^{(\al)})^{s-j'}\big) t^{n+\ell},\\
&S(\mathcal D_K (x^{(\al)})^s)=(-1)^s
(1{-}e(k')t)^{-s(\al_{k'}-\al_{-k'})}(1{-}e(k)t)^{-s(\al_k-\al_{-k})}\times\\
&\qquad\qquad\qquad\quad \times \sum\limits_{n,\ell=0}^{\infty}
d_k^{(n)}\big(d_{k'}^{(\ell)} \mathcal D_K (x^{(\al)})^{s}\big) h(k')_1^{\lg
\ell \rg} h(k)_1^{\lg n \rg} t^{n+\ell}.
\end{split}
\end{equation*}
\end{lemm}
\begin{proof} Let us calculate the following
\begin{eqnarray*}
&&\Delta(\mathcal D_K (x^{(\al)})^s)=\mathcal{F}(k)(\mathcal D_K (x^{(\al)}) \otimes 1 +
1 \otimes \mathcal D_K (x^{(\al)})^s \mathcal{F}^{-1}(k)\\
&&=\mathcal{F}(k) \sum_{0\le j'\le s\atop  \ell\ge 0}\binom{s}{j'}
({-}1)^{\ell}\mathcal D_K (x^{(\al)})^{j'} h(k')^{\lg \ell \rg} \otimes
(1{-}e(k')t)^{j'(\al_{k'}-\al_{-k'})-\ell} d_{k'}^{(\ell)}
 \mathcal D_K (x^{(\al)})^{s-j'} t^\ell F(k)\\
&&=\mathcal{F}(k) \sum_{0\le j'\le s\atop  \ell\ge 0}\binom{s}{j'}
({-}1)^{\ell}\Big(\mathcal D_K (x^{(\al)})^{j'} \otimes 1\Big)\Big( h(k')^{\lg \ell
\rg} \otimes (1{-}e(k')t)^{j'(\al_{k'}-\al_{-k'})-\ell}\Big)\cdot\\
&&\qquad\cdot\Big(1
\otimes d_{k'}^{(\ell)}
 \mathcal D_K (x^{(\al)})^{s-j'}\Big)F(k)t^\ell\\
&&=\mathcal{F}(k) \sum_{0\le j'\le s\atop  \ell\ge 0}\binom{s}{j'}
({-}1)^{\ell}\Big(\mathcal D_K (x^{(\al)})^{j'} \otimes 1\Big)\Big( h(k')^{\lg \ell
\rg} \otimes (1{-}e(k')t)^{j'(\al_{k'}-\al_{-k'})-\ell}\Big) \cdot \\
&&\qquad\sum\limits_{n=0}^{\infty}(-1)^n F(k)_n h(k)^{\lg n
\rg} \otimes d_k^{(n)} \big(d_{k'}^{(\ell)}
\mathcal D_K (x^{(\al)})^{s-j'}\big)t^n t^\ell\\
&&=\mathcal{F}(k) \sum_{0\le j'\le s\atop  \ell\ge 0}\binom{s}{j'}
({-}1)^{\ell}\Big(\mathcal D_K (x^{(\al)})^{j'} \otimes 1\Big) F(k)_n\\
&&\qquad\sum\limits_{n=0}^{\infty}(-1)^n  h(k')^{\lg \ell \rg} h(k)^{\lg n
\rg}\otimes (1{-}e(k')t)^{j'(\al_{k'}-\al_{-k'})-\ell} d_k^{(n)}
\big(d_{k'}^{(\ell)} \mathcal D_K (x^{(\al)})^{s-j'}\big)t^{n+\ell}\\
&&= \sum_{0\le j'\le s\atop  \ell\ge 0}\binom{s}{j'}
({-}1)^{\ell}\mathcal{F}(k)
F(k)_{n+j'(\al_{-k}-\al_k)}\Big(\mathcal D_K (x^{(\al)})^{j'} \otimes 1
\Big)\\
&&\qquad\sum\limits_{n=0}^{\infty}(-1)^n  h(k')^{\lg \ell \rg}
h(k)^{\lg n \rg}\otimes (1{-}e(k')t)^{j'(\al_{k'}-\al_{-k'})-\ell}
d_k^{(n)} \big(d_{k'}^{(\ell)} \mathcal D_K (x^{(\al)})^{s-j'}\big)t^{n+\ell}\\
&&=\sum_{0\le j'\le s\atop n, \ell\ge
0}\binom{s}{j'}(-1)^{n+\ell}\mathcal D_K (x^{(\al)})^{j'} h(k')^{\lg \ell \rg}
h(k)^{\lg n \rg} \otimes\\
&&\qquad(1{-}e(k)t)^{j'(\al_k-\al_{-k})-n}
(1{-}e(k')t)^{j'(\al_{k'}-\al_{-k'})-\ell}
d_k^{(n)}\big(d_{k'}^{(\ell)}\mathcal D_K (x^{(\al)})^{s-j'}\big) t^{n+\ell},
\\
&&S((\mathcal D_K (x^{(\al)}))^s)=v(k) v(k')
S_0(\mathcal D_K (x^{(\al)})^s)u(k')u(k)\\
&&=(-1)^s
v(k)(1{-}e(k')t)^{-s(\al_{k'}-\al_{-k'})}\sum\limits_{\ell=0}^{\infty}
d_{k'}^{(\ell)}(\mathcal D_K (\al^{(\al)})^s)  h(k')_1^{\lg \ell \rg} t^\ell u(k)\\
&&=({-}1)^s v(k)(1{-}e(k')t)^{-s(\al_{k'}-\al_{-k'})}
u(k)_{s(\al_k-\al_{-k})}
\sum\limits_{\ell,n=0}^{\infty}
d_k^{(n)}\big(d_{k'}^{(\ell)}\mathcal D_K (x^{(\al)})^{s}\big) h(k)_1^{\lg n
\rg} h(k')_1^{\lg \ell \rg} t^{n+\ell}\\
&&=({-}1)^s
(1{-}e(k')t)^{-s(\al_{k'}-\al_{-k'})}(1{-}e(k)t)^{-s(\al_k-\al_{-k})}\sum\limits_{n,\ell}
d_k^{(n)}\big(d_{k'}^{(\ell)}\mathcal D_K (x^{(\al)})^{s}\big) h(k')_1^{\lg
\ell \rg} h(k)_1^{\lg n \rg} t^{n+\ell}.
\end{eqnarray*}

This completes the proof.
\end{proof}

To describe $\mathbf{u}_{t,q}(\mathbf{K}(2n{+}1;\underline{1}))$ explicitly, we
still need an auxiliary Lemma.

\begin{lemm}\label{mrelation2}
Set $e(k)=2\mathcal D_K (x^{(2\ep_k+\ep_{-k})}),~e(k')=2\mathcal D_K (x^{(2\ep_{k'}+\ep_{-k'})})$, $d_k^{(n)}=\frac{1}{n!}(\ad
 e(k))^n$, $d_{(k')}^\ell=\frac{1}{\ell!} (\ad e(k'))^\ell$, Then
 \begin{gather*}
d_k^{(n)}
d_{k'}^{(\ell)}(\mathcal D_K (x^{(\al)})=\sum\limits_{j'=0}^{\ell}\sum\limits_{j=0}^{n}\overline{A(k',\ell)}_j
\overline{A(k,n)}_j \overline{B(k,k')}_{n,j}^{\ell,j'}\times\qquad\qquad\qquad\tag{\text{\rm i}}\\
 \qquad\qquad\qquad\qquad\qquad\quad \times \mathcal D_K \big(x^{(\al+\ell(2\ep_{k'}+\ep_{-k'})+n(2\ep_k+\ep_{-k})-j'(\ep_{k'}+\ep_{-k'})
 -j(\ep_k+\ep_{-k})-(\ell+n-j-j')\ep_0)}\big).\nonumber\\
 d_k^{(n)} d_{k'}^{(\ell)}(\mathcal D_K (x^{(\ep_i+\ep_{-i})})=\delta_{\ell
0}\delta_{n 0}\mathcal D_K (x^{(\ep_i+\ep_{-i})})- \delta_{\ell 1}\delta_{n
0}\delta_{ik'}e(k') -\delta_{\ell 0}
\delta_{n 1} \delta_{i k} e(k),\;\tag{\text{\rm ii}} \\
 d_k^{(n)} d_{k'}^{(\ell)}(\mathcal D_K (x^{(\ep_0)}))=\delta_{\ell
0}\delta_{n 0}\mathcal D_K (x^{(\ep_0)})- \delta_{\ell 0}\delta_{n 1}e(k)
-\delta_{\ell 1}
\delta_{n 0}  e(k').\quad\qquad\\
  d_k^{(n)} d_{k'}^{(\ell)}((\mathcal D_K (x^{(\al)}))^p)=\delta_{\ell
 0}\delta_{n0} \mathcal D_K ((x^{(\al)})^p)-\delta_{\ell
 0}\delta_{n1}(\delta_{\al,\ep_k+\ep_{-k}}+\delta_{\al,\ep_0})e(k)\qquad\tag{\text{\rm iii}}\\
\qquad-\delta_{\ell 1}\delta_{n
0}(\delta_{\al,\ep_{k'}+\ep_{-k'}}+\delta_{\al,\ep_0})e(k').
\end{gather*}
\end{lemm}
\begin{proof}
(i) follows from the proof of Lemma \ref{vrelation13}.

(ii) By Lemma \ref{vrelation13} (ii), we have
\begin{equation*}
\begin{split}
d_k^{(n)}d_{k'}^{(\ell)}(\mathcal D_K (x^{(\ep_i+\ep_{-i})})&=d_k^{(n)}
\big( \delta_{\ell 0 }\mathcal D_K (x^{(\ep_i+\ep_{-i})})-\delta_{\ell
1}\delta_{ik'}e(k')\big)\\
&=\delta_{\ell 0}d_{k}^{(n)}(\mathcal D_K  (x^{(\ep_i+\ep_{-i})}))-\delta_{\ell
1}\delta_{ik'}d_{k}^{(n)}\cdot(e(k'))\\
&=\delta_{\ell 0}\delta_{n0}\mathcal D_K (x^{(\ep_i+\ep_{-i})})-\delta_{\ell
0}\delta_{n1}\delta_{ik}e(k)-\delta_{\ell
1}\delta_{n0}\delta_{ik'}e(k'),\\
d_k^{(n)}d_{k'}^{(\ell)}\mathcal D_K (x^{(\ep_0)})&=d_k^{(n)}(\delta_{\ell
0}\mathcal D_K (x^{(\ep_0)})-\delta_{\ell 1}e(k'))\\
&=\delta_{\ell 0}\delta_{n 0}\mathcal D_K (x^{(\ep_0)})-\delta_{\ell 0}\delta_{n
1}e(k)-\delta_{n0}\delta_{\ell 1}e(k').
\end{split}
\end{equation*}
(iii) By (iii) of Lemma  \ref{vrelation13}, we get
\begin{equation*}
\begin{split}
& d_k^{(n)}d_{k'}^{(\ell)} (\mathcal D_K (x^{(\al)}))^p=d_k^{(n)}
\Bigl(\delta_{\ell 0}(\mathcal D_K (x^{(\al)}))^p-\delta_{\ell
1}(\delta_{\al,\ep_{k'}+\ep_{-k'}}{+}\delta_{\al,\ep_0})e(k)\Bigr)\\
&\qquad\qquad=\delta_{\ell 0}\delta_{n 0} (\mathcal D_K (x^{(\al)})^p) -\delta_{\ell
0}\delta_{n
1}\big(\delta_{\al,\ep_{k'}+\ep_{-k'}}{+}\delta_{\al,\ep_0}\big)e(k)-\delta_{n0}\delta_{\ell
1}\big(\delta_{\al,\ep_{k}{+}\ep_{-k}}{+}\delta_{\al,\ep_0}\big)e(k').
\end{split}
\end{equation*}

We complete the proof.
\end{proof}

Using Lemmas \ref{mHopf1}, \ref{mHopf2}, \ref{mrelation1} and \ref{mrelation2}, we get a new Hopf algebra structure over the same
restricted universal enveloping algebra
$\mathbf{u}(\mathbf{K}(2n{+}1;\underline{1}))$ over $\mathcal{K}$ by the
product of two different and commutative basic Drinfel'd twists.

\begin{theorem}\label{mHopf}
Fix distinguished elements $h(k):=\mathcal D_K (x^{(\ep_k+\ep_{-k})}),~e(k)=2\mathcal D_K  (x^{(2\ep_k+\ep_{-k})})$, $h(k')=\mathcal D_K (x^{(\ep_{k'}+\ep_{-k'})}), \;e(k')=2\mathcal D_K (x^{(2\ep_{k'}+\ep_{-k'})})$, with  $1 \leq k \neq k' \leq n$, there is a noncommutative and noncocommutative Hopf algebra $({\bf u}_{t,q}({\bf K})(2n{+}1;\underline{1}),m,\iota,\Delta,S,\varepsilon)$ over $\mathcal{K}[t]_{p}^{(q)}$ with the product undeformed, whose coalgebra structure is given by
\begin{equation*}
\begin{split}
&\Delta(\mathcal D_K (x^{(\al)}))=\mathcal D_K (x^{(\al)}) \otimes (1{-}e(k)t)^{\al_k-\al_{-k}}(1{-}e(k')t)^{\al_{k'}-\al_{-k'}}\\
&\qquad\qquad +\sum\limits_{n,\ell=0}^{p-1} (-1)^{\ell+n}h(k')^{\lg \ell \rg} h(k)^{\lg n \rg} \otimes (1{-}e(k')t)^{-\ell}(1{-}e(k)t)^{-n} d_k^{(n)}d_{k'}^{(\ell)}(\mathcal D_K (x^{(\al)}))t^{n+\ell},
\\
&S(x^{(\al)})=-(1{-}e(k')t)^{\al_{-k'}-\al_{k'}}(1{-}e(k)t)^{\al_{-k}-\al_k}\sum\limits_{n,\ell=0}^{p-1} d_k^{(n)}d_{k'}^{(\ell)} (\mathcal D_K (x^{(\al)}))h(k)_1^{\lg n \rg} h(k')_1^{\lg \ell \rg} t^{n+\ell},\\
&\varepsilon(\mathcal D_K (x^{(\al)}))=0,
\end{split}
\end{equation*}
for $0 \leq \al \leq \tau$, and
$ \dim_{\mathcal{K}}\mathbf{u}_{t,q}(\mathbf{K}(2n{+}1;\underline{1}))=p^{p^{2n+1}+1}$,
if $2n{+}4 \nequiv 0 \;(\mathrm{mod}\, p)$; for $0 \leq \al < \tau$,
$\dim_{\mathcal{K}}\mathbf{u}_{t,q}(\mathbf{K}(2n{+}1;\underline{1}))=p^{p^{2n+1}}$,
if $2n{+}4 \equiv 0\; (\mathrm{mod}\, p)$.
\end{theorem}
\begin{proof}
Let $I_{t,q}$ denote the ideal of $(U_{t,q}({\bf K}(2n{+}1;\underline{1})))$ over the ring
$\mathcal{K}[t]_{p}^{(q)}$ generated by the same generators as in $I~(q \in \mathcal{K})$.
Observe that the result in Lemma \ref{relation}, via the base change with $\mathcal{K}[t]$
replaced by $\mathcal{K}[t]_{p}^{(q)}$, the deformation is still valid for $U_{t,q}({\bf K}(2n{+}1;\underline{1}))$.

In what follows, we shall show that $I_{t,q}$ is a Hopf ideal of $U_{t,q}(K(2n{+}1;\underline{1}))$.
To this end, it suffices to verify that $\Delta$ and $S$ preserve the generators of $I_{t,q}$.

(I) By Lemmas \ref{relation}, \ref{mrelation1} \& \ref{mrelation2}, we have
\begin{equation}\label{mtheoremd}
\begin{split}
&\Delta(\mathcal D_K (x^{(\al)})^p)=\sum_{0\le j'\le p \atop n, \ell\ge
0}\binom{p}{j'} ({-}1)^{n+\ell}\mathcal D_K (x^{(\al)})^{j'} h(k')^{\lg \ell \rg} h(k)^{\lg n \rg} \otimes(1{-}e(k)t)^{j'(\al_k-\al_{-k})-n}\times \\
&\qquad\qquad\qquad\qquad\times(1{-}e(k')t)^{j'(\al_{k'}-\al_{-k'})-\ell} d_{k}^{(n)}d_{k'}^{(\ell)} \mathcal D_K (x^{(\al)})^{p-j'}t^{n+\ell}\\
&\equiv \sum\limits_{n,\ell=0}^{\infty}({-}1)^{n+\ell}h(k')^{\lg \ell \rg}h(k)^{\lg n \rg} {\otimes} (1{-}e(k)t)^{-n}(1{-}e(k')t)^{-\ell}d_k^{(n)} d_{k'}^{(\ell)}
\mathcal D_K (x^{(\al)})^p t^{n+\ell}{+}\mathcal D_K (x^{(\al)})^p {\otimes} 1\\
&=\mathcal D_K (x^{(\al)})^p \otimes 1+\sum\limits_{n,\ell=0}^{p-1}(-1)^{n+\ell}h(k')^{\lg \ell \rg}h(k)^{\lg n \rg}\otimes (1{-}e(k)t)^{-n} (1{-}e(k')t)^{-\ell}\cdot\\
&\qquad \cdot\Bigl(\delta_{\ell 0}\delta_{n0}\mathcal D_K (x^{(\al)})^p -\delta_{\ell 0}\delta_{n1}(\delta_{\al,\ep_k+\ep_{-k}}{+}\delta_{\al,\ep_0})e(k)-\delta_{\ell 1}\delta_{n0}(\delta_{\al,\ep_{k'}+\ep_{-k'}}{+}\delta_{\al,\ep_0})e(k')\Bigr)t^{n+\ell}\\
&=\mathcal D_K (x^{(\al)})^p \otimes 1+1 \otimes \mathcal D_K (x^{(\al)})^p  +h(k)^{\lg 1 \rg} \otimes (1{-}e(k)t)^{-1}(\delta_{\al,\ep_{k}+\ep_{-k}}{+}\delta_{\al,\ep_0})e(k)t\\
&\qquad+h(k')^{\lg 1 \rg} \otimes (1{-}e(k')t)^{-1}(\delta_{\al,\ep_{k'}+\ep_{-k'}}{+}\delta_{\al,\ep_0})e(k')t\\
&=\mathcal D_K (x^{(\al)})^p \otimes 1+1 \otimes \mathcal D_K (x^{(\al)})^p
\\
&\qquad +\delta_{\al,\ep_i+\ep_{-i}}\delta_{ik}h(k)^{\lg 1 \rg} \otimes (1{-}e(k)t)^{-1}e(k)t+\delta_{\al,\ep_0} h(k)^{\lg 1 \rg}\otimes (1{-}e(k)t)^{-1} e(k)t\\
&\qquad+\delta_{\al,\ep_i+\ep_{-i}}\delta_{ik'}  h(k')^{\lg 1 \rg} \otimes (1{-}e(k')t)^{-1}e(k')t +\delta_{\al,\ep_0}h(k')^{\lg 1 \rg} \otimes (1{-}e(k')t)^{-1}e(k')t.
\end{split}
\end{equation}
Hence, when $\al \neq \ep_i+\ep_{-i}$, $\ep_0$, we get
\begin{equation*}
\begin{split}
\Delta((\mathcal D_K (x^{(\al)}))^p)& \equiv (\mathcal D_K (x^{(\al)}))^p \otimes 1 + 1 \otimes (\mathcal D_K (x^{(\al)}))^p\\
&\in I_{t,q} \otimes U_{t,q}({\bf K}(2n{+}1;\underline{1})) + U_{t,q}({\bf K}(2n{+}1;\underline{1})) \otimes I_{t,q}.
\end{split}
\end{equation*}\label{}
And when $\al=\ep_i+\ep_{-i},~1 \leq i \leq n$,
\begin{equation*}
\begin{split}
\Delta (\mathcal D_K  (x^{(\ep_i+\ep_{-i})}))& = \mathcal D_K  (x^{(\ep_i+\ep_{-i})}) \otimes 1 +\sum\limits_{n,\ell=0}^{p-1} (-1)^{n+\ell} h(k')^{\lg \ell \rg} h(k)^{\lg n \rg} \otimes(1{-}e(k')t)^{-\ell} (1{-}e(k)t)^{-n}\\
& \hskip2cm \times\Bigl(\delta_{\ell 0} \delta_{n0} \mathcal D_K (x^{(\ep_i+\ep_{-i})}) - \delta_{\ell 0} \delta_{n1}\delta_{ik}e(k)-\delta_{\ell 1}\delta_{n0}\delta_{ik'}e(k')\Bigr)t^{n+\ell}\\
&= \mathcal D_K  (x^{(\ep_i+\ep_{-i})}) \otimes 1 + 1 \otimes \mathcal D_K  (x^{(\ep_i+\ep_{-i})})+\delta_{ik}h(k)^{\lg 1 \rg} \otimes (1{-}e(k)t)^{-1} e(k)t\\
&\hskip5.3cm+\delta_{ik'}h(k')^{\lg 1 \rg} \otimes (1{-}e(k')t)^{-1}e(k')t.\\
\end{split}
\end{equation*}

Combing this with (3.23), we obtain
\begin{eqnarray*}
\begin{split}
\Delta\big(&\mathcal D_K  (x^{(\ep_i+\ep_{-i})})^p-\mathcal D_K  (x^{(\ep_i+\ep_{-i})})\big)\\
&=(\mathcal D_K  (x^{(\ep_i+\ep_{-i})})^p-\mathcal D_K  (x^{(\ep_i+\ep_{-i})})) \otimes 1+ 1 \otimes (\mathcal D_K  (x^{(\ep_i+\ep_{-i})})^p-\mathcal D_K  (x^{(\ep_i+\ep_{-i})}))\\
&\in I_{t,q} \otimes U_{t,q}({\bf K}(2n{+}1;\underline{1})) + U_{t,q}({\bf K}(2n{+}1;\underline{1})) \otimes I_{t,q}.
\end{split}
\end{eqnarray*}

When $\al=\ep_0$,
\begin{eqnarray*}
\begin{split}
&\Delta\big(\mathcal D_K (x^{(\ep_0)})\big)=\mathcal D_K (x^{(\ep_0)}) \otimes 1 + \sum\limits_{n,\ell=0}^{p-1}(-1)^{n+\ell} h(k')^{\lg \ell \rg}h(k)^{\lg n \rg} \otimes (1{-}e(k')t)^{-\ell}(1{-}e(k)t)^{-n} \\
&\hskip4.5cm \times\Big( \delta_{\ell 0} \delta_{n 0}\mathcal D_K (x^{(\ep_0)})-\delta_{\ell 0}\delta_{n1} e(k)-\delta_{n0}\delta_{\ell 1}e(k')\Big)t^{n+\ell}
\end{split}
\end{eqnarray*}
\begin{eqnarray*}
\begin{split}
&=\mathcal D_K (x^{(\ep_0)}) {\otimes} 1 {+} 1 {\otimes} \mathcal D_K (x^{(\ep_0)}){+}h(k)^{\lg 1 \rg} {\otimes} (1{-}e(k)t)^{-1}e(k)t{+}h(k')^{\lg 1 \rg} {\otimes} (1{-}e(k')t)^{-1}e(k')t.
\end{split}
\end{eqnarray*}

Combing this with (3.23), we obtain
\begin{eqnarray*}
\begin{split}
\Delta\big(\mathcal D_K (x^{(\ep_0)})^p-\mathcal D_K (x^{(\ep_0)})\big)&=\big(\mathcal D_K (x^{(\ep_0)})^p-\mathcal D_K ( x^{(\ep_0)})\big) {\otimes} 1 + 1 {\otimes} \big(\mathcal D_K (x^{(\ep_0)})^p-\mathcal D_K (x^{(\ep_0)})\big)\\
& \in I_{t,q} \otimes U_{t,q}({\bf K}(2n{+}1;\underline{1})) + U_{t,q}({\bf K}(2n{+}1;\underline{1})) \otimes I_{t,q}. \qquad \qquad \qquad \qquad
\end{split}
\end{eqnarray*}

Thus, we show the ideal $I_{t,q}$ is also a coideal of the Hopf algebra $U_{t,q}({\bf K}(2n{+}1;\underline{1}))$.

(II)  By Lemmas \ref{relation}, \ref{mrelation1} \& \ref{mrelation2}, we have
\begin{equation}\label{mtheorema}
\begin{split}
S(\mathcal D_K (x^{(\al)})^p)
&=(-1)^p (1{-}e(k')t)^{-p(\al_{k'}-\al_{-k'})}(1{-}e(k)t)^{-p(\al_k-\al_{-k})}\times \\
&\quad \times \sum\limits_{n,\ell=0}^{\infty} \big(\delta_{\ell 0} \delta_{n0} \mathcal D_K (x^{(\al)})^p-\delta_{\ell 0}\delta_{n1}(\delta_{\al,\ep_k+\ep_{-k}}{+}\delta_{\al,\ep_0})e(k)\\
&\quad -\delta_{\ell 1} \delta_{n 0}(\delta_{\al,\ep_{k'}+\ep_{-k'}}{+}\delta_{\al,\ep_0})e(k')\big)h(k)_1^{\lg n \rg}h(k')_{1}^{\lg \ell \rg} t^{n+\ell}\\
&\equiv -\mathcal D_K (x^{(\al)})^p +(\delta_{\al,\ep_{k}+\ep_{-k}}{+}\delta_{\al,\ep_0})e(k)h(k)_{1}^{\lg 1 \rg}t\\
&\hskip2.05cm +(\delta_{\al,\ep_{k'}+\ep_{-k'}}{+}\delta_{\al,\ep_0})e(k')h(k')_{1}^{\lg 1 \rg} t.
\end{split}
\end{equation}
Thus, when $\al \neq \ep_i+\ep_{-i}$, $\ep_0$, for $1 \leq i \leq n$, we have $S(\mathcal D_K (x^{(\al)})^p)\equiv-\mathcal D_K (x^{(\al)})^p \in I_{t,q}$.

When $\al=\ep_i+\ep_{-i}$, for $1 \leq i \leq n$,
\begin{eqnarray*}
\begin{split}
&S(\mathcal D_K (x^{(\ep_i+\ep_{-i})}))=-\sum\limits_{n,\ell=0}^{\infty} d_{k}^{(n)}d_{k'}^{(\ell)} (
\mathcal D_K (x^{(\ep_i+\ep_{-i})}))\,h(k)_{1}^{\lg n \rg}\, h(k')_{1}^{\lg \ell \rg} t^{n+\ell}\\
&\qquad=-\sum\limits_{n,\ell=0}^{\infty} \Big(\delta_{\ell 0} \delta_{n0}
\mathcal D_K (x^{(\ep_i+\ep_{-i})})-\delta_{\ell 0}\delta_{n 1}\delta_{ik}e(k) -\delta_{\ell 1}
\delta_{n0}\delta_{ik'} e(k')\Big)\,h(k')_{1}^{\lg \ell \rg} h(k)_{1}^{\lg n \rg} t^{n+\ell}\\
&\qquad=-\mathcal D_K (x^{(\ep_i+\ep_{-i})})+\delta_{ik}e(k)h(k)_{1}^{\lg 1\rg} t+\delta_{ik'}e(k')h(k')_{1}^{\lg 1 \rg}t.
\end{split}
\end{eqnarray*}
Combing this with (\ref{mtheorema}), we have
$$
S\Big(\mathcal D_K (x^{(\ep_i+\ep_{-i})})^p-\mathcal D_K (x^{(\ep_i+\ep_{-i})})\Big)
=-\Big(\mathcal D_K (x^{(\ep_i+\ep_{-i})})^p-\mathcal D_K (x^{(\ep_i+\ep_{-i})})\Big) \in I_{t,q}.$$
For $\al=\ep_0,\; S(\mathcal D_K (x^{(\ep_0)}))=-\mathcal D_K (x^{(\ep_0)})+e(k)h(k)_{1}^{\lg 1 \rg}t+e(k')h(k')_{1}^{\lg 1 \rg}t$,
so
$$
S\Big(\mathcal D_K (x^{(\ep_0)})^p-\mathcal D_K (x^{(\ep_0)})\Big)=-\Big(\mathcal D_K (x^{(\ep_0)})^p-\mathcal D_K (x^{(\ep_0)})\Big) \in I_{t,q}.
$$

Hence, we proved that the ideal $I_{t,q}$ is preserved by the antipode $S$ of the quantization $U_{t,q}({\bf K}(2n{+}1;\underline{1}))$.

(III) It is obvious to notice that $\varepsilon(\mathcal D_K (x^{(\al)}))=0$, for $0 \leq \al \leq \tau$.

This completes the proof.
\end{proof}
\section{Quantization of horizontal type for Lie bialgebra of Cartan type {\bf K}}

In this section, we assume that $n\geq 2$. Take $h:=\mathcal D_K (x^{\ep_k+\ep_{-k}})$ and $e:=\mathcal D_K ( x^{\epsilon_k+\epsilon_m})$, $(1\leq
k\neq |m| \leq n)$ and denote by $\mathcal{F}(k,m)$ the
corresponding Drinfel'd twist.  Set  $d^{(\ell)}=\frac{1}{\ell!} (\ad e)^\ell$. For $m \in \{-1,\cdots,-n,1,\cdots, n\}$,
 set $\sigma (m):= \begin{cases} -1,&m >0,\\
                1,&m<0.
                \end{cases}$ Using the horizontal Drinfel'd twists $\mathcal{F}(k,m)$
and the same discussion in Sections 2, 3, we obtain some new
quantizations of horizontal type for the universal enveloping
algebra of the special algebra $\mathbf{K}(2n{+}1;\underline{1})$.

To simplify the formulas, let us introduce
the operator $d^{(\ell)}$ on $U(\mathbf{K})$ defined by
$d^{(\ell)}=\frac{1}{\ell!} (\ad e)^\ell$.  Use induction on $\ell$,
we can get $$d^{(\ell)}(\mathcal D_K (x^\al))=\sum\limits_{j=0}^{\ell} A_j B_{\ell-j}\mathcal D_K (x^{\al+\ell(2\ep_k+\ep_{-k})-j(\ep_k+\ep_{-k})-(\ell-j)\ep_0}),$$ where ~$A_j=(-1)^j\frac{1}{j!}\prod\limits_{i=0}^{j-1}(\al_{-k}-j),
~B_{\ell-j}=\sigma(m)^{\ell-j}\frac{1}{(\ell-j)!}\prod\limits_{i=0}^{\ell-j-1}(\al_{-m}-i).$

Recall the vertical basic twist of Cartan type $\mathbf{H}$ Lie
algebra in \cite{HTW} is given by For $h:=D_H (x^{\ep_k+\ep_{-k}})$ and $e:=D_H
(x^{\epsilon_k+\epsilon_m})$, $(1\leq k,|m| \leq n,\  m \neq \pm
k)$, (and the twists of the Hamiltonian algebra in characteristic $p>0$ is given by
$h=D_H (x^{\ep_k+\ep_{-k}}),~e=2D_H
(x^{\epsilon_k+\epsilon_m})~(1\leq k,|m| \leq n,\  m \neq \pm
k)$) using the quantizations
  of Cartan type $\mathbf{H}$ Lie algebra and the formulas $\mathcal D_K
(x^{\alpha})=\big(2-\sum\limits_{i=1}^{n}(\al_i+\al_{-i})\big)x^{\al-\ep_0}\p_0+
\sum\limits_{i=1}^{n}\al_0x^{\al-\ep_0}(\p_i+\p_{-i})+D_H (x^\al)$. We have
the following theorem which gives the quantization of $U({\bf K})$
by Drinfel'd twist $\mathcal{F}(k,m)$.

\begin{lemm}\label{hHopf1}
Fix two distinguished elements  $h:=\mathcal D_K (x^{\ep_k+\ep_{-k}})$ and
$e:=\mathcal D_K (x^{\epsilon_k+\epsilon_m})$, $(1\leq k \neq |m| \leq n)$, the corresponding horizontal quantization of
$U(\mathbf{K}^+_{\mathbb{Z}})$ over
$U(\mathbf{K}^+_{\mathbb{Z}})[[t]]$ by Drinfel'd twist
$\mathcal{F}(k,m)$ with the product undeformed is given by
\begin{eqnarray}
&&\Delta(\mathcal D_K (x^\al))
=\mathcal D_K ( x^\al) \otimes
(1{-}et)^{\al_k-\al_{-k}}+\sum\limits_{\ell=0}^{\infty}\sum\limits_{j=0}^{\ell} (-1)^\ell A_j
B_{\ell-j} h^{\lg \ell
\rg} \otimes (1{-}et)^{-\ell}\times\\
&& \hskip6cm \times\,\mathcal D_K (x^{\al+(\ell-j)(\ep_k-\ep_{-m})+j(\ep_m-\ep_{-k})})t^\ell, \nonumber\\
&& S(\mathcal D_K (
x^\al))=-(1{-}et)^{\al_{-k}-\al_k}\sum\limits_{\ell=0}^{\infty}\sum\limits_{j=0}^{\ell}
A_j B_{\ell-j} \mathcal D_K (x^{\al+(\ell-j)(\ep_k-\ep_{-m})+j(\ep_m-\ep_{-k})})h_{1}^{\lg
\ell \rg}t^\ell,\\
&&\varepsilon(\mathcal D_K (
x^{\al}))=0,
\end{eqnarray}
where $A_j=(-1)^j\frac{1}{j!}\prod\limits_{i=0}^{j-1}(\al_{-k}-j),
~B_{\ell-j}=\sigma(m)^{\ell-j}\frac{1}{(\ell-j)!}\prod\limits_{i=0}^{\ell-j-1}(\al_{-m}-i)$, $A_0=B_0=1$.
\end{lemm}

We firstly make {\it the modulo $p$ reduction} for the quantizations
of $U(\mathbf{K}^{+}_\mathbb{Z})$ in Lemma \ref{hHopf1} to yield the
horizontal quantizations of $U(\mathbf{K}(2n{+}1;\underline{1}))$ over
$U_t(\mathbf{K}(2n{+}1;\underline{1}))$.

\begin{theorem}\label{hHopf2}
Fix distinguished elements $h=\mathcal D_K ( x^{(\epsilon_k+\epsilon_{-k})})$,
$e=\mathcal D_K (x^{(\epsilon_k+\epsilon_m)})$ $(1\leq |m| \neq k \leq n)$,
the corresponding horizontal quantization of
$U(\mathbf{K}(2n{+}1;\underline{1}))$ over
$U_t(\mathbf{K}(2n{+}1;\underline{1}))$ with the product undeformed is
given by
\begin{gather}
\Delta(\mathcal D_K ( x^{(\al)}))=\mathcal D_K ( x^{(\al)})\otimes (1{-}et)^{\al_k-\al_{-k}}+\hskip4.75cm \label{hHopf2d}
\\
\qquad\sum\limits_{\ell=0}^{p-1}\sum\limits_{j=0}^{\ell}(-1)^\ell \bar{A}_j \bar{B}_{\ell-j} h^{\lg \ell \rg} \otimes
(1{-}et)^{-\ell}\mathcal D_K \big(x^{(\al+(\ell-j)(\ep_k-\ep_{-m})+j(\ep_m-\ep_{-k}))}\big)t^\ell,\nonumber\\
S(\mathcal D_K (
x^{(\al)}))=-(1{-}et)^{\al_{-k}-\al_k}\sum\limits_{\ell=0}^{p-1}\sum\limits_{j=0}^{\ell}
\bar{A}_j\bar{B}_{\ell-j}\mathcal D_K \big( x^{(\al+(\ell-j)(\ep_k-\ep_{-m})+j(\ep_m-\ep_{-k}))}\big)t^\ell,\label{hHopf2a}\\
\varepsilon(\mathcal D_K (x^{(\al)}))=0,\hskip7.25cm\label{hHopf2u}
\end{gather}
where $0 \leq \al \leq \tau$ for $2n+4\nequiv 0\;(\mathrm{mod} p)$ and $0 \leq \al < \tau$, for $2n+4\equiv 0\; (\mathrm{mod} p)$. And $\bar{A}_j \equiv (-1)^j \binom{\al_m+j}{j}
(\mathrm{mod} p)$, for $0 \leq j \leq \al_{-k}, B_{\ell-j}\equiv
\sigma(m)^{\ell-j}\binom{\al_k+\ell-j}{\ell-j}(\mathrm{mod} p)$
 for $0 \leq \ell{-}j \leq \al_{-m}$, otherwise, $\overline{A}_j=\overline{B}_{\ell-j}=0$.
\end{theorem}



To describe $\mathbf{u}_{t,q}(\mathbf{K}(2n{+}1;\underline{1}))$
explicitly, we still need an auxiliary Lemma.

\begin{lemm}\label{hrelation11}
Set $e=\mathcal D_K ( x^{(\ep_k+\ep_m)})~,d^{(\ell)}=\frac{1}{\ell !}\ad
e$. Then

\smallskip
$(\text{\rm i})$ \quad\, $d^{(\ell)}(\mathcal D_K (x^{(\al)}))=
\sum\limits_{j=0}^{\ell} \bar{A}_j \bar{B}_{\ell-j}\mathcal D_K (x^{(\al+(\ell-j)(\ep_k-\ep_{-m})+j(\ep_m-\ep_{-k}))}),$

\smallskip
$(\text{\rm ii})$ \quad $d^{(\ell)}(\mathcal D_K (x^{(\ep_i+\ep_{-i})})=\delta_{\ell 0} \mathcal D_K (x^{(\ep_i+\ep_{-i})})+\delta_{\ell 1}(\delta_{i,-m}{-}\delta_{im}{-}\delta_{ik})e$,  \textit{ for } $1 \leq i \leq n$,

$\qquad \;\, d^{(\ell)} (\mathcal D_K (x^{(\ep_0)}))=\delta_{\ell 0}\mathcal D_K (x^{(\ep_0)})$.

\smallskip
$(\text{\rm iii})$ \ \, $d^{(\ell)} (\mathcal D_K (
x^{(\al)})^p)=\delta_{\ell 0}\mathcal D_K (
x^{(\al)})^p-\delta_{\ell 1}\big(\delta_{\al,\ep_k+\ep_{-k}}{-}\sigma(m)\delta_{\al,\ep_m{+}\ep_{-m}}\big)e$.
\end{lemm}

Based on Theorem \ref{hHopf2} and Lemma \ref{hrelation11}, we
arrive at
\begin{theorem}\label{hHopf}
Fix distinguished elements $h:=\mathcal D_K (x^{(\ep_k+\ep_{-k})}),~e:=\mathcal D_K (
x^{(\ep_k+\ep_m)})$, with $1 \leq k \neq |m| \leq n$; there exists a
noncommutative and noncocommutative Hopf algebra (of horizontal
type) $(\mathbf{u}_{t,q}(\mathbf{K}(2n{+}1;\underline{1})),m,\iota,\Delta,S,\varepsilon)$ over
$\mathcal{K}[t]_p^{(q)}$ with the product undeformed, whose
coalgebra structure is given by
\begin{gather}
\Delta(\mathcal D_K (x^{(\al)}))=\mathcal D_K (x^{(\al)}) \otimes
(1{-}et)^{\al_k-\al_{-k}}+\sum\limits_{\ell=0}^{p-1}(-1)^\ell h^{\lg
\ell \rg}\otimes
(1{-}et)^{-\ell} d^{(\ell)}(\mathcal D_K ( x^{(\al)}))t^\ell ,\\
S(\mathcal D_K (x^{(\al)}))=-(1{-}et)^{\al_{-k}-\al_k}\sum\limits_{\ell=0}^{p-1}
d^{(\ell)}(\mathcal D_K (x^{(\al)}))h_{1}^{\lg \ell \rg} t^\ell,\qquad \qquad \qquad\qquad \qquad \quad\\
\varepsilon(\mathcal D_K ( x^{(\al)}))=0, \qquad \qquad \qquad\qquad \qquad\qquad
\qquad \qquad \qquad \qquad \qquad \qquad \quad
\end{gather}
where $0 \leq \al \leq \tau$, $\dim_{\mathcal{K}}\mathbf{u}_{t,q}(\mathbf{K}(2n{+}1;\underline{1}))=p^{p^{2n+1}+1}$ in case of $2n+4 \nequiv 0\; (\mathrm{mod} p)$, and $0 \leq \al < \tau$,
 $\dim_{\mathcal{K}}\mathbf{u}_{t,q}(\mathbf{K}(2n{+}1;\underline{1}))=p^{p^{2n+1}}$ otherwise.
\end{theorem}

\section{Other Drinfel'd twists and some Open Questions}

Besides the Drinfel'd twists (i), (ii) we interested in Section 2, there are another two kinds of interesting Drinfel'd twists:

\smallskip
$\mathrm{(iii)}\ h=\mathcal D_K (x^{\ep_k +\ep_{-k}}), \;e=\mathcal D_K ( x^{
\ep_{k}+\ep_{0}}) \quad (1 \leq k \leq n);$

\smallskip
$\mathrm{(ix)}\ h=\mathcal D_K (x^{\ep_0}), \;e=\mathcal D_K ( x^{\ep_k +
\ep_0}) \quad (1 \leq |k| \leq n).$

We can see that the elements in (iii) and (ix) also satisfies $[h,e]=e$,
so we can get triangular Lie bialgebra structure on $\mathcal{K}$ given by the classical
Yang-Baxter $r$-matrix $r:=h \otimes e -e \otimes h$. We can also get some new Drinfel'd twists.
We still use $\mathcal{F}(k)$'s to denote these Drinfel'd twists, and  call them the vertical twists due to the degree $||\ep_k+\ep_0||=1$.

The basic theory is given in  Sections 3 and 4, here we only list some main calculations.

\subsection{Construction of Drinfel'd twist }

To simplify the formulas, let us introduce the operator $d^{(\ell)}$ on $U({\bf K})$ defined by $d^{(\ell)}=\frac{1}{\ell !}(\ad e)^\ell$. Then we can get
\begin{lemm}\label{olemma1}
For $ \mathcal D_K (x^\al) \in U({\bf K})$, the following equalities hold
\begin{eqnarray}
&&d^{(\ell)}(\mathcal D_K (x^\al))=\sum\limits_{j=0}^{\ell} A_j B_{\ell-j}\mathcal D_K (x^{\al+(\ell-j)\ep_k+j(\ep_0-\ep_{-k})}) ,   \label{orelation1}\\
&&\mathcal D_K (x^\al) \cdot e^m=\sum\limits_{\ell=0}^{m}(-1)^\ell \binom{m}{\ell}\ell !\sum\limits_{j=0}^{\ell} A_j B_{\ell-j}\mathcal D_K (x^{\al+(\ell-j)\ep_k+j(\ep_0-\ep_{-k})}),\label{orelation3}\\
&&(\ad \mathcal D_K (x^\al))^\ell(e)=\sum\limits_{j=0}^{\ell}\binom{\ell}{j}\prod\limits_{i=0}^{\ell-j-1}((i{-}1)||\al||+\al_0)\prod\limits_{i=0}^{j-1}\big(i\al_k{-}(i{-}1)\al_{-k}\big
)\cdot\label{orelation4}\\
&&\qquad \qquad\qquad \qquad \cdot\mathcal D_K (x^{\ell\al+\ep_k-j(\ep_k+\ep_{-k})-(\ell{-}j{-}1)\ep_0}), \nonumber
\end{eqnarray}
where $A_j=\frac{1}{j!}\prod\limits_{i=0}^{j-1}(i-\al_{-k}),
~B_{\ell-j}=\frac{1}{(\ell-j)!}\prod\limits_{i=0}^{\ell-j-1}(||\al||-\al_0+i)$, $A_0=B_0=1$.
\end{lemm}
\begin{proof}
For (\ref{orelation1}), use induction on $\ell$. When $\ell=1$,
\begin{equation*}
\begin{split}
d(\mathcal D_K (x^\al))&=[\mathcal D_K (x^{\ep_0+\ep_{k}}),\mathcal D_K (x^\al)]\\
&=\bigl(\al_0{-}(2{-}\sum\limits_{i=1}^{n}(\al_i{+}\al_{-i}))\bigr)\mathcal D_K (x^{\al+\ep_k})-\al_{-k}\mathcal D_K (x^{\al+\ep_0-\ep_{-k}})\\
&=(||\al||-\al_0)\mathcal D_K (x^{\al+\ep_k})-\al_{-k}\mathcal D_K (x^{\al+\ep_0-\ep_{-k}}).
\end{split}
\end{equation*}

When $\ell \geq 1$, we have
{\setlength{\arraycolsep}{0pt}
\begin{eqnarray*}
&&d^{(\ell+1)}(\mathcal D_K (x^{\al}))=\frac{d}{\ell{+}1}d^{(\ell)}(\mathcal D_K (x^\al))
=\frac{1}{\ell{+}1}\sum\limits_{j=0}^{\ell}A_j B_{\ell-j} d (\mathcal D_K (x^{\al+(\ell-j)\ep_k+j(\ep_0-\ep_{-k})}))\\
&&=\frac{1}{\ell{+}1}\sum\limits_{j=0}^{\ell}A_jB_{\ell-j} \Big((\al_0{+}j){-}(2{-}(\sum\limits_{i=1}^{n}\al_i{+}\al_{-i}{+}\ell{-}j{-}j))\mathcal D_K (x^{\al+(\ell-j)\ep_k+j(\ep_0-\ep_{-k})+\ep_k})\\
&&\hskip6cm-(\al_{-k}{-}j)\mathcal D_K (x^{\al+(\ell-j)\ep_k+j(\ep_0-\ep_{-k})+\ep_0-\ep_{-k}})\Big)\\
&&=\frac{1}{\ell{+}1} \sum\limits_{j=0}^{\ell}A_jB_{\ell-j}\Big((||\al||{-}\al_0{+}\ell{-}j) \mathcal D_K (x^{\al+(\ell-j)\ep_k+j(\ep_0-\ep_{-k})+\ep_k})+\\
&&\hskip3.5cm+(j{-}\al_{-k})\mathcal D_K (x^{\al+(\ell-j)\ep_k+(j+1)(\ep_0-\ep_{-k})})\Big)\\
&&=\frac{1}{\ell{+}1} \sum\limits_{j=0}^{\ell}(\ell{-}j{+}1) A_j B_{\ell-j+1} \mathcal D_K (x^{\al+(\ell-j+1)\ep_k+j(\ep_0-\ep_{-k})})\\
&&\quad +\frac{1}{\ell{+}1} \sum\limits_{j=0}^{\ell}(j{+}1) A_{j+1} B_{\ell-j} \mathcal D_K (x^{\al+(\ell-j)\ep_k+(j+1)(\ep_0-\ep_{-k})})\\
&&=\sum\limits_{j=0}^{\ell{+}1} A_j B_{\ell+1-j} \mathcal D_K (x^{\al+(\ell+1-j)\ep_k+j(\ep_0-\ep_{-k})}).
\end{eqnarray*}
}

(5.2) follows from (5.1).

For (\ref{orelation4}), use induction on $\ell$.
When $\ell=1$, we have
\begin{equation*}
\begin{split}
\ad \mathcal D_K (x^\al)(e)&=[\mathcal D_K (x^{\al}),\mathcal D_K (x^{\ep_k+\ep_0})]\\
&=(2{-}\sum\limits_{i=1}^{n}(\al_i{+}\al_{-i}){-}\al_0)\mathcal D_K (x^{\al+\ep_k})+\al_{-k}\mathcal D_K (x^{\al+\ep_0-\ep_{-k}})\\
&=(-||\al||{+}\al_0) \mathcal D_K (x^{\al+\ep_k})+\al_{-k} \mathcal D_K (x^{\al+\ep_0-\ep_{-k}}).
\end{split}
\end{equation*}

When $\ell \geq 1$, we have
\begin{equation*}
\begin{split}
&(\ad \mathcal D_K (x^\al))^{\ell+1}(e) \\
&=\ad \mathcal D_K (x^\al) \sum\limits_{j=0}^{\ell}\binom{\ell}{j}\prod\limits_{i=0}^{\ell-j-1}((i{-}1)||\al||{+}\al_0)\prod\limits_{i=0}^{j-1}\big(i\al_k{-}(i{-}1)\al_{-k}\big
) \mathcal D_K (x^{\ell\al+\ep_k-j(\ep_k+\ep_{-k})-(\ell-j-1)\ep_0}) \\
&=\sum\limits_{j=0}^{\ell}\binom{\ell}{j} \prod\limits_{i=0}^{\ell-j}((i{-}1)||\al||{+}\al_0)\prod\limits_{i=0}^{j-1}(i\al_k{-}(i{-}1)\al_{-k})\mathcal D_K (x^{(\ell+1)\al+\ep_k-(\ell-j)\ep_0-j(\ep_k+\ep_{-k})})\\
&\quad +\sum\limits_{j=0}^{\ell}\binom{\ell}{j} \prod\limits_{i=0}^{\ell-j-1}((i{-}1)||\al||{+}\al_0)\prod\limits_{i=0}^{j}(i\al_k{-}(i{-}1)\al_{-k})\mathcal D_K (x^{(\ell+1)\al+\ep_k-(\ell-j-1)\ep_0-(j+1)(\ep_k+\ep_{-k})})\\
&=\sum\limits_{j=0}^{\ell{+}1} \binom{\ell{+}1}{j}\prod\limits_{i=0}^{\ell-j} ((i{-}1)||\al||{+}\al_0)\prod\limits_{i=0}^{j-1} (i\al_k{-}(i{-}1)\al_{-k})\mathcal D_K (x^{(\ell+1)\al+\ep_k-(\ell-j)\ep_0-j(\ep_k+\ep_{-k})}).
\end{split}
\end{equation*}
Thus we can get (\ref{orelation4}).
\end{proof}

As the proof in Section 3, we can also get
\begin{lemm}\label{olemma2}
For $a \in \mathbb{F}$, $\al \in \mathbb{Z}^{2n+1}$, and $\mathcal D_K (x^\al) \in {\bf K}$, the following equalities hold
\begin{eqnarray}
&&((\mathcal D_K (x^\al))^s \otimes 1) \cdot F_a=F_{a+s(\al_{-k}-\al_{k})} \cdot \big((\mathcal D_K (x^\al))^s \otimes 1\big),\label{orelation5}\\
&&(\mathcal D_K (x^\al))^s \cdot u_a=u_{a+s(\al_k-\al_{-k})}\sum\limits_{\ell=0}^{\infty} d^{(\ell)}\bigl((\mathcal D_K (x^\al))^s\bigr) h_{1-a}^{\lg \ell \rg} t^\ell,\label{orelation6}\\
&&(1 \otimes (\mathcal D_K (x^\al))^s) \cdot F_a =\sum\limits_{\ell=0}^{\infty} (-1)^{\ell} F_{a+\ell} \big(h_{a}^{\lg\ell \rg} \otimes d^{(\ell)}(\mathcal D_K (x^{\al})^s)\big) t^\ell.\label{orelation7}
\end{eqnarray}
\end{lemm}

\begin{lemm}\label{olemma3}
For $s \geq 1$, we have
\begin{eqnarray}
&&\Delta\big((\mathcal D_K (x^\al))^s\big)
=\sum\limits_{0 \leq j \leq s \atop \ell \geq 0} (-1)^\ell\binom{s}{j} (\mathcal D_K (x^\al))^j h^{\lg \ell \rg} \otimes (1{-}et)^{j(\al_k-\al_{-k})-\ell}\big(d^{(\ell)}(\mathcal D_K (x^\a))^{s-j}\big)\, t^\ell , \label{orelation8}\\
&&S((\mathcal D_K (x^\al))^s)=(-1)^s (1{-}et)^{-s(\al_k-\al_{-k})}\sum\limits_{\ell=0}^{\infty} d^{(\ell)}(\mathcal D_K (x^\al)^s)\, h_{1}^{\lg \ell \rg}\, t^\ell.\label{orelation9}
\end{eqnarray}
\end{lemm}

\begin{lemm}\label{oHopf1}
Fix two distinguished elements  $h:=\mathcal D_K (x^{\ep_k+\ep_{-k}})$ and
$e:=\mathcal D_K (x^{\epsilon_k+\epsilon_0})$, $(1\leq k \leq n$, the corresponding vertical quantization of
$U(\mathbf{K}^+_{\mathbb{Z}})$ over
$U(\mathbf{K}^+_{\mathbb{Z}})[[t]]$ by Drinfel'd twist
$\mathcal{F}(k,0)$ with the product undeformed is given by
\begin{gather}
\Delta(\mathcal D_K (x^\al))=\mathcal D_K ( x^\al) \otimes
(1{-}et)^{\al_k-\al_{-k}}\qquad\qquad\qquad\qquad\qquad\quad\\
\qquad\qquad\qquad\qquad\quad+\sum\limits_{\ell=0}^{\infty}
\sum\limits_{j=0}^{\ell} (-1)^\ell A_j
B_{\ell-j} h^{\lg \ell \rg} \otimes (1{-}et)^{-\ell}\mathcal D_K (x^{\al+(\ell-j)\ep_k+j(\ep_0-\ep_{-k})}) t^\ell,\nonumber\\
 S(\mathcal D_K (
x^\al))=-(1{-}et)^{\al_{-k}-\al_k}\sum\limits_{\ell=0}^{\infty}\sum\limits_{j=0}^{\ell}
A_j B_{\ell-j} \mathcal D_K (x^{\al+(\ell-j)\ep_k+j(\ep_0-\ep_{-k})})h_{1}^{\lg
\ell \rg} t^\ell, \\
\varepsilon(\mathcal D_K(x^{\al}))=0,\qquad\qquad\qquad\qquad\qquad\qquad\qquad \quad\qquad\qquad\qquad\
\end{gather}
where $A_j$, $B_{\ell-j}$ as defined in \ref{olemma1}.
\end{lemm}

Firstly, we make {\it the modulo $p$ reduction} for the quantizations
of $U(\mathbf{K}^+_\mathbb{Z})$ in Lemma \ref{oHopf1} to yield the
vertical quantizations of $U(\mathbf{K}(2n{+}1;\underline{1}))$ over
$U_t(\mathbf{K}(2n{+}1;\underline{1}))$.

\begin{theorem}\label{oHopf2}
Fix distinguished elements $h=\mathcal D_K ( x^{(\epsilon_k+\epsilon_{-k})})$,
$e=\mathcal D_K (x^{(\epsilon_k+\epsilon_0)})$ $(1\leq  k \leq n)$,
the corresponding vertical quantization of
$U(\mathbf{K}(2n{+}1;\underline{1}))$ over
$U_t(\mathbf{K}(2n{+}1;\underline{1}))$ with the product undeformed is
given by
\begin{gather}
\Delta(\mathcal D_K (x^{(\al)}))=\mathcal D_K ( x^{(\al)})\otimes (1{-}et)^{\al_k-\al_{-k}}+\qquad\qquad\qquad\qquad\qquad\qquad\quad \label{oHopf2d} \\
\qquad\qquad\qquad\qquad\qquad\qquad\sum\limits_{\ell=0}^{p-1}\sum\limits_{j=0}^{\ell} (-1)^\ell \bar{A}_j \bar{B}_{\ell-j} h^{\lg \ell \rg} \otimes
(1{-}et)^{-\ell}\mathcal D_K \big(x^{(\al+(\ell-j)\ep_k+j(\ep_0-\ep_{-k}))}\big)t^\ell,\nonumber\\
S(\mathcal D_K (
x^{(\al)}))=-(1{-}et)^{\al_{-k}-\al_k}\sum\limits_{\ell=0}^{p-1}\sum\limits_{j=0}^{\ell}
\bar{A}_j\bar{B}_{\ell-j}\mathcal D_K \big(x^{(\al+(\ell-j)\ep_k+j(\ep_0-\ep_{-k}))}\big)h_1^{\lg \ell \rg}t^\ell,\label{oHopf2a}\\
\varepsilon(\mathcal D_K (x^{(\al)}))=0,\qquad\qquad\qquad\qquad\qquad\qquad\qquad \qquad\qquad\qquad\qquad\label{oHopf2u}
\end{gather}
where $0 \leq \al \leq \tau$, if $2n{+}4\nequiv 0\; (\mathrm{mod}\, p);$ $0 \leq \al < \tau$, if $2n{+}4\equiv 0\; (\mathrm{mod}\, p)$. And $\bar{A}_j \equiv (-1)^j \binom{\al_0{+}j}{j}\,
(\mathrm{mod}\, p)$, for $0 \leq j \leq \al_{-k}$. Otherwise, $\overline{A}_j=0$, $B_{\ell-j}\equiv
(\ell{-}j)!\binom{\al_k{+}\ell{-}j}{j}B_{\ell-j}\,(\mathrm{mod}\, p)$.
\end{theorem}

\begin{proof}
As the same argument of Theorem \ref{hHopf2}, we can prove this theorem once we notice that
\begin{equation*}
\begin{split}
&\frac{A_j B_{\ell-j}(\al{+}(\ell{-}j)\ep_k{+}j(\ep_0{-}\ep_{-k}))!}{\al !}\\
&=\frac{\prod\limits_{i=0}^{j-1}(i{-}\al_{-k})}{j!}\frac{\prod\limits_{i=0}^{\ell-j-1}(||\al||{-}\al_0{+}j)}{(\ell{-}j)!}
\frac{(\al_k{+}(\ell{-}j))!(\al_{-k}{-}j)!(\al_0{+}j)!}{\al_k!\al_{-k}!\al_0!}
\\
&=(-1)^j \frac{\al_{-k}(\al_{-k}{-}1)\cdots (\al_{-k}{-}j{+}1)(\al_{-k}{-}j)!}{\al_{-k}!}\frac{(\al_0{+}j)!}{\al_0!j!} \frac{(\al_k{+}(\ell{-}j))!}{\al_k!}B_{\ell-j}\\
&=\bar{A}_j \binom{\al_k{+}\ell{-}j}{\ell{-}j}(\ell{-}j)!B_{\ell-j}=\bar{A}_j \bar{B}_{\ell-j}.
\end{split}
\end{equation*}

This completes the proof.
\end{proof}

Here we also need the following auxiliary Lemma
\begin{lemm}\label{olemma4}
Set $e=\mathcal D_K ( x^{(\ep_k+\ep_0)}),~d^{(\ell)}=\frac{1}{\ell !}\ad e$. Then

\smallskip
$(\text{\rm i})$ \quad\, $d^{(\ell)} (\mathcal D_K (x^{(\al)}))=
\sum\limits_{j=0}^{\ell} \bar{A}_j \bar{B}_{\ell-j}\mathcal D_K (x^{(\al+(\ell-j)\ep_k+j(\ep_0-\ep_{-k}))})$,

\smallskip
$(\text{\rm ii})$ \quad $d^{(\ell)}(\mathcal D_K (
x^{(\ep_i+\ep_{-i})}))=\delta_{\ell 0} \mathcal D_K (
x^{(\ep_i+\ep_{-i})})-\delta_{\ell
1}\delta_{i,k}e$,

$\qquad \ \, d^{(\ell)}  (\mathcal D_K (x^{(\ep_0)}))=\delta_{\ell 0}\mathcal D_K (x^{(\ep_0)})-\delta_{\ell 1}e$,

\smallskip
$(\text{\rm iii})$ \ \; $d^{(\ell)} (\mathcal D_K (
x^{(\al)})^p)=\delta_{\ell 0}(\mathcal D_K (
x^{(\al)})^p)-\delta_{\ell
1}\big(\delta_{\al,\ep_0}+\delta_{\al,\ep_k+\ep_{-k}}\big)e$.

\end{lemm}
\begin{proof}
(i) follows from the proof of Lemma \ref{olemma1}.

(ii) We can see that $\bar{A}_0=\bar{B}_0=1$.

For $\al=\ep_i+\ep_{-i}$, $\bar{A}_1=\delta_{ik}(-1)\binom{\al_0+1}{1}=-\delta_{ik}$;
$\bar{B}_1=\binom{\al_k+1}{1}B_1=(\al_k+1)(||\ep_i+\ep_{-i}||)=0$.
For $\al=\ep_0,~\bar{A}_1=0$, $\bar{B}_1=1!\binom{\al_k+1}{1}(||\ep_0||-1)=-1$. Then we can get (ii) by (i).

(iii) By (\ref{vrelation7}) and Lemma \ref{olemma1}, we have
\begin{equation*}
\begin{split}
&d (\mathcal D_K (x^{(\al)})^p)\equiv(-1)^p (\ad \mathcal D_K (x^{(\al)})^p)(e) \\
&=(-1)\frac{1}{(\al !)^p} \sum\limits_{j=0}^{p} \binom{p}{j} \prod\limits_{i=0}^{p-j-1} ((i{-}1)||\al||{+}\al_0)\prod\limits_{i=0}^{j-1} (i\al_k{-}(i{-}1)\al_{-k})\mathcal D_K (x^{p\al+\ep_k-(p-j-1)\ep_0- p(\ep_k+\ep_{-k})})\\
&\equiv-\frac{1}{\al!} \prod\limits_{i=0}^{p-1} ((i{-}1)||\al||{+}\al_0)\mathcal D_K (x^{p\al+\ep_k+\ep_0-p\ep_0})\\
 &\qquad -\frac{1}{\al!} \prod\limits_{i=0}^{p-1} (i\al_k{-}(i{-}1)\al_{-k})\mathcal D_K (x^{p\al+\ep_k+\ep_0-p(\ep_k+\ep_{-k})}) \qquad(\mod\, p)\\
&\equiv -\delta_{\al,\ep_0}e-\delta_{\al,\ep_k+\ep_{-k}}e\qquad(\mod\, p,~J).
\end{split}
\end{equation*}

This completes the proof.
\end{proof}

\begin{theorem}\label{oHopf}
Fix distinguished elements $h:=\mathcal D_K (x^{(\ep_k+\ep_{-k})}),~e:=\mathcal D_K (
x^{(\ep_k+\ep_0)})$, with $1 \leq k  \leq n$, there is a
noncommutative and noncocommutative Hopf algebra (of vertical
type) $(\mathbf{u}_{t,q}(\mathbf{K}(2n{+}1;\underline{1})),m,\iota,\Delta,S,\varepsilon)$ over
$\mathcal{K}[t]_p^{(q)}$ with the product undeformed, whose
coalgebra structure is given by
\begin{gather}
\Delta(\mathcal D_K (x^{(\al)}))=\mathcal D_K (x^{(\al)}) \otimes
(1{-}et)^{\al_k-\al_{-k}}+\sum\limits_{\ell=0}^{p-1}(-1)^\ell h^{\lg
\ell \rg}\otimes
(1{-}et)^{-\ell} d^{(\ell)} (\mathcal D_K ( x^{(\al)}))t^\ell, \\
S(\mathcal D_K (x^{(\al)}))=-(1{-}et)^{\al_{-k}-\al_k}\sum\limits_{\ell=0}^{p-1}
d^{(\ell)}(\mathcal D_K (x^{(\al)}))h_{1}^{\lg \ell \rg} t^\ell,\hskip4cm\\
\varepsilon(\mathcal D_K ( x^{(\al)}))=0,\hskip8.9cm
\end{gather}
for $0 \leq \al \leq \tau$, which is finite-dimensional with $\dim_{\mathcal{K}}\mathbf{u}_{t,q}(\mathbf{K}(2n{+}1;\underline{1}))=p^{p^{2n+1}+1}$, when $2n+4 \nequiv 0 \;(\mathrm{mod}\, p)$. And $0 \leq \al < \tau$, $\dim_{\mathcal{K}}\mathbf{u}_{t,q}(\mathbf{K}(2n{+}1;\underline{1}))=p^{p^{2n+1}}$, when $2n+4 \equiv 0 \;(\mathrm{mod}\, p)$.
\end{theorem}

\begin{proof}
(I) By (\ref{orelation8}) and Lemma \ref{olemma4}, we have
\begin{equation}
\begin{split}
\Delta(\mathcal D_K (x^{(\al)})^p)&=1 \otimes \mathcal D_K (x^{(\al)})^p+\mathcal D_K (x^{(\al)})^p \otimes 1 -h \otimes d (\mathcal D_K (x^{(\al)})^p)\, t \\
&=1 \otimes \mathcal D_K (x^{(\al)})^p+\mathcal D_K (x^{(\al)})^p \otimes 1+h \otimes (1{-}et)^{-1} (\delta_{\al,\ep_0}+\delta_{\al,\ep_i+\ep_{-i}}\delta_{ik})et.
\end{split}
\end{equation}
Thus, when $\al \neq \ep_i+\ep_{-i}$, $\al \neq \ep_0$, we have
\begin{equation*}
\begin{split}
\Delta((\mathcal D_K (x^{(\al)}))^p)&=1 \otimes (\mathcal D_K (x^{(\al)}))^p+(\mathcal D_K (x^{(\al)}))^p \otimes 1\\
&
 \in I_{t,q} \otimes U_{t,q}({\bf K}(2n+1,\underline{1}))+U_{t,q}({\bf K}(2n{+}1,\underline{1})) \otimes I_{t,q}.
\end{split}
\end{equation*}

If $\al = \ep_i+\ep_{-i}$, we obtain
$$\Delta(\mathcal D_K (x^{(\ep_i+\ep_{-i})}))=\mathcal D_K (x^{(\ep_i+\ep_{-i})}) \otimes 1 + 1\otimes \mathcal D_K (x^{(\ep_i+\ep_{-i})})+ h \otimes (1{-}et)^{-1} \delta_{i,k}et.$$
Then
\begin{equation*}
\begin{split}
\Delta\big(\mathcal D_K (x^{(\ep_i+\ep_{-i})})^p-\mathcal D_K (x^{(\ep_i+\ep_{-i})})\big)
&=\big(\mathcal D_K (x^{(\ep_i+\ep_{-i})})^p-\mathcal D_K (x^{(\ep_i+\ep_{-i})})\big) \otimes 1\\
&\quad+ 1 \otimes \big(\mathcal D_K (x^{(\ep_i+\ep_{-i})})^p-\mathcal D_K (x^{(\ep_i+\ep_{-i})})\big) \\
&\in  I_{t,q} \otimes U_{t,q}({\bf K}(2n{+}1,\underline{1}))+U_{t,q}({\bf K}(2n{+}1,\underline{1})) \otimes I_{t,q}.
\end{split}
\end{equation*}

If $\al = \ep_0$, we have
$$\Delta(\mathcal D_K (x^{(\ep_0)}))=\mathcal D_K (x^{(\ep_0}) \otimes 1 + 1\otimes \mathcal D_K (x^{(\ep_0})+ h \otimes (1{-}et)^{-1} et.$$
Thus,
\begin{equation*}
\begin{split}
\Delta\big(\mathcal D_K (x^{(\ep_0)})^p-\mathcal D_K (x^{(\ep_0)})\big)
&=\big(\mathcal D_K (x^{(\ep_0)})^p-\mathcal D_K (x^{(\ep_0)})\big) \otimes 1+ 1 \otimes \big(\mathcal D_K (x^{(\ep_0)})^p-\mathcal D_K (x^{(\ep_0)})\big) \\
&\in  I_{t,q} \otimes U_{t,q}({\bf K}(2n{+}1,\underline{1}))+U_{t,q}({\bf K}(2n{+}1,\underline{1})) \otimes I_{t,q}.
\end{split}
\end{equation*}

Thus, we show the ideal $I_{t,q}$ is also a coideal of the Hopf algebra $U_{t,q}({\bf K}(2n{+}1;\underline{1}))$.

(II) By (\ref{orelation9}) and Lemma \ref{olemma4}, we have
\begin{equation*}
\begin{split}
S(\mathcal D_K (x^{(\al)})^p)&=(-1)\mathcal D_K (x^{(\al)})^p -(-\delta_{\al,\ep_k+\ep_{-k}}-\delta_{\al,\ep_0})\,e\,h_{1}^{\lg 1 \rg}t\\
&=-\mathcal D_K (x^{(\al)})^p +(\delta_{\al,\ep_0}+\delta_{\al,\ep_i+\ep_{-i}}\delta_{ik})\,e\,h_{1}^{\lg 1 \rg}t.
\end{split}
\end{equation*}
Thus, when $\al \neq \ep_i+\ep_{-i}$, $\al \neq \ep_0$, we have
$$
S(\mathcal D_K (x^{(\al)})^p)=-\mathcal D_K (x^{(\al)})^p \in I_{t,q}.
$$
And then
\begin{equation*}
\begin{split}
S\big(\mathcal D_K (x^{(\ep_i+\ep_{-i})})^p-\mathcal D_K (x^{(\ep_i+\ep_{-i})})\big)
&
=-\mathcal D_K (x^{(\ep_i+\ep_{-i})})^p+\delta_{ik}eh_{1}^{\lg 1 \rg }t+\mathcal D_K (x^{(\ep_i+\ep_{-i})})-\delta_{ik}eh_{1}^{\lg 1 \rg}t\\
&=-\big(\mathcal D_K (x^{(\ep_i+\ep_{-i})})^p-\mathcal D_K (x^{(\ep_i+\ep_{-i})})\big) \in I _{t,q},\\
S\big(\mathcal D_K (x^{(\ep_0)})^p-\mathcal D_K (x^{(\ep_0)})\big)&=-\mathcal D_K (x^{(\ep_0)})^p+eh_{1}^{\lg 1 \rg}t+\mathcal D_K (x^{\ep_0})-eh_{1}^{\lg 1 \rg}t\\
&=-\big(\mathcal D_K (x^{(\ep_0)})^p-\mathcal D_K (x^{(\ep_0)})\big) \in I _{t,q}.
\end{split}
\end{equation*}

Therefore, we proved that $I_{t,q}$ is indeed preserved by antipode $S$ of
$U_{t,q}(\mathbf{K}(2n{+}1;\underline{1}))$.

(III) It is obvious that $\varepsilon(\mathcal D_K (x^{(\al)}))=0$,  for $0 \leq \al \leq \tau$, or $0 \leq \al < \tau$.

This completes the proof.
\end{proof}

\subsection{More quantizations.}
In this subsection, the arguments are similar to Section 4, thus we only give some calculation results.
Let $A(k)_j$ and $A(k')_{j'}$ denote the coefficients of the
corresponding quantizations of $U(\mathbf{K}^+_{\mathbb{Z}})$ over
$U(\mathbf{K}^+_{\mathbb{Z}})[[t]]$ given by Drinfel'd twists
$\mathcal{F}(k)$ and $\mathcal{F}(k')$ as in Corollary \ref{vHopf2},
respectively. Note that $A(k)_0=A(k')_0$ $=1$,
$A(k)_{-1}=A(k')_{-1}=0$.

\begin{lemm}\label{omHopf1}
Fix distinguished elements $h(k)=\mathcal D_K ( x^{\ep_k+\ep_{-k}}),~e(k)=\mathcal D_K ( x^{\ep_k+\ep_0})$
$(1 \leq k \leq n)$ and $h(k')=\mathcal D_K (x^{\ep_{k'}+\ep_{-k'}}),~e(k')=\mathcal D_K (x^{\ep_{k'}+\ep_0}) \quad (1
\leq k' \leq n)$ with $k \neq k'$. the corresponding quantization of
 $U(\mathbf{K}_{\mathbb{Z}}^{+})[[t]]$ by the Drinfel'd twist $\mathcal{F}=\mathcal{F}(k)\mathcal{F}(k')$
with the product undeformed is given by
{\setlength{\arraycolsep}{0pt}
\begin{eqnarray*}
&\Delta(\mathcal D_K (x^\al))&\,=\mathcal D_K (x^\al) \otimes
(1{-}e(k)t)^{\al_k-\al_{-k}}(1{-}e(k')t)^{\al_{k'}-\al_{-k'}}\\
&&\qquad +\sum\limits_{n,\ell=0}^{\infty}\sum\limits_{j'=0}^{\ell}\sum\limits_{j=0}^{n}
(-1)^{\ell+n}A(k)_jA(k')_{j'}B(k')_{\ell-j'}
 C_{n-j}^{\ell,j'}h(k')^{\lg
\ell \rg} h(k)^{\lg n\rg} \otimes\\
&& \qquad \qquad(1{-}e(k)t)^{-n}(1{-}e(k')t)^{-\ell}
\mathcal D_K (x^{\al+(\ell-j')\ep_k+j'(\ep_0-\ep_{-k'})+(n-j)\ep_k+j(\ep_0-\ep_{-k})})t^{n+\ell},\\
&S(\mathcal D_K (x^\al))&\,=-(1{-}e(k)t)^{\al_{-k}-\al_{k}}(1{-}e(k')t)^{\al_{-k'}-\al_{k'}}
\sum\limits_{n,\ell=0}^{\infty}\sum\limits_{j'=0}^{\ell} \sum\limits_{j=0}^{n}A(k)_j A(k')_{j'} B(k')_{\ell-j'}C_{n-j}^{\ell,j'}\times\\
&&\qquad \times \mathcal D_K (x^{\al+(\ell-j')\ep_{k'}+j'(\ep_0-\ep_{-k'})+(n-j)\ep_k+j(\ep_0-\ep_{-k})})
h(k)_1^{\lg n \rg} h(k')_{1}^{\lg \ell \rg} t^{\ell+n},
\end{eqnarray*}
}
and $\varepsilon(\mathcal D_K (x^\al))=0$, where $C_{n-j}^{\ell,j'}=\frac{1}{(n{-}j)!}\prod\limits_{i=0}^{n-j-1}(||\al||{-}\al_0{+}(\ell{-}j'){+}i)$, for $\mathcal D_K (x^\alpha) \in \mathbf{K}_{\mathbb{Z}}^{+}$.
\end{lemm}
\begin{proof} First of all, let us consider
\begin{equation*}
\begin{split}
\Delta(\mathcal D_K (x^\al))&=\mathcal{F}(k)\mathcal{F}(k')\Delta_0(\mathcal D_K (x^\al))\mathcal{F}(k')^{-1}\mathcal{F}(k)^{-1}\\
&=\mathcal{F}(k) \mathcal D_K (x^\al) \otimes
(1{-}e(k')t)^{\al_{k'}-\al_{-k'}}F(k)\\
&\quad+\mathcal{F}(k)\Big(\sum\limits_{\ell=0}^{\infty}({-}1)^\ell
h(k')^{\lg \ell \rg} \otimes (1{-}e(k')t)^{-\ell} d_{k'}^{(\ell)}(\mathcal D_K (x^{\al})) t^\ell \Big)F(k).
\end{split}
\end{equation*}
By Lemma 2.1 \& Corollary \ref{vHopf2}, we can get
\begin{equation*}
\begin{split}
\mathcal{F}(k) \big(\mathcal D_K (x^\al)& \otimes
\big(1{-}e(k')t\big)^{\al_{k'}-\al_{-k'}}\big)F(k)\\
&=\mathcal{F}(k)(\mathcal D_K (x^\al) \otimes 1)(1 \otimes
(1{-}e(k')t)^{\al_{k'}-\al_{-k'}})F(k)\\
&=\mathcal{F}(k)(\mathcal D_K (x^\al) \otimes 1) F(k)\big( 1 \otimes
(1{-}e(k')t)^{\al_{k'}-\al_{-k'}}\big)\\
&=\mathcal{F}(k)F(k)_{\al_k-\al_{-k}} (\mathcal D_K (x^\al) \otimes 1)(1
\otimes (1{-}e(k')t)^{\al_{k'}-\al_{-k'}})\\
&=\mathcal D_K (x^\al) \otimes
(1{-}e(k)t)^{\al_k-\al_{-k}}(1{-}e(k')t)^{\al_{k'}-\al_{-k'}},
\end{split}
\end{equation*}
\begin{equation*}
\begin{split}
&\mathcal{F}(k)\big(\sum\limits_{\ell=0}^{\infty}(-1)^\ell h(k')^{\lg
\ell \rg} \otimes (1{-}e(k')t)^{-\ell} d_{k'}^{(\ell)} \mathcal D_K (x^{\al})
t^\ell \big)F(k)\\
&=\mathcal{F}(k)\sum\limits_{\ell=0}^{\infty} ({-}1)^\ell h(k')^{\lg
\ell \rg} {\otimes} (1{-}e(k')t)^{-\ell}\sum\limits_{j'=0}^{\ell}
A(k')_{j'}
B(k')_{\ell{-}j'}\mathcal D_K (x^{\al{+}(\ell{-}j')\ep_{k'}{+}j'(\ep_0{-}\ep_{-k'})})
t^\ell F(k) \\
&=\Big(\sum\limits_{\ell=0}^{\infty} (-1)^\ell h(k')^{\lg \ell \rg}
\otimes (1{-}e(k')t) ^{-\ell}\Big)\Big(\sum\limits_{j'=0}^{\ell}
A(k')_{j'} B(k')_{\ell-j'} \mathcal{F}(k)\big(1 \otimes
\mathcal D_K (x^{\al+(\ell-j')\ep_{k'}+j'(\ep_0-\ep_{-k'})})\big)
F(k) t^\ell \Big)
\end{split}
\end{equation*}
\begin{equation*}
\begin{split}
&=\bigg(\sum\limits_{\ell=0}^{\infty} (-1)^\ell h(k')^{\lg \ell \rg}
\otimes (1{-}e(k')t) ^{-\ell}\bigg) \cdot \\
&\quad\cdot \bigg(\sum\limits_{j'=0}^{\ell} A(k')_{j'} B(k')_{\ell-j'}
\mathcal{F}(k) \sum\limits_{n=0}^{\infty} (-1)^n F(k)_n
\big(h(k)^{\lg n \rg} \otimes d_{k}^{(n)}(\mathcal D_K (x^{\al+(\ell-j')\ep_{k'}+j'(\ep_0-\ep_{-k'})}))t^n
\big) t^\ell\bigg).
\end{split}
\end{equation*}
Set $\al(\ell,k',j'):=\al{+}(\ell{-}j')\ep_{k'}{+}j'(\ep_0{-}\ep_k)$. So, it is easy to see
\begin{equation*}
\begin{split}
&d_{k}^{(n)}(
\mathcal D_K (x^{\al+(\ell-j')\ep_{k'}+j'(\ep_0-\ep_k)}))=d_{k}^{(n)} (\mathcal D_K (x^{\al(\ell,k',j')})) \\
&\quad= \sum\limits_{j=0}^{n}
\frac{\prod\limits_{i=0}^{n-j-1}(||\al(\ell,k',j')||{-}\al(\ell,k',j')_0{+}i)}{(n{-}j)!}
\frac{\prod\limits_{i=0}^{j-1}(i{-}\al(\ell,k',j')_{-k})}{j!}
\mathcal D_K (x^{\al(\ell,k',j')+(n-j)\ep_k+j(\ep_0-\ep_{-k})})\\
&\quad=\sum\limits_{j=0}^{n} \frac{\prod\limits_{i=0}^{n-j-1}
(|\al|{+}(\ell{-}j'){+}\al_0{+}j'{-}2{-}\al_0{-}j'{+}i)}{(n{-}j)!} A(k)_j \mathcal D_K (x^{\al+(\ell-j')\ep_{k'}+j'(\ep_0-\ep_{-k'})+(n-j)\ep_k+j(\ep_0-\ep_{-k})})\\
&\quad=\sum\limits_{j=0}^{n} \frac{\prod\limits_{i=0}^{n-j-1}(||\al||{-}\al_0{+}(\ell{-}j'){+}i)}{(n{-}j)!}A(k)_j \mathcal D_K (x^{\al+(\ell-j')\ep_{k'}+j'(\ep_0-\ep_{-k'})+(n-j)\ep_k+j(\ep_0-\ep_{-k})})\\
&\quad=\sum\limits_{j=0}^{n}C_{n-j}^{\ell,j'} A(k)_j \mathcal D_K (x^{\al+(\ell-j')\ep_{k'}+j'(\ep_0-\ep_{-k'})+(n-j)\ep_k+j(\ep_0-\ep_{-k})}).
\end{split}
\end{equation*}
Thus we have
\begin{equation*}
\begin{split}
\mathcal{F}(k)&\Big(\sum\limits_{\ell=0}^{\infty}(-1)^\ell h(k')^{\lg
\ell \rg} \otimes (1{-}e(k')t)^{-\ell} d_{k'}^{(\ell)}(\mathcal D_K (x^{\al})
t^\ell )\Big)F(k)\\
&=\sum\limits_{\ell=0}^{\infty} (-1)^\ell h(k')^{\lg \ell \rg}
\otimes (1{-}e(k')t)^{-\ell} \cdot \sum\limits_{j'=0}^{\ell}
A(k')_{j'}B(k')_{\ell-j'} \sum\limits_{n=0}^{\infty}(-1)^n (1
\otimes (1{-}e(k)t)^{-n}) \\
&\quad\cdot\Big(h(k)^{\lg n \rg} \otimes \sum\limits_{j=0}^{n}
C_{n-j}^{\ell,j'} A(k)_j \mathcal D_K (x^{\al+(\ell-j')\ep_{k'}+j'(\ep_0-\ep_{-k'})+(n-j)\ep_k+j(\ep_0-\ep_{-k})})\\
&=\sum\limits_{n,\ell=0}^{\infty} \sum\limits_{j'=0}^{\ell} \sum\limits_{j=0}^{n} (-1)^{\ell+n} h(k')^{\lg \ell \rg}t^{n+\ell}
h(k)^{\lg n \rg} \otimes (1{-}e(k')t)^{-\ell} (1{-}e(k)t)^{-n}A(k)_j A(k')_{j'}\times\\
&\quad \times B(k')_{\ell-j'}
C_{n-j}^{\ell,j'}\mathcal D_K (x^{\al+(\ell-j')\ep_{k'}+j'(\ep_0-\ep_{-k'})+(n-j)\ep_k+j(\ep_0-\ep_{-k})})
t^{n+\ell}.
\end{split}
\end{equation*}
Thereby, we get the first formula. We can get the other formulas by using the similar argument as Lemma \ref{mHopf1}.

This completes the proof.
\end{proof}

\begin{lemm}\label{omHopf2}
Fix distinguished elements
 $h(k)=\mathcal D_K (x^{(\ep_k+\ep_{-k})}),~e(k)=\mathcal D_K (x^{(\ep_k+\ep_{0})})$, $h(k')=\mathcal D_K (x^{(\ep_{k'}+\ep_{-k'})}),~e(k')=\mathcal D_K (x^{(\ep_{k'}+\ep_{0})})$ with $1 \leq k
\neq k' \leq n$; the corresponding quantization of
 $U(\mathbf{K}(2n{+}1;\underline{1}))$ on $U_t(\mathbf{K}(2n{+}1;\underline{1}))$ $($also
on $U(\mathbf{K}(2n{+}1;\underline{1}))[[t]]$$)$ with the product
undeformed is given by
\begin{gather*}
\Delta(\mathcal D_K (x^{(\al)}))=\mathcal D_K (x^{(\al)}) \otimes
(1{-}e(k)t)^{\al_k-\al_{-k}}(1{-}e(k')t)^{\al_{k'}-\al_{-k'}}\hskip4cm\\
\quad+\sum\limits_{n,\ell=0}^{p-1}\sum\limits_{j'=0}^{\ell}
\sum\limits_{j=0}^{n}(-1)^{\ell+n}h(k')^{\lg
\ell \rg} h(k)^{\lg n\rg} \otimes (1{-}e(k')t)^{-\ell}(1{-}e(k)t)^{-n} \overline{A(k)}_j\overline{A(k')}_{j'}
\overline{B(k,k')}_{n,j}^{\ell,j'}\times\qquad \qquad \qquad \qquad  \\
\qquad \times \mathcal D_K \big(x^{(\al+(\ell-j')\ep_{k'}+j'(\ep_0-\ep_{-k'})+(n-j)\ep_k+j(\ep_0-\ep_{-k}))}\big) t^{n+\ell},\\
S(\mathcal D_K (x^{(\al)}))=-(1{-}e(k')t)^{\al_{-k'}-\al_{k'}}(1{-}e(k)t)^{\al_{-k}-\al_{k}} \sum\limits_{n,\ell=0}^{p-1}
\sum\limits_{j'=0}^{\ell}\sum\limits_{j=0}^{n}\overline{A(k')}_{j'}\overline{B(k')}_{\ell-j'}\overline{B(k,k')}_{n,j}^{\ell,j'} \times\qquad \qquad  \\
\quad \times
\mathcal D_K\big (x^{(\al+(\ell-j')\ep_{k'}+j'(\ep_0-\ep_{-k'})+(n-j)\ep_k+j(\ep_0-\ep_{-k}))}\big)
h(k)_1^{\lg n \rg} h(k')_{1}^{\lg \ell \rg} t^{\ell+n},\\
\varepsilon(\mathcal D_K (x^{(\al)}))=0,\hskip10cm
\end{gather*}
where $0 \leq \al \leq \tau$, if $2n+1 \nequiv
0\; (\mathrm{mod}\, p);$ and $0 \leq \al < \tau$, if $2n+1 \equiv 0\;(\mathrm{mod}\,
p)$. $\overline{A(k')}_{j'}=(-1)^{j'} \binom{\al_{k'}{+}\ell{-}j'}{\al_{k'}}$, for $j' \leq \al_{-k'}$, $\overline{A}(k)_j=(-1)^j \binom{\al_k{+}n{-}j}{\al_k}$, for $j \leq \al_{-k}$. Otherwise, $\overline{A(k')}_{j'}=\overline{A(k)}_{j}=0$, $\overline{B(k,k')}_{n,j}^{\ell,j'}= \binom{\al_0{+}j{+}j'}{j{+}j'}\binom{j{+}j'}{j} \prod\limits_{i=0}^{\ell+n-j'-j-1}(||\al||{-}\al_0{+}i)$.
\end{lemm}
\begin{proof} Let us consider
\begin{equation*}
\begin{split}
&\Delta(\mathcal D_K (x^{(\al)}))=\frac{1}{\al
!}\Delta(\mathcal D_K (x^\al))=\mathcal D_K (x^{(\al)}) \otimes
(1{-}e(k)t)^{\al_k-\al_{-k}}(1{-}e(k')t)^{\al_{k'}-\al_{-k'}}\\
&\quad+ \sum\limits_{n,\ell=0}^{p-1} h(k')^{\lg \ell \rg} h(k)^{\lg
n \rg} \otimes
(1{-}e(k')t)^{-\ell}(1{-}e(k)t)^{-n}\sum\limits_{j'=0}^{n}\sum\limits_{j=0}^{n}
A(k')_{j'}B(k')_{\ell-j'} C_{n-j}^{\ell,j'}A(k)_j \times \\
&\qquad \times
\frac{(\al{+}(\ell{-}j')\ep_{k'}{+}j'(\ep_0{-}\ep_{-k'}){+}(n{-}j)\ep_k{+}j(\ep_0{-}\ep_{-k}))!}{\al
!}\times \\
&\qquad \times
\mathcal D_K (x^{(\al+(\ell-j')\ep_{k'}+j'(\ep_0-\ep_{-k'})+(n-j)\ep_k+j(\ep_0-\ep_{-k}))}).
\end{split}
\end{equation*}
And we can see that
\begin{equation*}
\begin{split}
A(k')_{j'}&B(k')_{\ell-j'} C_{n-j}^{\ell,j'}A(k)_j
\frac{(\al{+}(\ell{-}j')\ep_{k'}{+}j'(\ep_0{-}\ep_{-k'}){+}(n{-}j)\ep_k{+}j(\ep_0{-}\ep_{-k}))!}{\al
!}\\
&=\frac{\prod\limits_{i=0}^{j'-1}(i{-}\al_{-k'})}{j'!}\cdot\frac{\prod\limits_{i=0}^{\ell-j'-1}(||\al||{-}\al_0{+}i)}{(\ell{-}j')!}\cdot
\frac{\prod\limits_{i=0}^{j-1}(i{-}\al_{-k})}{j!}\cdot\frac{\prod\limits_{i=0}^{n-j-1}(||\al||{-}\al_0{+}(\ell{-}j'){+}i)}{(n{-}j)!}\cdot\\
&\quad \cdot\frac{(\al_{k'}{+}(\ell{-}j'))!}{\al_{k'}!}\cdot\frac{(\al_k{+}(n{-}j))!}{\al_k!}\cdot\frac{(\al_{-k'}{-}j')!}{\al_{-k'}!}\cdot\frac{(\al_{-k}{-}j)!}{\al_{-k}!}\cdot
\frac{(\al_0{+}j{+}j')!}{\al_0!}\\
&=(-1)^{j'}\frac{\al_{-k'}(\al_{-k'}{-}1) \cdots (\al_{-k'}{-}j')!}{\al_{-k'}!}\cdot \frac{\big(\al_{k'}{+}\ell{-}j'\big)!}{\al_{k'}!(\ell{-}j')!}\cdot\\
&\quad \cdot(-1)^j \frac{\al_{-k}(\al_{-k}{-}1)\cdots (\al_{-k}{-}j)!}{\al_{-k}!}\cdot\frac{(\al_k{+}n{-}j)!}{\al_k! (n{-}j)!}\prod\limits_{i=0}^{\ell+n-j'-j-1}(||\al||{-}\al_0{+}i) \frac{(\al_0{+}j{+}j')!}{\al_0! j!j'!}\\
&=\overline{A(k')}_{j'} \overline{A(k)}_j \binom{\al_0{+}j{+}j'}{j{+}j'}\binom{j{+}j'}{j} \prod\limits_{i=0}^{\ell+n-j'-j-1}(||\al||{-}\al_0{+}i) \\
&=\overline{A(k')}_{j'} \overline{A(k)}_j\overline{B(k,k')}_{n,j}^{\ell,j'}.
\end{split}
\end{equation*}

We can get the other formulas by using a similar argument.
\end{proof}

\begin{lemm}
For $s \geq 1$, one has
\begin{equation*}\label{omrelation1}
\begin{split}
&\Delta(\mathcal D_K (x^{(\al)})^s)=\sum_{0\le j'\le s\atop n, \ell\ge
0}\binom{s}{j'}(-1)^{n+\ell}\mathcal D_K (x^{(\al)})^{j'} h(k')^{\lg \ell \rg}
h(k)^{\lg n \rg} \otimes(1{-}e(k)t)^{j'(\al_k-\al_{-k})-n}\cdot\\
&\hskip3.5cm\cdot (1{-}e(k')t)^{j'(\al_{k'}-\al_{-k'})-\ell}
d_k^{(n)}\big(d_{k'}^{(\ell)}(\mathcal D_K  (x^{(\al)}))^{s-j'}\big) t^{n+\ell},\\
&S((\mathcal D_K (x^{(\al)})^s))=(-1)^s
(1{-}e(k')t)^{-s(\al_{k'}-\al_{-k'})}(1{-}e(k)t)^{-s(\al_k-\al_{-k})}\times\\
&\hskip3.5cm  \times \sum\limits_{n,\ell=0}^{p-1}
d_k^{(n)}\big(d_{k'}^{(\ell)}(\mathcal D_K (x^{(\al)}))^{s}\big) h(k')_1^{\lg
\ell \rg} h(k)_1^{\lg n \rg} t^{n+\ell}.
\end{split}
\end{equation*}
\end{lemm}
\begin{proof}
By the same proof as that of Lemma \ref{mrelation1}.
\end{proof}

To describe $\mathbf{u}_{t,q}(\mathbf{K}(2n{+}1;\underline{1}))$ explicitly, we
still need an auxiliary Lemma.

\begin{lemm}\label{omrelation2}
Set $e(k)=\mathcal D_K (x^{(\ep_k+\ep_{0})}),~e(k')=\mathcal D_K (x^{(\ep_{k'}+\ep_{0})})$, $d_k^{(n)}=\frac{1}{n!}(\ad
 e(k))^n$, $d_{k'}^{(\ell)}=\frac{1}{\ell!} (\ad e(k'))^\ell$, Then
 \begin{gather*}
d_k^{(n)}
d_{k'}^{(\ell)}(\mathcal D_K (x^{(\al)})
=\sum\limits_{j'=0}^{\ell}\sum\limits_{j=0}^{n}\overline{A(k')}_j
~\overline{A(k)}_j \overline{B(k,k')}_{n,j}^{\ell,j'}
 \mathcal D_K (x^{\al+(\ell-j')\ep_{k'}+j'(\ep_0-\ep_{-k})+(n-j)\ep_k+j(\ep_0-\ep_{-k'})}),\tag{\text{\rm i}}\\
 d_k^{(n)} d_{k'}^{(\ell)}(\mathcal D_K (x^{(\ep_i+\ep_{-i})})=\delta_{\ell
0}\delta_{n 0}\mathcal D_K (x^{(\ep_i+\ep_{-i})})- \delta_{\ell 1}\delta_{n
0}\delta_{ik'}e(k') -\delta_{\ell 0}
\delta_{n 1} \delta_{i k} e(k),\tag{\text{\rm ii}} \ \\
 d_k^{(n)} d_{k'}^{(\ell)}(\mathcal D_K (x^{(\ep_0)}))=\delta_{\ell
0}\delta_{n 0}\mathcal D_K (x^{(\ep_0)})- \delta_{\ell 1}\delta_{n 0}e(k')
-\delta_{\ell 0}
\delta_{n 1}  e(k),\qquad\quad \\
  d_k^{(n)} d_{k'}^{(\ell)}(\mathcal D_K (x^{(\al)}))^p=\delta_{\ell
 0}\delta_{n0} (\mathcal D_K (x^{(\al)}))^p-\delta_{\ell
 0}\delta_{n1}(\delta_{\al,\ep_k+\ep_{-k}}{+}\delta_{\al,\ep_0})e(k)\quad\quad \tag{\text{\rm iii}}\\
-\delta_{\ell 1}\delta_{n
0}(\delta_{\al,\ep_{k'}{+}\ep_{-k'}}{+}\delta_{\al,\ep_0})e(k').
\end{gather*}
\end{lemm}

\begin{theorem}\label{omHopf}
Fix distinguished elements $h(k):=\mathcal D_K (x^{(\ep_k+\ep_{-k})})$, $e(k)=\mathcal D_K  (x^{(\ep_k+\ep_{0})})$, $h(k')=\mathcal D_K (x^{(\ep_{k'}+\ep_{-k'})})$, $e(k')=\mathcal D_K (x^{(\ep_{k'}+\ep_{0})})$, with
$1 \leq k \neq k' \leq n$; there is a noncommutative and noncocommutative Hopf algebra $({\bf u}_{t,q}({\bf K})(2n{+}1;\underline{1}),m,\iota,\Delta,S,\varepsilon)$~over ~$\mathcal{K}[t]_{p}^{(q)}$ with the product undeformed, whose coalgebra structure is given by
\begin{equation*}
\begin{split}
&\Delta(\mathcal D_K (x^{(\al)}))=\mathcal D_K (x^{(\al)}) \otimes (1{-}e(k)t)^{\al_k-\al_{-k}}(1{-}e(k')t)^{\al_{k'}-\al_{-k'}}\\
&\qquad +\sum\limits_{n,\ell=0}^{p-1} (-1)^{\ell+n}h(k')^{\lg \ell \rg} h(k)^{\lg n \rg} \otimes (1{-}e(k')t)^{-\ell}(1{-}e(k)t)^{-n} d_k^{(n)}d_{k'}^{(\ell)}(\mathcal D_K (x^{(\al)}))t^{n+\ell},\\
&S(x^{(\al)})=-(1{-}e(k')t)^{\al_{-k'}-\al_{k'}}(1{-}e(k)t)^{\al_{-k}-\al_k}\sum\limits_{n,\ell=0}^{p-1} d_k^{(n)}d_{k'}^{(\ell)} (\mathcal D_K (x^{(\al)}))h(k)_1^{\lg n \rg} h(k')_1^{\lg \ell \rg} t^{n+\ell},\\
&\varepsilon(\mathcal D_K (x^{(\al)}))=0,
\end{split}
\end{equation*}
where $0 \leq \al \leq \tau$, $\dim_{\mathcal{K}}\mathbf{u}_{t,q}(\mathbf{K}(2n{+}1;\underline{1}))=p^{p^{2n+1}+1}$ for $2n+1 \nequiv 0 \;(\mathrm{mod}\, p)$, and $0 \leq \al < \tau$, $\dim_{\mathcal{K}}\mathbf{u}_{t,q}(\mathbf{K}(2n{+}1;\underline{1}))=p^{p^{2n+1}}$ for $2n+1 \equiv 0\; (\mathrm{mod}\, p)$.
\end{theorem}

\begin{remark}
As for the the Drinfel'd twist associated with (ix), we have similar formulas. So we omit the statements here.
\end{remark}

\medskip
\subsection{Open Questions} As is well-known, the authors \cite{AS} gave a certain classification
for the finite-dimensional complex pointed Hopf algebras with abelian finite group algebras (whose orders satisfying some conditions) as the coradicals.
This is the most important achievement in Hopf algebra theory during the last decade.
However, the new Hopf algebras of prime-power dimensions we obtained above should be pointed over a field of positive characteristic.

Before concluding the paper, we prefer to propose the following
interesting questions for further considerations.

{\bf Question 1.} \ Assume $\mathcal K$ is an algebraically
closed field with $t, q\in \mathcal K$. How many non-isomorphic (pointed) Hopf
algebra structures can be equipped on the universal restricted
enveloping algebra $\mathbf u(K(2n{+}1;\underline{1}))$? How to classify the pointed Hopf algebras of the given prime-power dimension over a field of positive characteristic?

\medskip
{\bf Question 2.} \ What are the conditions for $\mathbf
u_{t,q}(K(2n{+}1;\underline{1}))$ to be
a ribbon Hopf algebra (see \cite{HW2} and references therein)?

\medskip
{\bf Question 3.} \ It might be interesting to consider the
tensor product structures of representations for $\mathbf
u_{t,q}(K(2n{+}1;\underline{1}))$. How
do their tensor categories behave?

\bibliographystyle{amsalpha}

\end{document}